\documentclass[12pt,a4paper, twoside]{article}

\usepackage{amsmath,amssymb,amsthm,amsfonts,mathrsfs,amscd,environ}
\usepackage{dsfont}
\usepackage{latexsym,enumerate,geometry,extarrows}
\usepackage{stmaryrd}

\usepackage{amsfonts}
\usepackage{xcolor}

\usepackage{verbatim,fancyhdr,enumitem}
 \newtheorem{thm}{Theorem}[section]
 \newtheorem{lem}[thm]{Lemma}
 
 \newtheorem{exa}[thm]{Example}
 
 \newtheorem{defn}{Definition}[section]
 \newtheorem{rem}{Remark}[section]
 
 \numberwithin{equation}{section}

\def\dif{{\mathord{{\rm d}}}}
\def\no{\nonumber}

\def\mR{{\mathbb R}}
\def\mE{{\mathbb E}}
\def\mM{{\mathbb M}}

\def\mE{{\mathbb E}}
\def\mF{{\mathbb F}}

\def\mI{{\mathbb I}}
\def\mJ{{\mathbb J}}

\def\mM{{\mathbb M}}

\def\mR{{\mathbb R}}
\def\mS{{\mathbb S}}

\def\mU{{\mathbb U}}

\def\sF{{\mathscr F}}

\def\sB{{\mathscr B}}

\def\sF{{\mathscr F}}

\def\sH{{\mathscr H}}

\def\sM{{\mathscr M}}

\def\sP{{\mathscr P}}

\def\sU{{\mathscr U}}

\def\bd{\begin{defn}}
\def\ed{\end{defn}}
\def\bl{\begin{lem}}
\def\el{\end{lem}}
\def\bt{\begin{thm}}
\def\et{\end{thm}}
\def\br{\begin{rem}}
\def\er{\end{rem}}

\allowdisplaybreaks

\title{{\bf Extended MF-AFBSDDEswCN   and stochastic optimal controls  of  linear system}
}

\author{
{\bf Hao Wu$^{a)}$ }\\
\footnotesize{$^{a)}$School of Mathematics and Statistics,
South-Central Minzu University}\\
\footnotesize{ Wuhan, Hubei 430000, P.R.China}\\
\footnotesize{Email: wuhaomoonsky@163.com},
}

\begin{document}

\maketitle

\begin{abstract}
The research  studies  the  extended   mean-field  anticipated forward-backward   stochastic  delayed differential equations with common noise (extended MF-AFBSDDEswCN) by expanding the domination-monotonicity conditions. We generalize the conditions from a linear setting to a nonlinear one by including nonlinear adjoint functions which are key to guaranteeing the well-posedness of extended MF-AFBSDDEswCN.  Utilizing this broader well-posedness framework for extended  MF-AFBSDDEswCN, combined with other refined analytical tools, we examine two classes of stochastic  optimal control problems. These include a linear-convex problem and a linear-quadratic problem with input constraints that can be both time-dependent and random.  Moreover,  the control elements regarding the initial values are path-dependent, which  will cause some essential difficulties.   For each problem, we establish the existence and uniqueness of optimal controls and provide their explicit, closed-form representations.
\end{abstract}\noindent

\noindent

Optimal controls;
  Anticipated; Delayed;  Linear-convex problem.

\section{Introduction}
Throughout this paper, the superscript $"\top"$ represents the transpose of a vector or a matrix.  Brownian motion $\{B^{1}(\cdot):=(B^{1}_{1}(\cdot),B^{1}_{2}(\cdot), \cdots, B^{1}_{d}(\cdot) )^{\top}, 0\leq t \leq T\}$ supported on $(\Omega^1, \sF^1,  \mathbb{P}^1)$  is referred to as a $d$-dimensional idiosyncratic noise with filtration $\mF^{1}:= \{\sF^{1}_{t}\}_{t\geq 0},$  where $\mF^{1}$ generated by Brownian motion $(B_t^1)_{t\geq 0}$, i.e., $ \sF^{1}_{t}=\sigma\{B_{s}^{1}, s\leq t \}.  $   Brownian motion $\{B^{0}(\cdot):=(B^{0}_{1}(\cdot),B^{0}_{2}(\cdot), \cdots, B^{0}_{d}(\cdot) )^{\top}, 0\leq t \leq T\}$, carried on $(\Omega^0, \sF^0,  \mathbb{P}^0)$ with filtration $\mF^{0}:= \{\sF^{0}_{t}\}_{t\geq 0}$  is referred to a $d$-dimensional common noise,  where $\mF^{0}$ generated by Brownian motion $(B_t^0)_{t\geq 0}$, i.e.,   $ \sF^0_{t}=\sigma\{B_{s}^{0}, s\leq t \} .$     We define the product space $(\Omega, \sF, \mathbb{P})$, where $\Omega = \Omega^0 \times \Omega^1$. The sigma-algebra $\mathscr{F}$ is the completion of $\mathscr{F}^0 \otimes \mathscr{F}^1$, and $\mathbb{P}$ is the completion of $\mathbb{P}^0 \otimes \mathbb{P}^1$. The filtration $(\mathscr{F}_t)_{t \ge 0}$ is the complete and right-continuous augmentation of $(\mathscr{F}^0_t \otimes \mathscr{F}^1_t)_{t\geq 0}$.   The expectations with respect to $\mathbb{P},$ $\mathbb{P}^{1},$ and $\mathbb{P}^0$ are denoted by $\mathbb{E},$  $\mathbb{E}^1$ and $\mathbb{E}^0$ respectively.
For a random variable $\varsigma$,  its probability distribution is denoted by $\mathcal{L}(\varsigma)$. For $k \ge 1$, let $L^{k}_{\sF}(\Omega, \mathbb{R}^n)$ be the set of random variables $\varsigma: \Omega \to \mathbb{R}^n$ satisfying  $\mathbb{E}[|\varsigma|^k] < \infty$. Let $L^{k}_{\sF}(\Omega^0, \mathcal{P}_k(\mathbb{R}^n))$ be the set of random variables $u: \Omega^0 \to \mathcal{P}_k(\mathbb{R}^n)$ (where $\mathcal{P}_k(\mathbb{R}^n)$ denotes the set of probability measures on $\mathbb{R}^n$ whose $k$-th moment exists), satisfying $\mathbb{E}^0[(W^k(u, \delta_0))] < \infty$, where $\delta_0$ denotes the Dirac measure at $0 \in \mathbb{R}^n$. For a random variable $\varsigma: \Omega \to \mathbb{R}^n$, we define the mapping $\mathcal{L}^1(\varsigma)(\omega^0) = \mathcal{L}(\varsigma(\omega^0, \cdot))$ for any $\omega^0 \in \Omega^0$.
 According to \cite{CarmonaDelarue20182}, the mapping $\mathcal{L}^1(\varsigma)(\cdot)$ is not only a random variable from $\Omega^0$ to $\mathcal{P}(\mathbb{R}^d)$ but also represents the conditional distribution of $\xi$ given $\mathscr{F}^0$. With a slight abuse of notation, we do not distinguish a random variable $\varsigma$ on $\Omega^0$ from its natural extension on $\Omega$.

Stochastic delayed differential equations, often abbreviated as SDDEs, are a class of mathematical models that describe systems influenced by both randomness and time lags. In many real-world processes, the future state of a system does not depend solely on its current state, but also on its past states. This delay, or time lag, can arise from finite signal transmission speeds, memory effects, or biological and physical constraints. At the same time, these systems are frequently subject to random fluctuations, or noise, which can come from environmental disturbances, measurement errors, or inherent uncertainties. SDDEs combine these two features: they incorporate a deterministic or stochastic forcing term that depends on the history of the system, while also including a random component to capture unpredictable variations. Such equations appear in a wide range of fields, including population dynamics, where the reproduction of a species may depend on the population size at an earlier time; financial markets, where asset prices respond to past information and are driven by random shocks; and control theory, where delayed feedback and noisy measurements are common. The study of SDDEs is challenging because the presence of both delay and randomness introduces complex dynamics, such as oscillations, stability shifts, and non-Markovian behavior. Understanding these equations is essential for predicting the behavior of real-world systems that are inherently noisy and possess memory.

The duality between anticipated backward stochastic differential equations (ABSDEs) and SDDEs is rooted in the structural symmetry of their dependence on the state trajectory: ABSDEs incorporate future information through terms like $E[X(t+\theta)|\sF_{t}]$
while  SDDEs involve past values such as $X(t-\theta).$ This duality is particularly powerful in stochastic control, as it allows one to reformulate a control problem with delays (e.g., a system with a delayed state feedback) into an equivalent problem with anticipated dynamics, which often admits a more tractable solution via forward-backward stochastic differential equations (FBSDEs). For instance, in optimal control of delay systems, the adjoint process satisfies an anticipated backward stochastic differential equation, enabling the design of predictors or compensators for the delay. Conversely, in problems with anticipatory information (e.g., model predictive control with future reference signals), the dual delay representation facilitates the synthesis of robust controllers using classical delay-robust stability techniques. This interplay bridges the analysis of systems with memory and prediction, and has found applications in finance, engineering, and networked control systems. We refer the reader to \cite{ PengYang,  LvTaoWu, WangS} and references therein.

Stochastic differential equations of the McKean-Vlasov type, abbreviated as MVSDEs, constitute a category of SDEs where the coefficients depend not only on the state process but also on the probability distribution of the solution process \cite{CarmonaDelarue2018}. In this current study, many scholars  shift to McKean-Vlasov SDEs that incorporate common noise, which are also referred to in the existing literature as conditional McKean-Vlasov SDEs \cite{LackerShkolnikovzhang1990}. Regarding these types of SDEs, the system receives noise input from two distinct sources: one is an individualistic source of randomness, and the other is what is termed as common noise.

Comprehending large populations and intricate systems presents a key challenge in mathematical modeling, with applications spanning financial, economic, social, and electrical networks, as well as communication infrastructure, biological and ecological groups, and extensively engineered networks, among others. These systems generally develop from various sources of uncertainty—both individual (idiosyncratic noise) and shared (common noise) and the actions of each entity frequently rely on the collective state of the population. This inherently leads to McKean-Vlasov dynamics and the formalisms of mean-field games (MFG) and mean-field control (MFC), whose mathematical analysis is complex due to the inherent interdependence via the state's probability law. Specifically, when a common noise exists, the mean field ceases to be fixed and instead transforms into a probabilistic conditional distribution, giving rise to considerable further mathematical difficulties. For a comprehensive overview, we direct readers to \cite{LasryLions2007,Bensoussan2013,CarmonaDelarue2018,CarmonaDelarue20182}. Additionally, we point to \cite{CossoPerelli,TalbiTouzizhang2023,TalbiTouzizhang20231,TalbiTouzizhang20232} for the closely associated field of mean-field optimal stopping.

A significant portion of current research on mean field models concentrates on Brownian common noise. For problems involving Brownian common noise, notable contributions include the dynamic programming principle established by Pham and Wei \cite{PhamWei2017}
  for closed-loop controls, and by Djete et al. \cite{DjetePossamaoTan2022}
  in the context of non-Markovian frameworks and open-loop controls. Furthermore, Liang et al. \cite{DenkertKharroubiPham2024}
  investigated mean-field control problems with common noise. Lacker and Flem \cite{LasryFlem2023} present a thorough framework for closed-loop convergence in mean field games with common noise, introducing the notion of weak mean field equilibrium. Huang et al. \cite{HuangWang2023} utilized an innovative decoupling method to examine the control challenges of linear-quadratic (LQ) stochastic backward large population systems that involve partial information and common noise. Hua and Luo \cite{HuaLuo2024} investigate a type of LQ expanded mean field games with common noises, wherein the interactions between mean fields are solely mediated by the (conditional) expectations of states and controls, accommodating non-linear shifts in state coefficients and cost functions. Si and Shi \cite{SiShi}  applied a fixed-point methodology to scrutinize a range of general linear-quadratic mean-field Stackelberg games featuring common noise.

We note that a highly favored category within these models is the linear quadratic (LQ) type. This type is recognized for its ease of analysis and wide array of applications (cf.\cite{yong2013,HuangLiyong2015,LiSunxiong2019,JingruiSun,Pham2016}
),  prompting the investigation of appropriate Riccati equations. These typically involve Riccati ordinary differential equations when only idiosyncratic noises are present, and backward Stochastic Riccati equations in the scenario of common noise.
Furthermore,  investigating LQ problems have more  motivations. Theoretically, the LQ framework provides a manageable yet extensive structure that frequently allows for clear or easily obtainable solutions. It acts as an important reference point and a fundamental element for comprehending more complex nonlinear issues. By analyzing mean-field FBSDEs in the context of LQ problems, one can establish precise criteria for feasibility, acquire optimal controls in a straightforward format, and develop a deep understanding of how individual optimization, randomness, and mean-field interactions are connected. In practice, LQ mean-field control models can be directly applied across a wide range of fields. Examples include portfolio optimization in finance with price impacts, consensus control in multi-agent robotic systems, and macroeconomic models where agents' decisions depend on aggregate economic indicators(cf. \cite{Bardi2012,LiLiyu2020,LinJiangzhang2019}). Therefore, the study of linear-quadratic optimal control for mean-field FBSDEs not only advances the mathematical theory of stochastic control and mean-field games but also provides a versatile toolkit for designing and analyzing optimal strategies in complex, interconnected systems.

In this article,  we aim to examine the subsequent extended MF-FBSDEswCN, subject to nonlinear domination-monotonicity conditions, by establishing\\ \
\begin{align*}
\mathfrak{U}_{i}(t):&=\big(\big. (\mE[x_{1}(t)])^{\top}, (\mE[y_{1}(t)])^{\top}, (\mE[x_{2}(t)])^{\top}, (\mE[y_{2}(t)])^{\top}, \\
 &\quad\quad\quad\quad\quad(\mE^{\sF^{0}_{t}}[x_{1}(t)])^{\top}, (\mE^{\sF^{0}_{t}}[y_{1}(t)])^{\top}, (\mE^{\sF^{0}_{t}}[x_{2}(t)])^{\top}, (\mE^{\sF^{0}_{t}}[y_{2}(t)])^{\top},\\
&\quad\quad \quad(x_{1}(t-\theta))^{\top}, (\mE^{\sF_{t}}[y_{1}(t+\theta)])^{\top}, (x_{2}(t-\theta))^{\top}, (\mE^{\sF_{t}}[y_{2}(t+\theta)])^{\top}, \\
  &\quad\quad \quad\quad
  x_{i}(t)^{\top}, y_{i}(t)^{\top}, z^{0}_{i}(t)^{\top}, z^{1}_{i}(t)^{\top}\big)\big.^{\top}, i=1,2 :
\end{align*}

\begin{align}\label{1}
\begin{cases}
&\dif x_{1}(t)=b_{1}(t,\mathfrak{U}_{1}(t))\dif t
+\sum^{d}_{j=1}\sigma^{0}_{1j}(t,\mathfrak{U}_{1}(t))\dif B^{0}_{j}(t)+\sum^{d}_{j=1}\sigma^{1}_{1j}(t,\mathfrak{U}_{1}(t))\dif B^{1}_{j}(t), t\in [0,T],\\
&\dif y_{1}(t) =f_{1}(t,\mathfrak{U}_{1}(t))\dif t+\sum^{d}_{j=1} z^{0}_{1j}(t)\dif B^{0}_{i}(t)+\sum^{d}_{j=1} z^{1}_{1j}(t)\dif B^{1}_{j}(t), t\in [0,T],\\
&  y_{1}(T)=\Phi_{1}(x_{1}(T),x_{2}(T) ), y_{1}(t)=\eta_{1}(t), t\in (T, T+\theta], \\
& x_{1}(0)=\Psi_{1}(y_{1}(0),y_{2}(0)), x_{1}(t)=\ddot{\Psi}_{1}(\mE[y_{1}(t+\theta)],\mE[y_{2}(t+\theta)]), t\in [-\theta, 0).
\end{cases}
\end{align}

\begin{align}\label{2}
\begin{cases}
&\dif x_{2}(t)=b_{2}(t,\mathfrak{U}_{2}(t))\dif t
+\sum^{d}_{j=1}\sigma^{0}_{2j}(t,\mathfrak{U}_{2}(t))\dif B^{0}_{j}(t)+\sum^{d}_{j=1}\sigma^{1}_{2j}(t,\mathfrak{U}_{2}(t))\dif B^{1}_{j}(t), t\in [0,T],\\
&\dif y_{2}(t) =f_{2}(t,\mathfrak{U}_{2}(t))\dif t+\sum^{d}_{j=1} z^{0}_{2j}(t)\dif B^{0}_{i}(t)+\sum^{d}_{j=1} z^{1}_{2j}(t)\dif B^{1}_{j}(t), t\in [0,T],\\
&  y_{2}(T)=\Phi_{2}(x_{1}(T),x_{2}(T) ), y_{2}(t)=\eta_{2}(t), t\in (T, T+\theta], \\
& x_{2}(0)=\Psi_{2}(y_{1}(0),y_{2}(0)), x_{2}(t)=\ddot{\Psi}_{2}(\mE[y_{1}(t+\theta)],\mE[y_{2}(t+\theta)]), t\in [-\theta, 0).
\end{cases}
\end{align}
Denote $W_{i}(\cdot):=(x_{i}(\cdot)^{\top},$ $y_{i}(\cdot)^{\top},z^{0}_{i}(\cdot)^{\top},z^{1}_{i}(\cdot)^{\top})^{\top}, i=1,2,$   which represent the unknown processes  with    $z^{k}_{i}(\cdot):=(z^{k}_{i1}(\cdot)^{\top},z^{k}_{i2}(\cdot)^{\top},\cdots, ,z^{k}_{id}(\cdot)^{\top})^{\top}, \sigma^{k}_{i}(\cdot):=(\sigma^{k}_{i1}(\cdot)^{\top},\sigma^{k}_{i2}(\cdot)^{\top},\cdots,\sigma^{k}_{id}(\cdot)^{\top})^{\top},$  $i=1,2, k=0,1.$   The coefficients are denoted as $\Phi_{i}: \Omega\times\mR^{n}\times \mR^{n}\rightarrow \mR^{n}, $ $\Psi_{i}, \ddot{\Psi}_{i}: \mR^{n}\times \mR^{n}
  \rightarrow \mR^{n} $  and $f_{i},b_{i}, \sigma^{k}_{ij}: \Omega \times [0,T] \times \mR^{n+n+n+n+n+n+n+n+n+n+n+n+n+n+nd+nd}\rightarrow \mR^{n}, j=1,2,\cdots d, i=1,2, k=0,1.$ For simplicity, denote $\aleph_{i}(\cdot):=(f_{i}(\cdot)^{\top},b_{i}(\cdot)^{\top},\sigma^{0}_{i}(\cdot)^{\top},\sigma^{1}_{i}(\cdot)^{\top})^{\top}, i=1, 2. $   Besides,  the combination of Eq.\eqref{1} and Eq.\eqref{2} is denoted by $(\pi).$  Consequently, all coefficients of the System $(\pi)$ is encompassed by $(\Phi_{i}, \Psi_{i}, \ddot{\Psi}_{i},\aleph_{i}), i=1,2$. It is evident that there are two initial couplings within System $(\pi).$    Subsequently, when we examine the optimal control problems associated with System $(\pi),$ the two initial values will serve as the two initial controls. This article aims to establish results concerning the well-posedness of System $(\pi),$ specifically addressing existence, uniqueness, and  the applications in stochastic optimal control problems.

Compared with the research in    Tian and Yu \cite{Tianyu2023},   Liu et al. \cite{LiuNiuWangyu2026} and  Buckdahn et al. \cite{BuckdahnLiPeng2009},  we differ from their systems \cite{Tianyu2023,LiuNiuWangyu2026,BuckdahnLiPeng2009}  in several ways.    Firstly,  one of the main motivations for studying the  mean-field equation of the form of System $(\pi)$ comes from  Buckdahn et al. \cite{BuckdahnLiPeng2009}.  Under the certain  conditions,  they got the existence and uniqueness result of solution. Moreover, Buckdahn et al. \cite{BuckdahnLiPeng2009} prove that the unique solution for the  system they considered  is closely related to  a nonlocal partial differential equation.   But they  did not take into account the  couplings with respect to  the initial values and the terminal values, delayed terms, anticipated terms and common noise. Based on Buckdahn et al. \cite{BuckdahnLiPeng2009},    we address the inspired mean-field system containing  two couplings with respect to  the initial values and the terminal values of  the System $(\pi)$, delayed terms, anticipated terms and common noise. For subsequent control applications,  these two couplings will leads to  the initial values  become control elements.     Secondly, compared with  the system in Tian and Yu \cite{Tianyu2023},  we examine a generalized mean field system, not only  including a complicated  mean field  coupling structure, but also  including  common noise, delayed terms, anticipated terms.  Moreover,  we consider the system with nonlinear
domination-monotonicity conditions while Tian and Yu \cite{Tianyu2023} investigate their system under the linear
domination-monotonicity conditions. Due to this nonlinear domination-monotonicity conditions and the common noise, delayed terms, anticipated terms appeared in the  system, the system has different application scenarios. Thirdly,  although  Liu et al.   \cite{LiuNiuWangyu2026}  have considered their system with nonlinear domination-monotonicity conditions on coefficients, the system they considered is not a mean-field system   and  adjoint functions  didn't contain  common noise, delayed terms, anticipated terms  while  our system does.  Another different point is the most crucial one. When we study the corresponding control problem, our initial control element is path-dependent, rather than a point as  in  \cite{LiuNiuWangyu2026}. This difference makes the value function essentially different and more challenging to deal with.    Specifically, the controlled system is formulated as follows using linear coupled SDEs defined on $[0, T],$
\begin{align}\label{4}
&\dif x_{1}(t)=[F_{1}(t)x_{1}(t) +H_{1}(t)u_{1}(t)+\tau\bar{H}_{1}(t)\mE[u_{1}(t)]+\tau \ddot{H}_{1}(t)\mE^{\sF^{0}_{t}}[u_{1}(t)]+\tau\tilde{H}_{1}(t)u_{1}(t-\theta)\no\\
&\quad + \bar{F}_{2}(t)\mE[x_{2}(t)]+ \bar{F}_{1}(t)\mE[x_{1}(t)]+ \tilde{F}_{2}(t)x_{2}(t-\theta)+\tilde{F}_{1}(t) x_{1}(t-\theta)\no\\
&\quad+ \ddot{F}_{2}(t)\mE^{\sF^{0}_{t}}[x_{2}(t)] +\ddot{F}_{1}(t)\mE^{\sF^{0}_{t}}[ x_{1}(t)]+\upsilon_{1}(t)]\dif t\no\\
&\quad+\sum^{d}_{j=1}[K_{1j}(t)x_{1}(t) +D_{1j}(t)u_{1}(t)+\iota_{1j}(t)]\dif B^{0}_{j}(t), \no\\
&\quad+\sum^{d}_{j=1}[\bar{K}_{1j}(t)x_{1}(t) +\bar{D}_{1j}(t)u_{1}(t)+\bar{\iota}_{1j}(t)]\dif B^{1}_{j}(t), t\in [0,T],\no\\
\end{align}
\begin{align}\label{4d}
&\dif x_{2}(t)=[F_{2}(t)x_{2}(t)+H_{2}(t)u_{2}(t)+\tau\bar{H}_{2}(t)\mE[u_{2}(t)]+\tau\tilde{H}_{2}(t)u_{2}(t-\theta)\no\\
&\quad+\tau\ddot{H}_{2}(t)\mE^{\sF^{0}_{t}}[u_{2}(t)]+ \bar{F}_{2}(t)\mE[x_{1}(t)]+\tilde{F}_{2}(t)x_{1}(t-\theta) +\ddot{F}_{2}(t)\mE^{\sF^{0}_{t}}[ x_{1}(t)]+\upsilon_{2}(t) ]\dif t,\no\\
& \quad+\sum^{d}_{j=1}[K_{2j}(t)x_{2}(t) +D_{2j}(t)u_{2}(t)+\iota_{2j}(t)]\dif B^{0}_{j}(t), t\in [0,T],\no\\
&\quad +\sum^{d}_{j=1}[\bar{K}_{2j}(t)x_{2}(t) +\bar{D}_{2j}(t)u_{2}(t)+\bar{\iota}_{2j}(t)]\dif B^{1}_{j}(t), t\in [0,T],
\end{align}
with the initial values and terminal values
\begin{align*}
\begin{cases}
x_{1}(0)=\sH\xi_{1}+x_{0}, x_{2}(0)=\sH\xi_{2}+x_{0},\\
 x_{1}(t)=\kappa_{1}(t), x_{2}(t)=\kappa_{2}(t), t\in [-\theta, 0).
\end{cases}
\end{align*}
where $\tau>0$ is a sufficiently small constant, $\xi_{i}\in \mR^{m}, \kappa_{i}(\cdot)\in L^{2}_{\mF}(-\theta,0;\mR^{n}) $ $\upsilon_{i}(\cdot),\iota_{ij},\bar{\iota}_{ij}\in \\ L^{2}_{\mF}(0,T;\mR^{n}) ,i=1,2, j=1,2,\cdots, d,$       $\sH$ is an   suitable matrix and $F_{i}(\cdot), H_{i}(\cdot),\bar{H}_{i}(\cdot), \tilde{H}_{i}(\cdot), \\ \ddot{H}_{i}(\cdot), \bar{F}_{i}(\cdot),\tilde{F}_{i}(\cdot), \ddot{F}_{i}(\cdot), $  $ K_{ij}(\cdot),\bar{K}_{ij}(\cdot), D_{ij}(\cdot),\bar{D}_{ij}(\cdot), i=1,2, j=1,2,\cdots, d,$ are suitable\\ matrix-valued functions. The detailed definitions  will be provided below.  Set $(i=1,2)$
\begin{align}\label{501}
\nu_{i}(t)=
\begin{cases}
&\sH\xi_{i}+x_{0}, t=0,\\
&\kappa_{i}(t), t\in [-\theta, 0).
\end{cases}
\end{align}
  For any $(\nu_{i}(\cdot), u_{i}(\cdot)) \in L_{\mathbb{F}}^{2}(-\theta, 0;\mathbb{R}^n)\times L_{\mathbb{F}}^{2}(-\theta, T;\mathbb{R}^k)$,   as a special case of Theorem \ref{tt} in Section 4, there exists a unique solution $x_{i}(\cdot) \equiv x_{i}(\cdot; \nu_{i}, u_{i}) \in S_{\mathbb{F}}^{2}(-\theta,T;\mathbb{R}^n),i=1,2.$ for Eqs.\eqref{4}-\eqref{4d}.
 When later studying the control problem corresponding to the System \eqref{4}-\eqref{4d}, $\nu_{1}$ and $\nu_{2}$ will become two controlled elements. Obviously, our initial control elements are path-dependent, rather than a point as  in  \cite{LiuNiuWangyu2026}.

$\mathbf{Linear-convex\, problem}$($\mathbf{Problem\,(LC)}$):  When we study this linear-convex problem,  compared to conventional convex control problems,  our value functional  has the following form($\tau >0$ is a sufficiently small constant):
\begin{align}\label{5}
J(\nu_{1},\nu_{2}, u_{1},u_{2})&:=g_{11}(\xi_{1}+\tau\xi_{2})+ g_{12}(\xi_{2}+\tau\xi_{1})\no\\
&+\mE\bigg\{\bigg.g_{21}(x_{1}(T)+x_{2}(T) ) +g_{22}( x_{1}(T)+x_{2}(T) ) \no \\
 &  + \int^{T}_{0}g_{31}(t,x_{1}(t) )\dif t +\int^{T}_{0}g_{32}(t,x_{2}(t))\dif t\no\\
 &  + \int^{T}_{0}\tilde{g}_{31}(t-\theta,x_{1}(t-\theta) )\dif t +\int^{T}_{0}\tilde{g}_{32}(t-\theta,x_{2}(t-\theta))\dif t\no\\
 &  + \int^{T}_{0}\bar{g}_{31}(t,\mE[x_{1}(t)] )\dif t +\int^{T}_{0}\bar{g}_{32}(t,\mE[x_{2}(t)])\dif t\no\\
 &  + \int^{T}_{0}\ddot{g}_{31}(t,\mE^{\sF^{0}_{t}}[x_{1}(t)] )\dif t +\int^{T}_{0}\ddot{g}_{32}(t,\mE^{\sF^{0}_{t}}[x_{2}(t)])\dif t\no\\
 &  + \int^{T}_{T-\theta}\langle \tau \tilde{H}_{1}(t+\theta)^{\top}\mE^{\sF_{t}}[\eta_{1}(t+\theta)],  u_{1}(t)\rangle\dif t\no\\
&+ \int^{T}_{T-\theta}\langle \tau \tilde{H}_{2}(t+\theta)^{\top}\mE^{\sF_{t}}[\eta_{2}(t+\theta)],  u_{2}(t)\rangle\dif t\no\\
 &+\int^{T}_{0}g_{41}(t,u_{1} (t) )\dif t +\int^{T}_{0}g_{42}(t,  u_{2}(t))\dif t\no\\
  &+\int^{0}_{-\theta}\bar{g}_{41}(t,u_{1} (t) )\dif t +\int^{0}_{-\theta}\bar{g}_{42}(t,  u_{2}(t))\dif t\no\\
  &+\int^{T}_{T-\theta}\langle \tilde{F}_{2}(t+\theta)^{\top}\mE^{\sF_{t}}[\eta_{2}(t+\theta)]+\tilde{F}_{1}(t+\theta)^{\top}\mE^{\sF_{t}}[\eta_{1}(t+\theta)], x_{1}(t) \rangle\dif t\no\\
  &+\int^{T}_{T-\theta}\langle \tilde{F}_{2}(t+\theta)^{\top}\mE^{\sF_{t}}[\eta_{1}(t+\theta)], x_{2}(t) \rangle\dif t \bigg\}\bigg.,
\end{align}
where $g_{2i}: \Omega \times \mathbb{R}^n \rightarrow [0,\infty),$  $g_{3i},\bar{g}_{3i},\tilde{g}_{3i}, \ddot{g}_{3i}: \Omega \times [0,T] \times \mathbb{R}^n \rightarrow [0,\infty)$  and $g_{1i}: \mathbb{R}^m \rightarrow [0, \infty)$ and $g_{4i},\bar{g}_{4i}: \Omega \times [0, T] \times \mathbb{R}^k \rightarrow [0, \infty).$  We  require $g_{ji},\bar{g}_{4i},\bar{g}_{3i},\tilde{g}_{3i},\ddot{g}_{3i}, j=1,2,3,4, i=1,2$ to satisfy  convexity,  which will be presented  later.  We intend to find a quartet of $(\nu^{*}_{1},\nu^{*}_{2},u_{1}^{*},u_{2}^{*})$ such that
\begin{align}\label{5+}
J(\nu^{*}_{1},\nu^{*}_{2},u_{1}^{*},u_{2}^{*}):=\inf_{(\nu_{1},\nu_{2}, u_{1},u_{2})\in L^{2}_{\mF}(-\theta, 0; \mR^{n})\times L^{2}_{\mF}(-\theta, 0; \mR^{n}) \times L^{2}_{\mF}(-\theta, T;\mR^{k})\times L^{2}_{\mF}(-\theta, T; \mR^{k})}J(\nu_{1},\nu_{2}, u_{1},u_{2}).
\end{align}
In the following,  for simplicity, \eqref{5+} is denoted by
\begin{align}\label{5+-+}
J(\nu^{*}_{1},\nu^{*}_{2},u_{1}^{*},u_{2}^{*}):=\inf_{(\nu_{1},\nu_{2}, u_{1},u_{2})}J(\nu_{1},\nu_{2}, u_{1},u_{2}).
\end{align}
We establish the existence and uniqueness of optimal controls and derive their explicit closed-form representations. Further details for linear-convex problem will be elaborated in the following sections.

$\mathbf{Linear-quadratic\, problem\, with\, input\, constraints}\,( \mathbf{Problem\,(LQ-IC)}):$ In this case, we retain the linear controlled system \eqref{4}-\eqref{4d} and reformulate the performance criterion \eqref{5} into the following quadratic form:

\begin{align}\label{6}
 &\mJ(\nu_{1},\nu_{2}, u_{1},u_{2}):=\frac{1}{2}\langle A_{1}\xi_{1}, \xi_{1}\rangle+\frac{1}{2}\langle A_{2}\xi_{2}, \xi_{2}\rangle,\no\\
 &+\frac{1}{2}\mE\bigg\{\bigg. \langle C_{1}(x_{1}(T)+x_{2}(T)), x_{1}(T)+x_{2}(T)\rangle + \langle C_{2}(x_{1}(T)+x_{2}(T)), x_{1}(T)+x_{2}(T)\rangle  \no \\
 &  + \int^{T}_{0}\langle R_{1}(t)x_{1}(t), x_{1}(t)\rangle\dif t+ \int^{T}_{0}\langle R_{2}(t)x_{2}(t), x_{2}(t)\rangle\dif t\no\\
 &  + \int^{T}_{0}\langle \bar{R}_{1}(t)\mE[x_{1}(t)], \mE[x_{1}(t)]\rangle\dif t+ \int^{T}_{0}\langle \bar{R}_{2}(t)\mE[x_{2}(t)], \mE[x_{2}(t)]\rangle\dif t\no\\
 &  + \int^{T}_{0}\langle \ddot{R}_{1}(t)\mE^{\sF^{0}_{t}}[x_{1}(t)], \mE^{\sF^{0}_{t}}[x_{1}(t)]\rangle\dif t+ \int^{T}_{0}\langle \ddot{R}_{2}(t)\mE^{\sF^{0}_{t}}[x_{2}(t)], \mE^{\sF^{0}_{t}}[x_{2}(t)]\rangle\dif t\no\\
&  + \int^{T}_{0}\langle \tilde{R}_{1}(t-\theta)x_{1}(t-\theta), x_{1}(t-\theta)\rangle\dif t+ \int^{T}_{0}\langle \tilde{R}_{2}(t-\theta)x_{2}(t-\theta), x_{2}(t-\theta)\rangle\dif t\bigg\}\bigg.\no\\
 &  +\mE\bigg\{\bigg. \int^{T}_{T-\theta}\langle \tau \tilde{H}_{1}(t+\theta)^{\top}\mE^{\sF_{t}}[\eta_{1}(t+\theta)],  u_{1}(t)\rangle\dif t\no\\
&+ \int^{T}_{T-\theta}\langle \tau \tilde{H}_{2}(t+\theta)^{\top}\mE^{\sF_{t}}[\eta_{2}(t+\theta)],  u_{2}(t)\rangle\dif t\no\\
 &+\int^{T}_{0}\langle I_{1}(t) u_{1}(t), u_{1}(t)\rangle\dif t +\int^{T}_{0}\langle I_{2}(t)u_{2}(t), u_{2}(t)\rangle\dif t\no\\
 &+\int^{0}_{-\theta}\langle I_{1}(t) u_{1}(t), u_{1}(t)\rangle\dif t +\int^{0}_{-\theta}\langle I_{2}(t)u_{2}(t), u_{2}(t)\rangle\dif t\no\\
 &+\int^{T}_{T-\theta}\langle \tilde{F}_{2}(t+\theta)^{\top}\mE^{\sF_{t}}[\eta_{2}(t+\theta)]+\tilde{F}_{1}(t+\theta)^{\top}\mE^{\sF_{t}}[\eta_{1}(t+\theta)], x_{1}(t) \rangle\dif t\no\\
  &+\int^{T}_{T-\theta}\langle \tilde{F}_{2}(t+\theta)^{\top}\mE^{\sF_{t}}[\eta_{1}(t+\theta)], x_{2}(t) \rangle\dif t \bigg\}\bigg.,
\end{align}
where $A_{i}, C_{i},I_{i}(\cdot), R_{i}(\cdot), \bar{R}_{i}(\cdot), \tilde{R}_{i}(\cdot), \ddot{R}_{i}(\cdot)$ are suitable matrix-valued functions, which will be presented in the following.
Unlike Problem (LC), we now constrain the control variables (inputs) $\nu_{i},i=1,2$ and $u_{i},i=1,2$ to reside within nonempty, closed, convex constraint sets:
\begin{equation}\label{as1}
  \bar{\mU}_0 \subset \mathbb{R}^m,    \mU_{0}(\cdot) \equiv \{\mU_{0}(t) \subset \mathbb{R}^n, t \in  [-\theta, 0]\} \text{ and } \mU(\cdot) \equiv \{\mU(\omega,t) \subset \mathbb{R}^k, (\omega,t) \in \Omega \times [-\theta, T]\}
   \end{equation}
satisfying certain conditions, respectively. Let
\begin{equation}\label{as}
    \begin{cases}
    \sU_{0} := \{\nu(\cdot) \in L^2(-\theta, 0; \mathbb{R}^n) \mid \nu(0)=\sH \xi +x_{0}\in \bar{\mU}_{0}, \mbox{for}\,  \xi\in \mR^{m}; \\
     \quad\quad\quad\quad\quad\quad\quad\quad\quad\quad\quad\quad\quad\quad\quad\quad\quad\quad\quad\quad\quad\nu(t) \in \mU_{0}( t), \mbox{for}\, t \in  [-\theta,0)\},\\
    \sU := \{u(\cdot) \in L^2_{\mF}(-\theta, T;\mathbb{R}^k) \mid u(\omega,t) \in \mU(\omega, t) \text{ for almost all } (\omega, t) \in \Omega \times [-\theta,T]\}
   \end{cases}
    \end{equation}
and refer to $\sU_0\times \sU_0 \times \sU\times \sU$ as the admissible control set. Our linear-quadratic stochastic optimal control problem with input constraints (IC) is presented as follows.

\textbf{Problem (LQ-IC):} Determine a quartet of admissible control inputs $(\nu^*_{1},\nu^*_{2}, u^*_{1}(\cdot),u^*_{2}(\cdot))$ such that
\begin{equation}
    \mJ(\nu^*_{1},\nu^*_{2}, u^*_{1},u^*_{2})= \inf_{(\nu_{1},\nu_{2}, u_{1},u_{2}) \in \mU_0\times \mU_0 \times \sU\times \sU} \mJ(\nu_{1},\nu_{2}, u_{1},u_{2}).
    \label{eq:9}
\end{equation}

In this scenario, $(\nu^*_{1},\nu^*_{2}, u^*_{1}, u^*_{2})$, $x^*_{i}(\cdot) := x(\cdot; \nu^*_{i}, u^*_{i}), i=,2,$, and $(\nu^*_{1},\nu^*_{2}, u^*_{1},$ $u^*_{2}, x^{*}_{1}, x^{*}_{2} )$ are designated as a  quartet of optimal controls, the associated optimal state, and an optimal sextet for Problem (LQ-IC), respectively.

In summary, the main contributions are as follows:
\begin{itemize}
\item[$\bullet$] Based on the work of \cite{BuckdahnLiPeng2009}, we consider an extended mean-field system which has two couplings, pertaining to the initial and terminal values of   the System $(\pi)$ while the systems in \cite{BuckdahnLiPeng2009} lack these two couplings.

\item[$\bullet$] Compared to     Liu et al. \cite{LiuNiuWangyu2026}, which didn't study the mean-field system, and Tian and Yu \cite{Tianyu2023}, which didn't consider nonlinear
domination-monotonicity conditions and the system they considered didn't contain common noise, delayed terms, anticipated terms, this paper investigates an extended mean-field system including common noise , delayed terms, anticipated terms with  nonlinear
domination-monotonicity conditions. This type of system contains more difficulty in forward-backward coupling construction.

\item[$\bullet$] By utilizing the well-posedness result of the extended MF-AFBSDDEswCN along with other refined analytical techniques, we investigate two classes of stochastic optimal control problems: a linear-convex problem and a linear-quadratic problem subject to input constraints that may be time-dependent and random. For each problem, we establish the existence and uniqueness of optimal controls and derive their explicit closed-form expressions

\item[$\bullet$]  The last  contribution  is the most crucial one.  The two couplings, pertaining to the initial and terminal values of   the system $(\pi)$  are   path-dependent,  rather than a point.  When we study the corresponding control problems, our initial control element is path-dependent, rather than a point as  in  \cite{LiuNiuWangyu2026}. This difference makes the value functional essentially different and more challenging to deal with.

\end{itemize}

\section{Notations}
  For $x, y \in \mR^{n},$ we use $|x |$  to denote the Euclidean norm of
$x,$  $\langle x, y\rangle$ to denote the Euclidean inner product.
 For $A\in \mR^{ n\times d},$     $|A |$ represents  $\sqrt{\mathrm{Tr} (AA^{\top})}.$    Let $\sP$ stand for the $\mathbb{F}$-progressively measurable $\sigma$-field and
$\sB(\mathbb{R}^n)$ stand for the Borel $\sigma$-field on $\mathbb{R}^n$. $\mS^{n}\subset \mR^{n\times n}$ represents the set of all the symmetrical matrices.
Moreover,   $1_{A}$ denotes the indicator function on set $A.$  Next, We intend to define the following Banach spaces of random vectors or stochastic processes $(p\geq 2).$
\begin{itemize}
\item[1)] $L^{2}_{\sF_{T}}(\Omega,\mR^{n})$ is the set of $\sF_{T}-$measuurable random vectors $\zeta: \Omega\rightarrow \mR^{n}$ such that $\|\zeta\|_{L^{2}_{\sF_{T}}(\mR^{n})}:=\{\mE[|\zeta|^{2}]\}^{\frac{1}{2}}<\infty.$
\item[2)] $L^{\infty}_{\sF_{T}}(\Omega, \mR^{n})$ is the set of $\sF_{T}-$measuurable random vectors $\zeta: \Omega\rightarrow \mR^{n}$ such that $\|\zeta\|_{L^{\infty}_{\sF_{T}}(\mR^{n})}:=\operatorname{ess sup}_{\omega\in \Omega}|\zeta|<\infty.$
\item[3)] $L^{2}_{\mF}(0,T; \mR^{n})$ is the set of $\sP-$measurable stochastic processes $\varphi: \Omega\times [0,T]\rightarrow \mR^{n}$ such that $\|\varphi\|_{L^{2}_{\mF}(\mR^{n})}:=\bigg\{\bigg. \mE\int^{T}_{0}|\varphi(t)|^{2}\dif t   \bigg\}\bigg.^{\frac{1}{2}}<\infty.$

\item[4)] $L^{2}(0,T; \mR^{n})$ is the set of $\sB([0,T])-$measurable  deterministic function $\varphi:  [0,T]\rightarrow \mR^{n}$ such that $\|\varphi\|_{L^{2}(\mR^{n})}:=\bigg\{\bigg. \int^{T}_{0}|\varphi(t)|^{2}\dif t   \bigg\}\bigg.^{\frac{1}{2}}<\infty.$

\item[5)] $S^{2}_{\mF}(0,T;\mR^{n})$ is the set of   continuous  stochastic processes $\varphi: \Omega\times [0,T]\rightarrow \mR^{n}$ such that $\|\varphi\|_{S^{2}_{\mF}(\mR^{n})}:=\bigg\{\bigg. \mE[\sup_{0\leq t\leq T}|\varphi(t)|^{2}] \bigg\}\bigg.^{\frac{1}{2}}<\infty.$

\item[6)] $L^{\infty}_{\mF}(0,T;\mR^{n})$  consists of all  $\sP-$measurable stochastic processes  such that $\|\psi\|_{L^{\infty}_{\mF}(\mR^{n})}:=\,\mbox{esssup}_{ t\in [0,T]}|\psi|<\infty.$

\item[7)] $L^{\infty}(0,T;\mR^{n})$  consists of all deterministic function such that $\|\psi\|_{L^{\infty}(\mR^{n})}:=\sup_{ t\in [0,T]}|\psi|<\infty.$

\end{itemize}

Furthermore, for the  sake of simplicity, we set
\begin{itemize}
\item[1)] $\sM_{\mF}(\mR^{n+n+nd+nd}):=L^{2}_{\mF}(0,T;\mR^{n})\times L^{2}_{\mF}(0,T;\mR^{n})\times L^{2}_{\mF}(0,T;\mR^{nd})\times L^{2}_{\mF}(0,T;\mR^{nd}) $ equipped with the norm $\|\alpha(\cdot)\|_{\sM(\mR^{n+n+nd+nd})}:=\{\mE[\Xi_{\alpha}]\}^{\frac{1}{2}}$ for any\\ $\alpha(\cdot):=(\varphi_{1}(\cdot)^{\top}, \varphi_{2}(\cdot)^{\top},\varphi_{3}(\cdot)^{\top},\varphi_{4}(\cdot)^{\top})^{\top}\in \sM_{\mF}(\mR^{n+n+nd+nd})$ where
\begin{align}\label{99}
\Xi_{\alpha}:=\int^{T}_{0}|\varphi_{1}(t)|^{2}\dif t  +\int^{T}_{0}|\varphi_{2}(t)|^{2}\dif t+\int^{T}_{0}|\varphi_{3}(t)|^{2}\dif t.
\end{align}

\item[2)]$\mM_{\mF}(\mR^{n+n+nd+nd}):=S^{2}_{\mF}(0,T;\mR^{n})\times S^{2}_{\mF}(0,T;\mR^{n})\times L^{2}_{\mF}(0,T;\mR^{nd})\times L^{2}_{\mF}(0,T;\mR^{nd}) $ equipped with the norm $\|W(\cdot)\|_{\mM_{\mF}(\mR^{n+n+nd+nd})}:=\{\mE[\blacklozenge_{W}]\}^{\frac{1}{2}}$ for any\\ $W(\cdot):=(x(\cdot)^{\top}, y(\cdot)^{\top}, z^{0}(\cdot)^{\top}, z^{1}(\cdot)^{\top})^{\top}\in \mM_{\mF}(\mR^{n+n+nd+nd})$ where
\begin{align}\label{301+++}
\blacklozenge_{W}:=\sup_{0\leq t\leq T}|x(t)|^{2}  +\sup_{0\leq t\leq T}|y(t)|^{2}\dif t+\int^{T}_{0}|z^{0}(t)|^{2}\dif t+\int^{T}_{0}|z^{1}(t)|^{2}\dif t.
 \end{align}

\end{itemize}

 \section{Assumptions} To begin, in order to prove the existence and uniqueness result,  we introduce the following assumptions on the coefficients.  The inspiration for the following  assumptions partially  comes from \cite{LiuNiuWangyu2026}.\\
$\mathbf{Assumption\, 1}:$
\begin{itemize}
\item[ (i)]$\eta_{i}\in L^{2}(T,T+\theta;\mathbb{R}^{n }).$  Moreover, $\Psi_{i},$ $\ddot{\Psi}_{i},$ $ \Phi_{i}$ and $\aleph_{i},i=1,2$  are each measurable relative to the $\sigma-$algebras $\sB(\mR^{n})\times \sB(\mR^{n}),$ $\sB([-\theta,0])\times\sB(\mR^{n})\times \sB(\mR^{n}),$ $ \sF_{T}\times\sB(\mR^{n})\times \sB(\mR^{n})$ and \\ $\sP\times\sB(\mR^{n+n+n+n+n+n+n+n+n+n+n+n+n+n+nd+nd}),$ respectively. Furthermore, $\Phi_{i}(0,0)\in L^{2}_{\sF_{T}}(\Omega, \mR^{n}), \aleph_{i}(\cdot,0)\in \sM_{\mF}(\mR^{n+n+nd+nd}),i=1,2,$  and there exists a positive constant $L$ such that for any $t\in [0,T], \omega\in \Omega, ,x_{1},y_{1},x_{2},y_{2},x_{3},y_{3},x_{4},y_{4},x_{5},y_{5},x_{6},y_{6}\in \mR^{n}, $
   $$|\Psi_{1}(0,y_{2})|+|\Psi_{2}(y_{1},0)|+|l(t,x_{1},y_{1},x_{2},y_{2},x_{3},y_{3},x_{4},y_{4},x_{5},y_{5},x_{6},y_{6},0,0,0,0)|\leq L,$$
  where $l:=b_{i},\sigma^{k}_{i},f_{i},i=1,2, k=0, 1.$

\item[(ii)](Lipschitz conditions) there exist two  constants\\
 $L_{b}>0,L_{\sigma}>0,L_{f}>0,L_{\Phi}>0 ,L_{\Psi}>0,L_{\ddot{\Psi}}>0 $ and two small enough constants $\varepsilon>0, \epsilon>0$ such that, for\\
$t\in [0,T],$ $x_{1},x_{2}, \bar{x}_{1},\bar{x}_{2},x_{3},x_{4}, \bar{x}_{3},\bar{x}_{4},x_{5},x_{6}, \bar{x}_{5},\bar{x}_{6}, x,\bar{x}, y_{1}, y_{2},$ $y_{3}, y_{4}, y_{5},y_{6}, \bar{y}_{1},\bar{y}_{2},\bar{y}_{3},\bar{y}_{4},\\
\bar{y}_{5},\bar{y}_{6},y,\bar{y}\in \mR^{n},$  $z:=(z_{1}^{\top},z_{2}^{\top},\cdots,z_{d}^{\top})^{\top}\in \mR^{nd},\bar{z}:=(\bar{z}_{1}^{\top},\bar{z}_{2}^{\top},\cdots,\bar{z}_{d}^{\top})^{\top}\in \mR^{nd},$\\
$\tilde{z}:=(\tilde{z}_{1}^{\top},\tilde{z}_{2}^{\top},\cdots,\tilde{z}_{d}^{\top})^{\top}\in \mR^{nd},\tilde{\bar{z}}:=(\tilde{\bar{z}}_{1}^{\top},\tilde{\bar{z}}_{2}^{\top},\cdots,\tilde{\bar{z}}_{d}^{\top})^{\top}\in \mR^{nd},$
\begin{align*}
&|b_{i}(t,x_{1},y_{1},x_{2},y_{2},x_{3},y_{3},x_{4},y_{4},x_{5},y_{5},x_{6},y_{6},
x,y,z,\tilde{z} )\\
&-b_{i}(t,\bar{x}_{1},\bar{y}_{1},\bar{x}_{2},\bar{y}_{2},\bar{x}_{3},\bar{y}_{3},\bar{x}_{4},\bar{y}_{4},\bar{x}_{5},\bar{y}_{5},\bar{x}_{6},\bar{y}_{6},
\bar{x},\bar{y},\bar{z},\tilde{\bar{z}} )|\\
&\leq\varepsilon|x_{1}-\bar{x}_{1}|+\varepsilon|y_{1}-\bar{y}_{1}|+\varepsilon|x_{2}-\bar{x}_{2}|+\varepsilon|y_{2}-\bar{y}_{2}|\\
&+\varepsilon|x_{3}-\bar{x}_{3}|+\varepsilon|y_{3}-\bar{y}_{3}|+\varepsilon|x_{4}-\bar{x}_{4}|+\varepsilon|y_{4}-\bar{y}_{4}|
\\
&+\varepsilon|x_{5}-\bar{x}_{5}|+\varepsilon|y_{5}-\bar{y}_{5}|+\varepsilon|x_{6}-\bar{x}_{6}|+\varepsilon|y_{6}-\bar{y}_{6}|
\\
&+L_{b}|x-\bar{x}|+L_{b}|y-\bar{y}|+L_{b}|z-\bar{z}|+L_{b}|\tilde{z}-\tilde{\bar{z}}|, i=1,2.
\end{align*}

\begin{align*}
&|\sigma_{i}(t,x_{1},y_{1},x_{2},y_{2},x_{3},y_{3},x_{4},y_{4},x_{5},y_{5},x_{6},y_{6},
x,y,z,\tilde{z} )\\
&-\sigma_{i}(t,\bar{x}_{1},\bar{y}_{1},\bar{x}_{2},\bar{y}_{2},\bar{x}_{3},\bar{y}_{3},\bar{x}_{4},\bar{y}_{4},\bar{x}_{5},\bar{y}_{5},\bar{x}_{6},\bar{y}_{6},
\bar{x},\bar{y},\bar{z},\tilde{\bar{z}} )|\\
&\leq\varepsilon|x_{1}-\bar{x}_{1}|+\varepsilon|y_{1}-\bar{y}_{1}|+\varepsilon|x_{2}-\bar{x}_{2}|+\varepsilon|y_{2}-\bar{y}_{2}|\\
&+\varepsilon|x_{3}-\bar{x}_{3}|+\varepsilon|y_{3}-\bar{y}_{3}|+\varepsilon|x_{4}-\bar{x}_{4}|+\varepsilon|y_{4}-\bar{y}_{4}|
\\
&+\varepsilon|x_{5}-\bar{x}_{5}|+\varepsilon|y_{5}-\bar{y}_{5}|+\varepsilon|x_{6}-\bar{x}_{6}|+\varepsilon|y_{6}-\bar{y}_{6}|
\\
&+L_{\sigma}|x-\bar{x}|+L_{\sigma}|y-\bar{y}|+L_{\sigma}|z-\bar{z}|+L_{\sigma}|\tilde{z}-\tilde{\bar{z}}|, i=1,2.
\end{align*}

\begin{align*}
&|f_{i}(t,x_{1},y_{1},x_{2},y_{2},x_{3},y_{3},x_{4},y_{4},x_{5},y_{5},x_{6},y_{6},
x,y,z,\tilde{z} )\\
&-f_{i}(t,\bar{x}_{1},\bar{y}_{1},\bar{x}_{2},\bar{y}_{2},\bar{x}_{3},\bar{y}_{3},\bar{x}_{4},\bar{y}_{4},\bar{x}_{5},\bar{y}_{5},\bar{x}_{6},\bar{y}_{6},
\bar{x},\bar{y},\bar{z},\tilde{\bar{z}} )|\\
&\leq\varepsilon|x_{1}-\bar{x}_{1}|+\varepsilon|y_{1}-\bar{y}_{1}|+\varepsilon|x_{2}-\bar{x}_{2}|+\varepsilon|y_{2}-\bar{y}_{2}|\\
&+\varepsilon|x_{3}-\bar{x}_{3}|+\varepsilon|y_{3}-\bar{y}_{3}|+\varepsilon|x_{4}-\bar{x}_{4}|+\varepsilon|y_{4}-\bar{y}_{4}|
\\
&+\varepsilon|x_{5}-\bar{x}_{5}|+\varepsilon|y_{5}-\bar{y}_{5}|+\varepsilon|x_{6}-\bar{x}_{6}|+\varepsilon|y_{6}-\bar{y}_{6}|
\\
&+L_{f}|x-\bar{x}|+L_{f}|y-\bar{y}|+L_{f}|z-\bar{z}|+L_{f}|\tilde{z}-\tilde{\bar{z}}|, i=1,2.
\end{align*}

\begin{align*}
|\Phi_{1}(x_{1},x_{2})-\Phi_{1}(\bar{x}_{1},\bar{y}_{2})|\leq L_{\Phi}|x_{1}-\bar{x}_{1}|+\epsilon|x_{2}-\bar{x}_{2}|,\\
|\Phi_{2}(x_{1},x_{2})-\Phi_{2}(\bar{x}_{1},\bar{y}_{2})|\leq\epsilon|x_{1}-\bar{x}_{1}|+L_{\Phi}|x_{2}-\bar{x}_{2}|.
\end{align*}

\begin{align*}
|\Psi_{1}(y_{1},y_{2} )-\Psi_{1}(\bar{y}_{1},\bar{y}_{2})|\leq L_{\Psi}|y_{1}-\bar{y}_{1}|+\epsilon|y_{2}-\bar{y}_{2}|,\\
|\Psi_{2}(y_{1},y_{2})-\Psi_{2}(\bar{y}_{1},\bar{y}_{2})|\leq\epsilon|y_{1}-\bar{y}_{2}|+L_{\Psi}|y_{2}-\bar{y}_{2}|.
\end{align*}

\begin{align*}
|\ddot{\Psi}_{1}(t,y_{1},y_{2} )-\ddot{\Psi}_{1}(t,\bar{y}_{1},\bar{y}_{2})|\leq L_{\ddot{\Psi}}|y_{1}-\bar{y}_{1}|+L_{\ddot{\Psi}}|y_{2}-\bar{y}_{2}|,\\
|\ddot{\Psi}_{2}(t,y_{1},y_{2})-\ddot{t,\Psi}_{2}(\bar{y}_{1},\bar{y}_{2})|\leq L_{\ddot{\Psi}}|y_{1}-\bar{y}_{2}|+L_{\ddot{\Psi}}|y_{2}-\bar{y}_{2}|.
\end{align*}

\item[(iii)] There exist two constants $L_{1}>0,$ $L_{2}>0,$  a matrix $\sH\in \mR^{n\times m},$   several matrix-valued functions $H_{i}(\cdot),\bar{H}_{i}(\cdot),\ddot{H}_{i}(\cdot)\in L^{\infty}(0,T;\mR^{n\times k}), $ $\tilde{H}_{i}(\cdot)\in L^{\infty}(-\theta, T;\mR^{n\times k}) $   $M(\cdot) \in L^{\infty}(0,T+\theta;\mathbb{R}^{n \times n})$  and $D_{i}(\cdot), \bar{D}_{i}(\cdot)\in L^{\infty}(0,T;\mR^{n\times k})$ with \\ $D_{i}(\cdot):=(D_{i1}(t)^{\top},D_{i2}(t)^{\top},\cdots,D_{di}(t)^{\top})^{\top},$ $\bar{D}_{i}(\cdot):=(\bar{D}_{i1}(t)^{\top},\bar{D}_{i2}(t)^{\top},\cdots,\bar{D}_{di}(t)^{\top})^{\top},$\\ four $\sB(\mR^{m})$-measurable mappings $\bar{\chi}_{i1},\bar{\chi}_{i2}: \mR^{m} \rightarrow\mR^{m}, $  two $\sB([-\theta,0])\times\sB(\mR^{n})-$measu\\-rable mappings $\ddot{\chi}_{i}: [-\theta,0]\times \mR^{n} \rightarrow\mR^{n}, $ and two $\sP\times \sB(\mR^{k})-$measurable mappings
    $\chi_{i}:\Omega\times [0,T]\times \mR^{k}\rightarrow \mR^{k}, i=1,2$ such that the following holds for \\ $x^{\prime}_{i},x^{\prime\prime}_{i},x^{\prime\prime\prime}_{i},x^{\prime\prime\prime\prime}_{i},x^{\prime\prime\prime\prime\prime}_{i},
    x^{\prime\prime\prime\prime\prime\prime}_{i},y^{\prime}_{i},y^{\prime\prime}_{i},
    y^{\prime\prime\prime}_{i},y^{\prime\prime\prime\prime}_{i},y^{\prime\prime\prime\prime\prime}_{i},y^{\prime\prime\prime\prime\prime\prime}_{i},x_{i},y_{i},y\in \mR^{n},z_{i},\tilde{z}_{i}\in \mR^{nd},$ \\ $\bar{x}^{\prime}_{i},\bar{x}^{\prime\prime}_{i},\bar{x}^{\prime\prime\prime}_{i},\bar{x}^{\prime\prime\prime\prime}_{i},
    \bar{x}^{\prime\prime\prime\prime\prime}_{i},\bar{x}^{\prime\prime\prime\prime\prime\prime}_{i},
    \bar{y}^{\prime}_{i},\bar{y}^{\prime\prime}_{i},
    \bar{y}^{\prime\prime\prime}_{i},\bar{y}^{\prime\prime\prime\prime}_{i},\bar{y}^{\prime\prime\prime\prime\prime}_{i},
    \bar{y}^{\prime\prime\prime\prime\prime\prime}_{i},\bar{x}_{i},\bar{y}_{i},\bar{y}\in \mR^{n},\bar{z}_{i},\tilde{\bar{z}}_{i}\in \mR^{nd},i=1,2$.
\begin{itemize}
\item[1)]Adjoint function
\begin{align*}
&\chi_{i}(\cdot,0)\in L^{2}_{\mF}(0,T;\mR^{k}),i=1,2,\no\\
&|\bar{\chi}_{ik}( v_{1})-\bar{\chi}_{ik}( v_{2})|<L_{2}|v_{1}-v_{2}|,i=1,2,k=1,2,\no\\
&|\ddot{\chi}_{i}(t, u_{1})-\ddot{\chi}_{i}(t, u_{2})|<L_{2}|u_{1}-u_{2}|,i=1,2,\no\\
&|\chi_{i}(t, u_{1})-\chi_{i}(t, u_{2})|<L_{2}|u_{1}-u_{2}|,i=1,2,\no\\
&\langle \bar{\chi}_{ik}(  v_{1})-\bar{\chi}_{ik}( v_{2}), v_{1}-v_{2} \rangle\leq -L_{3}|\bar{\chi}_{ik}( v_{1})-\bar{\chi}_{ik}( v_{2})|^{2},i=1,2,k=1,2,\no\\
&\langle \chi_{i}(t, u_{1})-\chi_{i}(t, u_{2}), u_{1}-u_{2} \rangle\leq -L_{3}|\chi_{i}( t,u_{1})-\chi_{i}( t,u_{2})|^{2}, i=1,2
\end{align*}
for any $(\omega,t)\in \Omega\times [0,T], v_{1},v_{2}\in \mR^{m},u_{1},u_{2}\in \mR^{k}.$

\item[2)]Domination conditions( $a_{i}, c_{i}, i=1,2,$ are non-negative constants, $\tau$ is a sufficiently small positive constant):
\begin{align*}
&|\Psi_{1}(y_{1},y_{2})-\Psi_{1}(\bar{y}_{1},\bar{y}_{2})|\leq L_{2}|\bar{\chi}_{11}(\frac{\sH^{\top}y_{1}-\tau \sH^{\top}y_{2}}{1-\tau^{2}})-\bar{\chi}_{11}(\frac{\sH^{\top}\bar{y}_{1}-\tau \sH^{\top}\bar{y}_{2}}{1-\tau^{2}})|\no\\
&\quad \quad\quad\quad \quad\quad\quad \quad\quad\quad\quad+L_{2}|\bar{\chi}_{12}(\frac{\sH^{\top}y_{2}-\tau \sH^{\top}y_{1}}{1-\tau^{2}})-\bar{\chi}_{12}(\frac{\sH^{\top}\bar{y}_{2}-\tau \sH^{\top}\bar{y}_{1}}{1-\tau^{2}})|,\no\\
&|\Psi_{2}(y_{1},y_{2})-\Psi_{2}(\bar{y}_{1},\bar{y}_{2})|\leq L_{2}|\bar{\chi}_{21}(\frac{\sH^{\top}y_{2}-\tau \sH^{\top}y_{1}}{1-\tau^{2}})-\bar{\chi}_{11}(\frac{\sH^{\top}\bar{y}_{2}-\tau \sH^{\top}\bar{y}_{1}}{1-\tau^{2}})|\no\\
&\quad \quad\quad\quad \quad\quad\quad \quad\quad\quad\quad+L_{2}|\bar{\chi}_{22}(\frac{\sH^{\top}y_{1}-\tau \sH^{\top}y_{2}}{1-\tau^{2}})-\bar{\chi}_{22}(\frac{\sH^{\top}\bar{y}_{1}-\tau \sH^{\top}\bar{y}_{2}}{1-\tau^{2}})|.
\end{align*}

\begin{align*}
&|\ddot{\Psi}_{1}(t,y_{1},y_{2})-\ddot{\Psi}_{1}(t,\bar{y}_{1},\bar{y}_{2})|
\leq L_{2}|\ddot{\chi}_{1}(a_{1}M_{1}(t+\theta)^{\top}y_{1}+c_{1}M_{2}(t+\theta)^{\top}y_{2})\\
&\quad\quad\quad\quad\quad\quad\quad\quad\quad\quad\quad\quad\quad\quad\quad-\ddot{\chi}_{1}(a_{1}M_{1}(t+\theta)^{\top}\bar{y}_{1}]
+c_{1}M_{2}(t+\theta)^{\top}\bar{y}_{2})|,\\
&|\ddot{\Psi}_{2}(t,y_{1},y_{2})-\ddot{\Psi}_{2}(t,\bar{y}_{1},\bar{y}_{2})|
\leq L_{2}|\ddot{\chi}_{2}(a_{2}M_{1}(t+\theta)^{\top}y_{1}+c_{2}M_{2}(t+\theta)^{\top}y_{2})\\
&\quad\quad\quad\quad\quad\quad\quad\quad\quad\quad\quad\quad\quad\quad\quad-\ddot{\chi}_{2}(a_{2}M_{1}(t+\theta)^{\top}\bar{y}_{1}]
+c_{2}M_{2}(t+\theta)^{\top}\bar{y}_{2})|.
\end{align*}

\begin{align*}
&|l_{1}(t,x^{\prime}_{1},y^{\prime}_{1},x^{\prime\prime}_{1},y^{\prime\prime}_{1}
,x^{\prime\prime\prime}_{1},y^{\prime\prime\prime}_{1},x^{\prime\prime\prime\prime}_{1},y^{\prime\prime\prime\prime}_{1},
x^{\prime\prime\prime\prime\prime}_{1},y^{\prime\prime\prime\prime\prime}_{1},x^{\prime\prime\prime\prime\prime\prime}_{1},
y^{\prime\prime\prime\prime\prime\prime}_{1}
,x_{1},y_{1},z_{1},\tilde{z}_{1})\no\\
&-l_{1}(t,x^{\prime}_{1},y^{\prime}_{2},x^{\prime\prime}_{1},y^{\prime\prime}_{1},x^{\prime\prime\prime}_{1},y^{\prime\prime\prime}_{2}
,x^{\prime\prime\prime\prime}_{1},y^{\prime\prime\prime\prime}_{1},x^{\prime\prime\prime\prime\prime}_{1},y^{\prime\prime\prime\prime\prime}_{2},
x^{\prime\prime\prime\prime\prime\prime}_{1},y^{\prime\prime\prime\prime\prime\prime}_{1} ,x_{1},y_{2},z_{2},\tilde{z}_{2})|\no\\
&\leq L_{2}(|\chi_{1}(t,H_{1}(t)^{\top}y_{1}+\tau \bar{H}_{1}(t)^{\top}y^{\prime}_{1}+\tau \ddot{H}_{1}(t)^{\top}y^{\prime\prime\prime}_{1} +\tau \tilde{H}_{1}(t+\theta)^{\top}y^{\prime\prime\prime\prime\prime}_{1}\no\\
&\quad \quad\quad\quad \quad\quad\quad \quad\quad\quad \quad\quad\quad \quad\quad\quad \quad\quad \quad \quad\quad +D_{1}(t)^{\top}z_{1}+\bar{D}_{1}(t)^{\top}\tilde{z}_{1})\no\\
&\quad \quad\quad -\chi_{1}(t,H_{1}(t)^{\top}y_{2}+\tau \bar{H}_{1}(t)^{\top}y^{\prime}_{2}+\tau \ddot{H}_{1}(t)^{\top}y^{\prime\prime\prime}_{2} +\tau \tilde{H}_{1}(t+\theta)^{\top}y^{\prime\prime\prime\prime\prime}_{2}\no\\
&\quad \quad\quad\quad \quad\quad\quad \quad\quad\quad \quad\quad\quad \quad\quad\quad \quad\quad \quad \quad\quad +D_{1}(t)^{\top}z_{2}+\bar{D}_{1}(t)^{\top}\tilde{z}_{2}).
\end{align*}

\begin{align*}
&|l_{2}(t,x^{\prime}_{1},y^{\prime}_{1},x^{\prime\prime}_{1},y^{\prime\prime}_{1}
,x^{\prime\prime\prime}_{1},y^{\prime\prime\prime}_{1},x^{\prime\prime\prime\prime}_{1},y^{\prime\prime\prime\prime}_{1},
x^{\prime\prime\prime\prime\prime}_{1},y^{\prime\prime\prime\prime\prime}_{1},x^{\prime\prime\prime\prime\prime\prime}_{1},
y^{\prime\prime\prime\prime\prime\prime}_{1}
,x_{1},y_{1},z_{1},\tilde{z}_{1})\no\\
&-l_{2}(t,x^{\prime}_{1},y^{\prime}_{1},x^{\prime\prime}_{1},y^{\prime\prime}_{2},x^{\prime\prime\prime}_{1},y^{\prime\prime\prime}_{1},
x^{\prime\prime\prime\prime}_{1},y^{\prime\prime\prime\prime}_{2},
x^{\prime\prime\prime\prime\prime}_{1},y^{\prime\prime\prime\prime\prime}_{1},
x^{\prime\prime\prime\prime\prime\prime}_{1},y^{\prime\prime\prime\prime\prime\prime}_{2} ,x_{1},y_{2},z_{2},\tilde{z}_{2})|\no\\
&\leq L_{2}(|\chi_{2}(t,H_{2}(t)^{\top}y_{1}+\tau \bar{H}_{2}(t)^{\top}y^{\prime\prime}_{1}+\tau \ddot{H}_{2}(t)^{\top}y^{\prime\prime\prime\prime}_{1} +\tau \tilde{H}_{2}(t+\theta)^{\top}y^{\prime\prime\prime\prime\prime\prime}_{1}\no\\
&\quad \quad\quad\quad \quad\quad\quad \quad\quad\quad \quad\quad\quad \quad\quad\quad \quad\quad \quad \quad\quad +D_{2}(t)^{\top}z_{1}+\bar{D}_{2}(t)^{\top}\tilde{z}_{1})\no\\
&\quad \quad\quad -\chi_{2}(t,H_{2}(t)^{\top}y_{2}+\tau \bar{H}_{2}(t)^{\top}y^{\prime\prime}_{2}+\tau \ddot{H}_{2}(t)^{\top}y^{\prime\prime\prime\prime}_{2} +\tau \tilde{H}_{2}(t+\theta)^{\top}y^{\prime\prime\prime\prime\prime\prime}_{2}\no\\
&\quad \quad\quad\quad \quad\quad\quad \quad\quad\quad \quad\quad\quad \quad\quad\quad \quad\quad \quad \quad\quad +D_{2}(t)^{\top}z_{2}+\bar{D}_{2}(t)^{\top}\tilde{z}_{2}),\no\\
&l_{1}:=b_{1},\sigma^{0}_{1},\sigma^{1}_{1}, l_{2}:=b_{2},\sigma^{0}_{2},\sigma^{1}_{2}.
\end{align*}

\item[3)] Monotonicity conditions:  set
$$W_{1}=(x_{1}, y_{1},z_{1},\tilde{z}_{1}), W_{2}=(x_{2}, y_{2},z_{2},\tilde{z}_{2}),\bar{W}_{1}=(\bar{x}_{1}, \bar{y}_{1},\bar{z}_{1},\tilde{\bar{z}}_{1}), \bar{W}_{2}=(\bar{x}_{2}, \bar{y}_{2},\bar{z}_{2},\tilde{\bar{z}}_{2}),$$ $$\mathfrak{U}_{i}:=(x^{\prime}_{1},y^{\prime}_{1},x^{\prime}_{2},y^{\prime}_{2},x^{\prime}_{3},y^{\prime}_{3},x^{\prime}_{4},y^{\prime}_{4},
x^{\prime}_{5},y^{\prime}_{5},x^{\prime}_{6},y^{\prime}_{6},x_{i},y_{i},z_{i},\tilde{z}_{i}), $$  $$\bar{\mathfrak{U}}_{i}:=(\bar{x}^{\prime}_{1},\bar{y}^{\prime}_{1},\bar{x}^{\prime}_{2},\bar{y}^{\prime}_{2},\bar{x}^{\prime}_{3},
\bar{y}^{\prime}_{3},\bar{x}^{\prime}_{4},\bar{y}^{\prime}_{4},\bar{x}^{\prime}_{5},
\bar{y}^{\prime}_{5},\bar{x}^{\prime}_{6},\bar{y}^{\prime}_{6}
,\bar{x}_{i},\bar{y}_{i},\bar{z}_{i},\tilde{\bar{z}}_{i}),i=1,2. $$
 We assume that

\begin{align}
\begin{cases}
&\langle \Psi_{1}(y_{1},y)-\Psi_{1}(\bar{y}_{1},y), y_{1}-\bar{y}_{1} \rangle \leq -L_{3}|\bar{\chi}_{11}(\frac{\sH^{\top}y_{1}-\tau \sH^{\top}y}{1-\tau^{2}})-\bar{\chi}_{11}(\frac{\sH^{\top}\bar{y}_{1}-\tau \sH^{\top}y}{1-\tau^{2}} )|^{2}\\
&\quad \quad\quad\quad \quad\quad\quad \quad\quad\quad\quad\quad\quad\quad-L_{3}|\bar{\chi}_{12}(\frac{\sH^{\top}y-\tau \sH^{\top}y_{1}}{1-\tau^{2}})-\bar{\chi}_{12}(\frac{\sH^{\top}y-\tau \sH^{\top}\bar{y}_{1}}{1-\tau^{2}} )|^{2},\\
&\langle \Psi_{2}(y,y_{2})-\Psi_{2}(y,\bar{y}_{2}), y_{2}-\bar{y}_{2} \rangle
\leq -L_{3}|\bar{\chi}_{21}(\frac{\sH^{\top}y_{2}-\tau \sH^{\top}y}{1-\tau^{2}})-\bar{\chi}_{21}(\frac{\sH^{\top}\bar{y}_{2}-\tau \sH^{\top}y}{1-\tau^{2}} )|^{2}\\
&\quad \quad\quad\quad \quad\quad\quad \quad\quad\quad\quad\quad\quad\quad-L_{3}|\bar{\chi}_{22}(\frac{\sH^{\top}y-\tau \sH^{\top}y_{2}}{1-\tau^{2}})-\bar{\chi}_{22}(\frac{\sH^{\top}y-\tau \sH^{\top}\bar{y}_{2}}{1-\tau^{2}} )|^{2},\\
& \langle\Phi_{1}(x_{1},x_{2})-\Phi_{1}(\bar{x}_{1},x_{2}), x_{1}-\bar{x}_{1} \rangle
+\langle \Phi_{2}(x_{1},x_{2})-\Phi_{2}(x_{1},\bar{x}_{2}), x_{2}-\bar{x}_{2} \rangle \geq 0,\\
&\langle \aleph_{1}(t,\mathfrak{U}_{1})-\aleph_{1}(t,\bar{\mathfrak{U}}_{1}), W_{1}-\bar{W}_{1} \rangle+\langle \aleph_{2}(t,\mathfrak{U}_{2})-\aleph_{2}(t, \bar{\mathfrak{U}}_{2}), W_{2}-\bar{W}_{2} \rangle\\
& \leq -L_{3}|\chi_{1}(t,H_{1}(t)^{\top}y_{1}+\tau \bar{H}_{1}(t)^{\top}y^{\prime}_{1}+\tau \ddot{H}_{1}(t)^{\top}y^{\prime}_{3}+\tau \tilde{H}_{1}(t+\theta)^{\top}y^{\prime}_{5}\\
&\quad \quad\quad\quad \quad\quad\quad \quad\quad\quad\quad\quad\quad\quad\quad\quad\quad\quad\quad\quad+D_{1}(t)^{\top}z_{1}+\bar{D}_{1}(t)^{\top}\tilde{z}_{1})\\
&\quad\quad\quad-\chi_{1}(t,H_{1}(t)^{\top}\bar{y}_{1}+\tau \bar{H}_{1}(t)^{\top}\bar{y}^{\prime}_{1}+\tau \ddot{H}_{1}(t)^{\top}\bar{y}^{\prime}_{3}+ \tau \tilde{H}_{1}(t+\theta)^{\top}\bar{y}^{\prime}_{5}\\
&\quad \quad\quad\quad \quad\quad\quad \quad\quad\quad\quad\quad\quad\quad\quad\quad\quad\quad\quad\quad+D_{1}(t)^{\top}\bar{z}_{1}+\bar{D}_{1}(t)^{\top}\tilde{\bar{z}}_{1})|\\
&\quad -L_{3}|\chi_{2}(t,H_{2}(t)^{\top}y_{2}+\tau \bar{H}_{2}(t)^{\top}y^{\prime}_{2}+\tau \ddot{H}_{2}(t)^{\top}y^{\prime}_{4}+\tau \tilde{H}_{2}(t+\theta)^{\top}y^{\prime}_{6}\\
&\quad \quad\quad\quad \quad\quad\quad \quad\quad\quad\quad\quad\quad\quad\quad\quad\quad\quad\quad\quad+D_{2}(t)^{\top}z_{2}+\bar{D}_{2}(t)^{\top}\tilde{z}_{2})\\
&\quad\quad\quad-\chi_{2}(t,H_{2}(t)^{\top}\bar{y}_{2}+\tau \bar{H}_{2}(t)^{\top}\bar{y}^{\prime}_{2}+\tau \ddot{H}_{2}(t)^{\top}\bar{y}^{\prime}_{4}+ \tau \tilde{H}_{2}(t+\theta)^{\top}\bar{y}^{\prime}_{6}\\
&\quad \quad\quad\quad \quad\quad\quad \quad\quad\quad\quad\quad\quad\quad\quad\quad\quad\quad\quad\quad+D_{2}(t)^{\top}\bar{z}_{2}+\bar{D}_{2}(t)^{\top}\tilde{\bar{z}}_{2})|.
\end{cases}
\end{align}

\end{itemize}

\end{itemize}

When investigating  the Problem (LC), we will outline the strict assumptions necessary for the coefficients within the controlled system  \eqref{4}-\eqref{4d}.\\
$\mathbf{Assumption\, 2:}$\\
 $F_{i}(\cdot), \bar{F}_{i}(\cdot), \ddot{F}_{i}(\cdot), K_{ij}(\cdot), \bar{K}_{ij}(\cdot)\in L^{\infty}(0,T;\mathbb{R}^{n \times n}), H_{i}(\cdot), \bar{H}_{i}(\cdot), \ddot{H}_{i}(\cdot)\in L^{\infty}(0,T;\mathbb{R}^{n \times k}),$\\
  $\tilde{F}_{i}(\cdot)\in L^{\infty}(0,T+\theta;\mathbb{R}^{n \times n}),  \tilde{H}_{i}(\cdot)\in L^{\infty}(0,T+\theta;\mathbb{R}^{n \times k}),$$ D_{ij}(\cdot), \bar{D}_{ij}(\cdot)\in L^{\infty}(0,T;\mathbb{R}^{n \times k}),$\\ $\sH \in \mathbb{R}^{n \times m},$
 $\upsilon_{i}(\cdot) \in L_{\mathbb{F}}^{2}(0,T;\mathbb{R}^n),$ $\iota_{ij}(\cdot),  \bar{\iota}_{ij}(\cdot)\in L_{\mathbb{F}}^{2}(0,T;\mathbb{R}^n)$, and $x_0 \in \mathbb{R}^n$ for any $j = 1, 2, \dots, d, i=1,2$.  Moreover, there exists a sufficiently small constant $\tau_{1}> 0$ such that $\sup_{t\in [0,T]}|\bar{F}_{i}(t)|<\tau_{1},\sup_{t\in [0,T+K]}|\tilde{F}_{i}(t)|<\tau_{1},\sup_{t\in [0,T]}|\ddot{F}_{i}(t)|<\tau_{1},i=1,2. $

 We give the notations.
  $$K_{i}(\cdot) := (K_{i1}(\cdot)^{\top}, K_{i2}(\cdot)^{\top}, \dots, K_{id}(\cdot)^{\top})^{\top},$$
  $$\bar{K}_{i}(\cdot) := (\bar{K}_{i1}(\cdot)^{\top}, \bar{K}_{i2}(\cdot)^{\top}, \dots, \bar{K}_{id}(\cdot)^{\top})^{\top},$$

  $$D_{i}(\cdot):= (D_{i1}(\cdot)^{\top}, D_{i2}(\cdot)^{\top}, \dots, D_{id}(\cdot)^{\top})^{\top},$$

 $$\bar{D}_{i}(\cdot):= (\bar{D}_{i1}(\cdot)^{\top}, \bar{D}_{i2}(\cdot)^{\top}, \dots, \bar{D}_{id}(\cdot)^{\top})^{\top},$$

  $$\iota_{i}(\cdot) := (\iota_{i1}(\cdot)^{\top}, \iota_{i1}(\cdot)^{\top}, \dots, \iota_{id}(\cdot)^{\top})^{\top}.$$

  $$\bar{\iota}_{i}(\cdot) := (\bar{\iota}_{i1}(\cdot)^{\top}, \bar{\iota}_{i1}(\cdot)^{\top}, \dots, \bar{\iota}_{id}(\cdot)^{\top})^{\top}.$$

For a function $h:\mR^{n}\supset D\rightarrow \mR,$    we denote that
$$\nabla h(x) := \left( \frac{\partial h}{\partial x_1}(x), \frac{\partial h}{\partial x_2}(x), \dots, \frac{\partial h}{\partial x_n}(x) \right)^\top.$$
Now,  we give the assumptions for value functional. \\
$\mathbf{Assumption\, 3:}$
\begin{itemize}
\item[1)]  $g_{2i}(\cdot)$ and $g_{3i}( t,\cdot),$  $\bar{g}_{3i}( t,\cdot),$ $\ddot{g}_{3i}( t,\cdot),$  $i=1,2$ are convex, moreover $\tilde{g}_{3i}( t,\cdot),$ $g_{1i}(\cdot)$ and $g_{4i}(t,\cdot),\bar{g}_{4i}(t,\cdot),i=1,2$ are uniformly convex with parameter $\delta > 0$ (see Definition \ref{d1} in the Appendix) for almost all $(\omega, t)\in \Omega \times [0,T]$.
  \item[2)] $g_{1i}(\cdot)$, $g_{2i}(\cdot)$, $g_{3i}(t, \cdot),$ $\bar{g}_{3i}(t, \cdot),$ $\tilde{g}_{3i}(t, \cdot)$ $\ddot{g}_{3i}(t, \cdot)$ and $g_{4i}(t, \cdot)$ are continuously differentiable for almost all $(\omega, t)$. Moreover, $\nabla g_{2i}(\cdot),$  $\nabla g_{3i}(t,\cdot),$  $\nabla \bar{g}_{3i}(t,\cdot),$ $\nabla \ddot{g}_{3i}(t,\cdot),$  $\nabla \tilde{g}_{3i}(t,\cdot), i=1,2$ are uniform Lipschitz continuous in $(\omega, t)$ and the Lipschitz constants of   $\nabla g_{2i}(\cdot)$ are small enough.
 \item[3)] $g_{2i},$  $g_{4i},$  $\bar{g}_{4i}$ are $\sF_{T}\times\sB(\mathbb{R}^n)$-measurable, $\sP \times \sB(\mathbb{R}^k)$-measurable and $\sB([-\theta, 0]) \times \sB(\mathbb{R}^k)$-measurable, respectively.
$g_{3i}, \bar{g}_{3i}, \ddot{g}_{3i},\tilde{g}_{3i}$ are $\sP \times \mathcal{B}(\mathbb{R}^n)$-measurable,
\item[4)]$g_{2i}(0)\in L_{\sF_{T}}(\mR),$ $g_{3i}(\cdot,0),$ $\bar{g}_{3i}(\cdot,0),\ddot{g}_{3i}(\cdot,0)\in L^{2}_{\mF}(0, T;\mathbb{R}), \tilde{g}_{3i}(\cdot,0)  \in L^{2}_{\mF}(-\theta, T;\mathbb{R}),g_{4i}(\cdot,0)\\\in L^{2}_{\mF}(0, T;\mathbb{R}), \bar{g}_{4i}(\cdot,0)\in L^{2}_{\mF}(-\theta, 0;\mathbb{R}), \nabla g_{2i}(0) \in L_{\mathcal{F}_T}^2(\mathbb{R}^n),$ $\nabla g_{3i}(\cdot,0), \nabla \bar{g}_{3i}(\cdot,0),\nabla\ddot{g}_{3i}(\cdot,0)\\\in L^{2}_{\mF}(0, T;\mathbb{R}^{n}),\nabla \tilde{g}_{3i}(\cdot,0)\in L^{2}_{\mF}(-\theta, T;\mathbb{R}^{n}),\\ \text{and}\,  \nabla g_{4i}(\cdot,0)\in L_{\mathbb{F}}^2(0, T;\mathbb{R}^k), \nabla\bar{g}_{4i}(\cdot,0) \in L_{\mathbb{F}}^2(-\theta, 0;\mathbb{R}^k),i=1,2.$

\end{itemize}

When we study the linear-quadratic problem, we need the following assumptions.\\
$\mathbf{Assumption\, 4:}$
\begin{enumerate}
    \item[(1)]$A_{i}\in \mathbb{S}^m, $ $\tilde{R}_{i}(\cdot)  \in L^{\infty}(-\theta, T;\mathbb{S}^n)$, $C_{i} \in L_{\mathcal{F}_T}^{\infty}(\mathbb{S}^n)$, $R_{i}(\cdot),\bar{R}_{i}(\cdot),\ddot{R}_{i}(\cdot)  \in L^{\infty}(0,T;\mathbb{S}^n)$\\ and $I_{i}(\cdot) \in L^{\infty}(-\theta, T;\mathbb{S}^k), i=1,2$.
    \item[(2)] There exists a constant $\delta > 0$ such that $A_{i} - \delta \mI_m$, $C_{i},$ $R_{i}(t)$, $\ddot{R}_{i}(t)$, $\bar{R}_{i}(t),$ \\$\tilde{R}_{i}(t)- \delta \mI_n$, and   $I_{i}(t) - \delta \mI_k$ are positive semidefinite $\dif P\times \dif t-a.e.,$  $i=1,2.$
\end{enumerate}
In what follows, $\mI_n$ denotes the $(n \times n)$ identity matrix.

Concerning the constrained sets $\bar{\mU}_{0}, \mU_{0}(\cdot)$ and $\mU(\cdot)$ (see \cite{DDD8}), we introduce the following assumption.\\
$\mathbf{Assumption\, 5:}$
\begin{enumerate}
    \item[(i)]  $\mU_0, \mU_{0}(\cdot)$ and $\mU(\omega, t)$ are nonempty, closed, and convex $\dif P\times\dif t-$a.e.
    \item[(ii)] The indicator function $(\omega, t, u) \mapsto 1_{\mU(\omega,t)}(u)$ is $\sP \times \sB(\mathbb{R}^k)$-measurable.
    \item[(iii)] There exists a process $a(\cdot) \in L^2_\mathbb{F}(-\theta,T; \mathbb{R}^k)$ such that $a(\omega, t) \in \mU(\omega, t)$ $\dif P\times\dif t-$a.e.
\end{enumerate}

By convention, we shall use the letter $C$ to represent a positive constant. This constant's value is determined solely by the constants specified in the Assumption and may differ from one situation to another.
\section{Well-posedness of  extended MF-FBSDE}

In this section, we examine the well-posedness of  System $(\pi)$.
The main results of this section  will be given  as follows.
\subsection{Prior estimates}
\bl\label{t1}
Assume   $\mathbf{Assumption 1}.$  If \\ $W_{i}(\cdot) \in \mM_{\mF}(\mR^{{n+n+nd+nd}}), i=1,2 $  is a solution to System $(\pi),$ where $W_{i}(\cdot):=(x_{i}(\cdot)^{\top},$ $y_{i}(\cdot)^{\top},$ $z^{0}_{i}(\cdot)^{\top}, z^{1}_{i}(\cdot)^{\top})^{\top}, i=1,2.$  Consequently,  the following estimate holds:
\begin{align}\label{11}
&\mE[\blacklozenge_{W_{1}}]+\mE[\blacklozenge_{W_{2}}]\no\\
&\leq C(|\Psi_{1}(0,0)|^{2}+|\Psi_{2}(0,0)|^{2}+|\ddot{\Psi}_{1}(0,0)|^{2}+|\ddot{\Psi}_{2}(0,0)|^{2}
 +\mE[ |\Phi_{1}(0,0)|^{2}+ |\Phi_{2}(0,0)|^{2}\no\\
 &\quad  \sup_{T\leq t\leq T+\theta}|\eta_{1}(t)|^{2}+ \sup_{T\leq t\leq T+\theta}|\eta_{2}(t)|^{2}+\Xi_{\aleph_{1}(\cdot,0)}+\Xi_{\aleph_{2}(\cdot,0)}]),
\end{align}
where $\Xi_{\aleph_{i}(\cdot,0)}, \blacklozenge_{W_{i}}, i=1,2$  are given by  \eqref{99} and \eqref{301+++}, respectively, $C$ is a constant only depending on  the constants in $\mathbf{Assumption 1}.$ Consider another set of coefficients $(\tilde{\Psi}_{i},\tilde{\ddot{\Psi}}_{i},\tilde{\Phi}_{i}, \tilde{\aleph}_{i}),\\i=1,2$ and let $\tilde{W}_{i}(\cdot):=(\tilde{x}_{i}(t)^{\top}(\cdot),\tilde{y}_{i}(\cdot)^{\top},\tilde{z}^{0}_{i}(\cdot)^{\top},\tilde{z}^{1}_{i}(\cdot)^{\top})^{\top}   \in \mM_{\mF}(\mR^{n+n+nd+nd}), i=1,2 $ be a solution to System $(\pi)$ associated with $(\tilde{\Psi}_{i},\tilde{\ddot{\Psi}}_{i}, \tilde{\Phi}_{i}, \tilde{\aleph}_{i}).$ Then the following estimate holds:
\begin{align}\label{12}
&\mE[\blacklozenge_{\hat{W}_{1}}]+\mE[\blacklozenge_{\hat{W}_{2}}]\no\\
&\leq C(\sup_{-\theta\leq t\leq 0}|\hat{\ddot{\Psi}}_{1}(\mE[\tilde{y}_{1}(t+\theta)],\mE[\tilde{y}_{2}(t+\theta)] )|^{2} +\sup_{-\theta\leq t\leq 0}|\hat{\ddot{\Psi}}_{2}(\mE[\tilde{y}_{1}(t+\theta)],\mE[\tilde{y}_{2}(t+\theta)] )|^{2}\no\\
 &\quad +|\hat{\Psi}_{1}(\tilde{y}_{1}(0),\tilde{y}_{2}(0))|^{2}+|\hat{\Psi}_{2}(\tilde{y}_{1}(0),\tilde{y}_{2}(0))|^{2}\no\\ &\quad+\mE[|\hat{\Phi}_{1}(\tilde{x}_{1}(T),\tilde{x}_{2}(T))|^{2}+|\hat{\Phi}_{2}(\tilde{x}_{1}(T),\tilde{x}_{2}(T))|^{2}
+\Xi_{\hat{\aleph}_{1}(\cdot, \tilde{\mathfrak{U}}_{1}(\cdot))}+\Xi_{\hat{\aleph}_{2}(\cdot, \tilde{\mathfrak{U}}_{2}(\cdot))}]),
\end{align}
where $\hat{l}=l-\tilde{l}, l=W_{i},\mathfrak{U}_{i},\Psi_{i},\Phi_{i}, \ddot{\Psi}_{i},\aleph_{i},$
\begin{align*}
&\tilde{\mathfrak{U}}_{i}(t):=\big(\big.(\mE[\tilde{x}_{1}(t)])^{\top}, (\mE[\tilde{y}_{1}(t)])^{\top}, (\mE[\tilde{x}_{2}(t)])^{\top}, (\mE[\tilde{y}_{2}(t)])^{\top}, \\
&\quad\quad\quad\quad\quad(\mE^{\sF^{0}_{t}}[\tilde{x}_{1}(t)])^{\top}, (\mE^{\sF^{0}_{t}}[\tilde{y}_{1}(t)])^{\top}, (\mE^{\sF^{0}_{t}}[\tilde{x}_{2}(t)])^{\top}, (\mE^{\sF^{0}_{t}}[\tilde{y}_{2}(t)])^{\top},\\
&\quad\quad \quad(\tilde{x}_{2}(t-\theta))^{\top}, (\mE^{\sF_{t}}[\tilde{y}_{2}(t+\theta)])^{\top}, (\tilde{x}_{2}(t-\theta))^{\top}, (\mE^{\sF_{t}}[\tilde{y}_{2}(t+\theta)])^{\top},\\
&\quad\quad \quad\quad\tilde{ x}_{i}(t)^{\top}, \tilde{y}_{i}(t)^{\top}, \tilde{z}_{i}(t)^{\top}\big)\big.^{\top}, i=1,2.
\end{align*}
\el
\begin{proof}
  To begin, we aim to prove   \eqref{11}.   Using a conventional estimate for  stochastic differential equations (see, e.g., \cite{zhang2017}), we derive
\begin{align}\label{13}
&\mE[\sup_{0\leq t\leq T}|x_{1}(t)|^{2}] + \mE[\sup_{0\leq t\leq T}|x_{2}(t)|^{2}] \no\\
&\leq C\bigg\{\bigg.   |\Psi_{1}(y_{1}(0),y_{2}(0))-\Psi_{1}(0,y_{2}(0))|^{2}
+|\Psi_{1}(0,y_{2}(0))-\Psi_{1}(0,0)|^{2}\no\\
&\quad \quad \quad\quad \quad +|\Psi_{2}(y_{1}(0),y_{2}(0))-\Psi_{2}(y_{1}(0), 0)|^{2}
+|\Psi_{2}(y_{1}(0),0)-\Psi_{2}(0,0)|^{2}\no\\
&\quad \quad \quad\quad \quad +|\Psi_{1}(0,0)|^{2}+ |\Psi_{2}(0,0)|^{2}+|\ddot{\Psi}_{1}(0,0)|^{2}+ |\ddot{\Psi}_{2}(0,0)|^{2}\no\\
&\quad \quad \quad\quad \quad + \varepsilon\mE[\sup_{0\leq t\leq T}|y_{1}(t)|^{2}]+ \varepsilon\mE[\sup_{0\leq t\leq T}|y_{2}(t)|^{2}]\no\\
&   +\mE\bigg[\bigg.\bigg(\bigg.\int^{T}_{0}|b_{1}(t,0)|\dif t  \bigg)\bigg.^{2}+\bigg(\bigg.\int^{T}_{0}|b_{2}(t,0)|\dif t  \bigg)\bigg.^{2}  +\bigg.\bigg(\bigg.\int^{T}_{0}|\sigma_{1}(t, 0)|\dif t  \bigg)\bigg.^{2}+\bigg(\bigg.\int^{T}_{0}|\sigma_{2}(t,0)|\dif t  \bigg)\bigg.^{2} \no \\
&   +\mE\bigg[\bigg.\bigg(\bigg.\int^{T}_{0}|b_{1}(t,0, \mE[y_{2}(t)],0,\mE^{\sF^{0}_{t}}[y_{2}(t)],0,\mE^{\sF_{t}}[y_{2}(t+\theta)], 0)-b_{1}(t,0)|\dif t  \bigg)\bigg.^{2} \no\\ &+\bigg(\bigg.\int^{T}_{0}|b_{2}(t,0,  \mE[y_{1}(t)],0,\mE^{\sF^{0}_{t}}[y_{1}(t)],0,\mE^{\sF_{t}}[y_{1}(t+\theta)], 0)-b_{2}(t,0)|\dif t  \bigg)\bigg.^{2}  \no\\
&    +\bigg.\bigg(\bigg.\int^{T}_{0}|\sigma_{1}(t,0, \mE[y_{2}(t)],0,\mE^{\sF^{0}_{t}}[y_{2}(t)],0,\mE^{\sF_{t}}[y_{2}(t+\theta)], 0)-\sigma_{1}(t, 0)|\dif t  \bigg)\bigg.^{2}\no\\ &+\bigg(\bigg.\int^{T}_{0}|\sigma_{2}(t,0,  \mE[y_{1}(t)],0,\mE^{\sF^{0}_{t}}[y_{1}(t)],0,\mE^{\sF_{t}}[y_{1}(t+\theta)], 0)-\sigma_{2}(t, 0)|\dif t  \bigg)\bigg.^{2} \no \\
 &   +\bigg(\bigg.\int^{T}_{0}|b_{1}(t, \Upsilon_{1}(t) )-\mE[b_{1}(t,0, \mE[y_{2}(t)],0,\mE^{\sF^{0}_{t}}[y_{2}(t)],0,\mE^{\sF_{t}}[y_{2}(t+\theta)], 0)]|\dif t  \bigg)\bigg.^{2}\no\\
 &    +\bigg(\bigg.\int^{T}_{0}|b_{2}(t,\Upsilon_{2}(t)) -\mE[b_{2}(t,0,  \mE[y_{1}(t)],0,\mE^{\sF^{0}_{t}}[y_{1}(t)],0,\mE^{\sF_{t}}[y_{1}(t+\theta)], 0)]|\dif t  \bigg)\bigg.^{2}   \no\\
&  +\int^{T}_{0}|\sigma_{1}(t,\Upsilon_{1}(t))-\sigma_{1}(t,0, \mE[y_{2}(t)],0,\mE^{\sF^{0}_{t}}[y_{2}(t)],0,\mE^{\sF_{t}}[y_{2}(t+\theta)], 0)|^{2}\dif t \no\\
 &   +\int^{T}_{0}|\sigma_{2}(t,\Upsilon_{2}(t))-\sigma_{2}(t,0,  \mE[y_{1}(t)],0,\mE^{\sF^{0}_{t}}[y_{1}(t)],0,\mE^{\sF_{t}}[y_{1}(t+\theta)], 0)]|^{2}\dif t  \bigg]\bigg.  \bigg\}\bigg.,
\end{align}
where
\begin{align*}
\Upsilon_{i}(t):=(0^{\top}, &(\mE[y_{1}(t)])^{\top}, 0^{\top}, (\mE[y_{2}(t)])^{\top},0^{\top},(\mE^{\sF^{0}_{t}}[y_{1}(t)])^{\top},0^{\top},(\mE^{\sF^{0}_{t}}[y_{2}(t)])^{\top},0^{\top}, \\
&(\mE^{\sF_{t}}[y_{1}(t+\theta)])^{\top}, 0^{\top},(\mE^{\sF_{t}}[y_{2}(t+\theta)])^{\top}, 0^{\top}, y_{i}(t)^{\top}, z^{0}_{i}(t)^{\top},z^{1}_{i}(t)^{\top} )^{\top},i=1,2.
\end{align*}
Consequently, by virtue of the domination and Lipschitz conditions established in Assumption 1 (and observing that the values  $\epsilon>0, \varepsilon>0$ are small enough), it follows that
\begin{align}\label{14}
&\mE[\sup_{0\leq t\leq T}|x_{1}(t)|^{2}] + \mE[\sup_{0\leq t\leq T}|x_{2}(t)|^{2}]\no \\
&\leq K\bigg\{\bigg.\diamond_{1}+ |\Psi_{1}(0,0)|^{2}+ |\Psi_{2}(0,0)|^{2}+ |\ddot{\Psi}_{1}(0,0)|^{2}+ |\ddot{\Psi}_{2}(0,0)|^{2}+ (\varepsilon+\epsilon) \mE[\sup_{0\leq t\leq T}|y_{1}(t)|^{2}]\no\\
&\quad \quad+ (\varepsilon+\epsilon)\mE[\sup_{0\leq t\leq T}|y_{2}(t)|^{2}]+\sup_{T<t\leq T+\theta}|\eta_{1}(t)| + \sup_{T<t\leq T+\theta}|\eta_{2}(t)|\no\\
&\quad \quad   +\mE\bigg[\bigg.\bigg(\bigg.\int^{T}_{0}|b_{1}(t,0)|\dif t  \bigg)\bigg.^{2} +\bigg(\bigg.\int^{T}_{0}|b_{2}(t, 0)|\dif t  \bigg)\bigg.^{2}  \bigg]\bigg.\no \\
& \quad \quad   +\mE\bigg[\bigg.\bigg(\bigg.\int^{T}_{0}|\sigma_{1}(t,0)|\dif t  \bigg)\bigg.^{2} +\bigg(\bigg.\int^{T}_{0}|\sigma_{2}(t, 0)|\dif t  \bigg)\bigg.^{2}  \bigg]\bigg. \bigg\}\bigg.,\no \\
\end{align}
where \begin{align}\label{15}
\diamond_{1}&:=\bigg|\bigg.\bar{\chi}_{11}(\frac{\sH^{\top}y_{1}(0)-\tau \sH^{\top}y_{2}(0)}{1-\tau^{2}})-\bar{\chi}_{11}(\frac{-\tau \sH^{\top}y_{2}(0)}{1-\tau^{2}})\bigg|\bigg.^{2}\no\\
&\quad+\bigg|\bigg.\bar{\chi}_{12}(\frac{\sH^{\top}y_{2}(0)-\tau \sH^{\top}y_{1}(0)}{1-\tau^{2}})-\bar{\chi}_{12}(\frac{ \sH^{\top}y_{2}(0)}{1-\tau^{2}})\bigg|\bigg.^{2}\no\\
&\quad+\bigg|\bigg.\bar{\chi}_{21}(\frac{\sH^{\top}y_{2}(0)-\tau \sH^{\top}y_{1}(0)}{1-\tau^{2}})-\bar{\chi}_{21}(\frac{ -\tau \sH^{\top}y_{1}(0)}{1-\tau^{2}})\bigg|\bigg.^{2}\no\\
&\quad+\bigg|\bigg.\bar{\chi}_{22}(\frac{\sH^{\top}y_{1}(0)-\tau \sH^{\top}y_{2}(0)}{1-\tau^{2}})-\bar{\chi}_{22}(\frac{ \sH^{\top}y_{1}(0)}{1-\tau^{2}})\bigg|\bigg.^{2}\no\\
&+\mE\int^{T}_{0}|\chi_{1}(t,H_{1}(t)^{\top}y_{1}(t)+\tau\bar{H}_{1}(t)^{\top}\mE[y_{1}(t)]+\tau\ddot{H}_{1}(t)^{\top}\mE^{\sF^{0}_{t}}[y_{1}(t)]\no\\
&\quad\quad\quad\quad\quad+\tau\tilde{H}_{1}(t+\theta)^{\top}\mE^{\sF_{t}}[y_{1}(t+\theta)]+D_{1}(t)^{\top}z^{0}_{1}(t)+D_{1}(t)^{\top}z^{1}_{1}(t))-\chi_{1}(t,0)|^{2}\dif t\no\\
&+\mE\int^{T}_{0}|\chi_{2}(t,H_{2}(t)^{\top}y_{2}(t)+\tau\bar{H}_{2}(t)^{\top}\mE[y_{2}(t)]+\tau\ddot{H}_{2}(t)^{\top}\mE^{\sF^{0}_{t}}[y_{2}(t)]\no\\
&\quad\quad\quad\quad\quad+\tau\tilde{H}_{2}(t+\theta)^{\top}\mE^{\sF_{t}}[y_{2}(t+\theta)]
+D_{2}(t)^{\top}z^{0}_{2}(t)+\bar{D}_{2}(t)^{\top}z^{1}_{2}(t))
-\chi_{2}(t,0)|^{2}\dif t.
\end{align}
 In addition, by employing the fundamental estimate for BSDEs (see, e.g., \cite{PardouxPeng1990}) and the Lipschitz conditions,  it holds that
\begin{align}\label{16}
&\mE\bigg[\bigg. \sup_{0\leq t\leq T+\theta}|y_{1}(t)|^{2}+\sup_{0\leq t\leq T+\theta}|y_{2}(t)|^{2}+ \int^{T}_{0}|z_{1}(t)|^{2}\dif t+\int^{T}_{0}|z_{2}(t)|^{2}\dif t\bigg]\bigg. \no\\
&\leq C\mE\bigg[\bigg.
\sup_{T<t\leq T+\theta}|\eta_{1}(t)| + \sup_{T<t\leq T+\theta}|\eta_{2}(t)|\no\\
&\quad \quad+|\Phi_{1}(0,0)|^{2}+ |\Phi_{2}(0,0)|^{2}
   +\bigg[\bigg.\bigg(\bigg.\int^{T}_{0}|f_{1}(t,0)|\dif t  \bigg)\bigg.^{2} \no\\
   & \quad \quad \quad\quad \quad \quad+\bigg(\bigg.\int^{T}_{0}|f_{2}(t, 0)|^{2}\dif t  \bigg)\bigg.^{2}+   \sup_{0\leq t\leq T}|x_{1}(t)|^{2}+\sup_{0\leq t\leq T}|x_{2}(t)|^{2}\bigg]\bigg..
\end{align}
The combination of \eqref{14} and \eqref{16} yield that
\begin{align}\label{17}
&\mE[\blacklozenge_{W_{1}}+\blacklozenge_{W_{2}}]\no\\
&\leq C\{ \diamond_{1}+|\Psi_{1}(0,0)|^{2}+ |\Psi_{2}(0,0)|^{2}+\ddot{|\Psi}_{1}(0,0)|^{2}+ |\ddot{\Psi}_{2}(0,0)|^{2}\no\\
  &\quad +\mE[ |\Phi_{1}(0,0)|^{2}+ |\Phi_{2}(0,0)|^{2}+   \Xi_{\aleph_{1}(\cdot, 0)}+ \Xi_{\aleph_{2}(\cdot, 0)}]\}.
\end{align}
Next, employing It\^{o}'s formula to $\langle x_{1}(\cdot), y_{1}(\cdot)\rangle, \langle x_{2}(\cdot), y_{2}(\cdot)\rangle$ and subsequently summing the results alongside the monotonicity conditions from Assumptions 1, allows us to obtain

\begin{align}\label{18}
&\mE[\langle \Phi_{1}(x_{1}(T), x_{2}(T)), x_{1}(T) \rangle]+\mE[\langle \Phi_{2}(x_{1}(T), x_{2}(T)), x_{2}(T) \rangle]]\no\\
&\leq\mE[\langle \Psi_{1}(y_{1}(0), y_{2}(0)), y_{1}(0) \rangle]+\mE[\langle \Psi_{2}(y_{1}(0), y_{2}(0)), y_{2}(0) \rangle]\no\\
& \quad + \mE\bigg[\bigg. \int^{T}_{0}\langle \aleph_{1}(t, \mathfrak{U}_{1}(t)), W_{1}(t)\rangle \dif t +\int^{T}_{0}\langle \aleph_{2}(t, \mathfrak{U}_{2}(t)), W_{2}(t)\rangle \dif t \no\\
&\leq \mE[\langle \Psi_{1}(y_{1}(0), y_{2}(0)), y_{1}(0) \rangle]+\mE[\langle \Psi_{2}(y_{1}(0), y_{2}(0)), y_{2}(0) \rangle]\no\\
& \quad + \mE\bigg[\bigg. \int^{T}_{0}\langle \aleph_{1}(t, 0), W_{1}(t)\rangle \dif t +\int^{T}_{0}\langle \aleph_{2}(t,0), W_{2}(t)\rangle \dif t   \bigg]\bigg.\no\\
&-L_{3}\mE\bigg[\bigg.\int^{T}_{0}|\chi_{1}(t,H_{1}(t)^{\top}y_{1}(t)+\tau\bar{H}_{1}(t)^{\top}\mE[y_{1}(t)]+\tau\ddot{H}_{1}(t)^{\top}\mE^{\sF_{t}^{0}}[y_{1}(t)]\no\\
&\quad\quad\quad\quad\quad\quad\quad+\tau\tilde{H}_{1}(t+\theta)^{\top}\mE^{\sF_{t}}[y_{1}(t+\theta)]+D_{1}(t)^{\top}z^{0}_{1}(t)+\bar{D}_{1}(t)^{\top}z^{1}_{1}(t))\no\\
&\quad\quad\quad\quad\quad\quad\quad\quad\quad\quad\quad\quad\quad\quad\quad\quad\quad\quad\quad\quad\quad\quad\quad\quad\quad\quad\quad\quad\quad\quad-\chi_{1}(t,0)|^{2}\dif t\no\\
&+\int^{T}_{0}|\chi_{2}(t,H_{2}(t)^{\top}y_{2}(t)+\tau\bar{H}_{2}(t)^{\top}\mE[y_{2}(t)]+\tau\ddot{H}_{2}(t)^{\top}\mE^{\sF_{t}^{0}}[y_{2}(t)]\no\\
&\quad\quad\quad\quad\quad\quad\quad+\tau\tilde{H}_{2}(t+\theta)^{\top}\mE^{\sF_{t}}[y_{2}(t+\theta)]+D_{2}(t)^{\top}z^{0}_{2}(t)+\bar{D}_{2}(t)^{\top}z^{1}_{2}(t))\no\\
&\quad\quad\quad\quad\quad\quad\quad\quad\quad\quad\quad\quad\quad\quad\quad\quad\quad\quad\quad\quad\quad\quad\quad\quad\quad\quad\quad\quad\quad\quad-\chi_{2}(t,0)|^{2}\dif t\bigg]\bigg..
\end{align}
where\\
 \begin{align*}
\mathfrak{U}_{i}(t):&=\big(\big. (\mE[x_{1}(t)])^{\top}, (\mE[y_{1}(t)])^{\top}, (\mE[x_{2}(t)])^{\top}, (\mE[y_{2}(t)])^{\top}, \\
 &\quad\quad\quad\quad\quad(\mE^{\sF^{0}_{t}}[x_{1}(t)])^{\top}, (\mE^{\sF^{0}_{t}}[y_{1}(t)])^{\top}, (\mE^{\sF^{0}_{t}}[x_{2}(t)])^{\top}, (\mE^{\sF^{0}_{t}}[y_{2}(t)])^{\top},\\
&\quad\quad \quad(x_{1}(t-\theta))^{\top}, (\mE^{\sF_{t}}[y_{1}(t+\theta)])^{\top}, (x_{2}(t-\theta))^{\top}, (\mE^{\sF_{t}}[y_{2}(t+\theta)])^{\top}, \\
  &\quad\quad \quad\quad
  x_{i}(t)^{\top}, y_{i}(t)^{\top}, z^{0}_{i}(t)^{\top}, z^{1}_{i}(t)^{\top}\big)\big.^{\top}, i=1,2.
\end{align*}
Applying the monotonicity conditions in  Assumption 1,   we have
\begin{align}\label{19}
&\mE[\langle \Phi_{1}(0, x_{2}(T)), x_{1}(T) \rangle]+\mE[\langle \Phi_{2}(x_{1}(T),0), x_{2}(T) \rangle]\no\\
&\leq \langle \Psi_{1}(0, y_{2}(0)), y_{1}(0) \rangle+\mE[\langle \Psi_{2}(y_{1}(0), 0), y_{2}(0) \rangle\no\\
& \quad + \mE\bigg[\bigg. \int^{T}_{0}\langle \aleph_{1}(t, 0), W_{1}(t)\rangle \dif t +\int^{T}_{0}\langle \aleph_{2}(t,0), W_{2}(t)\rangle \dif t \bigg]\bigg.\no\\
& \quad - L_{3}I_{1}.
\end{align}
By solving the above inequality and  substituting $I_{1}$ in to $\eqref{17},$    it holds that
\begin{align}\label{20}
&\mE[\blacklozenge_{W_{1}}+\blacklozenge_{W_{2}}]\no\\
&\leq C\{ |\Psi_{1}(0,0)|^{2}+ |\Psi_{2}(0,0)|^{2}+|\ddot{\Psi}_{1}(0,0)|^{2}+ |\ddot{\Psi}_{2}(0,0)|^{2}\no\\
&\quad   +\mE[ |\Phi_{1}(0,0)|^{2}+ |\Phi_{2}(0,0)|^{2}+   \Xi_{\aleph_{1}(\cdot, 0)}+ \Xi_{\aleph_{2}(\cdot, 0)}]\no\\
&\quad+\mE[\langle \Psi_{1}(0, y_{2}(0)), y_{1}(0) \rangle]+\mE[\langle \Phi_{2}(y_{1}(0), 0), y_{2}(0) \rangle]]\no\\
& \quad + \mE\bigg[\bigg. \int^{T}_{0}\langle \aleph_{1}(t, 0), W_{1}(t)\rangle \dif t +\int^{T}_{0}\langle \aleph_{2}(t, 0), W_{2}(t)\rangle \dif t\bigg]\bigg. \no\\
&\quad -\mE[\langle \Phi_{1}(0, x_{2}(T)), x_{1}(T) \rangle]-\mE[\langle \Phi_{2}(x_{1}(T),0), x_{2}(T) \rangle].
\end{align}
By  a standard calculus(noting that $\epsilon,\varepsilon$ are small enough),  we derived
\begin{align}\label{20}
&\mE[\blacklozenge_{W_{1}}+\blacklozenge_{W_{2}}]\no\\
&\leq K_{1}\{ |\Psi_{1}(0,0)|^{2}+ |\Psi_{2}(0,0)|^{2}+|\ddot{\Psi}_{1}(0,0)|^{2}+ |\ddot{\Psi}_{2}(0,0)|^{2}\no\\
 &\quad\quad\quad\quad\quad\quad + \mE[ |\Phi_{1}(0,0)|^{2}+ |\Phi_{2}(0,0)|^{2}+ \Xi_{\aleph_{1}(\cdot, 0)} +\Xi_{\aleph_{2}(\cdot, 0)} ]\}\no\\
&\quad +\frac{1}{2}\mE[\blacklozenge_{W_{1}}+\blacklozenge_{W_{2}}].
\end{align}
The desired result for \eqref{11} is derived  .

Next, we intend to prove the second result \eqref{12}.  Set
$$\grave{\Psi}_{i}(t,y_{1},y_{2}):=\Psi_{i}(y_{1}+\tilde{y}_{1}(0), y_{2}+\tilde{y}_{2}(0))-\tilde{\Psi}_{i}(\tilde{y}_{1}(0), \tilde{y}_{2}(0)), i=1,2,$$
$$\grave{\ddot{\Psi}}_{i}(t,y_{1},y_{2}):=\ddot{\Psi}_{i}(t,y_{1}+\mE[\tilde{y}_{1}(t+\theta)], y_{2}+\mE[\tilde{y}_{2}(t+\theta)])-\tilde{\ddot{\Psi}}_{i}(t,\mE[\tilde{y}_{1}(t+\theta)], \mE[\tilde{y}_{2}(t+\theta)]), i=1,2,$$
$$\grave{\Phi}_{i}(x_{1},x_{2}):=\Phi_{i}(x_{1}+\tilde{x}_{1}(T), x_{2}+\tilde{x}_{2}(T))-\tilde{\Phi}_{i}(\tilde{x}_{1}(T), \tilde{x}_{2}(T)), i=1,2, $$
$$\grave{\aleph}_{i}(t,\mathfrak{U}):=\aleph_{i}(t,\mathfrak{U}+\tilde{\mathfrak{U}}_{i}(t))-\tilde{\aleph}_{i}(t,\tilde{\mathfrak{U}}_{i}(t)), i=1,2, $$
for any $\mathfrak{U}\in \mR^{n+n+n+n+n+n+n+n+n+n+n+n+n+n+nd+nd}.$ One can readily verify that Assumption 1  remains valid under the new set of coefficients $(\grave{\Psi}_{i},\grave{\Phi}_{i},\grave{\aleph}_{i})$ with the same constants, $H, H_{i}, \bar{H}_{i},\tilde{H}_{i} , D_{i}, \bar{D}_{i},i=1,2$ in Assumption 1  and the new adjoint functions:

$$\bar{\grave{\chi}}_{11}(v):=\bar{\chi}_{11}(v+\frac{\sH^{\top}\tilde{y}_{1}(0)-\tau \sH^{\top}\tilde{y}_{2}(0)}{1-\tau^{2}}),
\bar{\grave{\chi}}_{12}(v):=\bar{\chi}_{12}(v+\frac{\sH^{\top}\tilde{y}_{2}(0)-\tau \sH^{\top}\tilde{y}_{1}(0)}{1-\tau^{2}}),$$
$$\bar{\grave{\chi}}_{21}(v):=\bar{\chi}_{21}(v+\frac{\sH^{\top}\tilde{y}_{2}(0)-\tau \sH^{\top}\tilde{}y_{1}(0)}{1-\tau^{2}}),
\bar{\grave{\chi}}_{22}(v):=\bar{\chi}_{22}(v+\frac{\sH^{\top}\tilde{y}_{1}(0)-\tau \sH^{\top}\tilde{y}_{2}(0)}{1-\tau^{2}}),$$

$$
\ddot{\grave{\chi}}_{1}(v):=\bar{\ddot{\chi}}_{1}(v+a_{1}M_{1}(t+\theta)^{\top}\mE[\tilde{y}_{1}(t+\theta)]+c_{1}M_{2}(t+\theta)^{\top}\mE[\tilde{y}_{2}(t+\theta)]),
$$

$$
\ddot{\grave{\chi}}_{2}(v):=\bar{\ddot{\chi}}_{2}(v+a_{2}M_{1}(t+\theta)^{\top}\mE[\tilde{y}_{1}(t+\theta)]+c_{2}M_{2}(t+\theta)^{\top}\mE[\tilde{y}_{2}(t+\theta)]),
$$

\begin{align*}
\grave{\chi}_{1}(t,u):=\chi_{1}(t,u+H_{1}(t)^{\top}\tilde{y}_{1}(t)&+\tau\bar{H}_{1}(t)^{\top}\mE[\tilde{y}_{1}(t)]+\tau\ddot{H}_{1}(t)^{\top}\mE^{\sF^{0}_{t}}[\tilde{y}_{1}(t)]\\+\tau\tilde{H}_{1}(t+\theta)^{\top}\mE^{\sF_{t}}[\tilde{y}_{1}(t+\theta)]
&+D_{1}(t)^{\top}\tilde{z}^{0}_{1}(t))+D_{1}(t)^{\top}\tilde{z}^{1}_{1}(t)),
\end{align*}

\begin{align*}
\grave{\chi}_{2}(t,u):=\chi_{2}(t,u+H_{2}(t)^{\top}\tilde{y}_{2}(t)&+\tau\bar{H}_{2}(t)^{\top}\mE[\tilde{y}_{2}(t)]+\tau\ddot{H}_{1}(t)^{\top}\mE^{\sF^{0}_{t}}[\tilde{y}_{1}(t)]\\+\tau\tilde{H}_{2}(t+\theta)^{\top}\mE^{\sF_{t}}[\tilde{y}_{2}(t+\theta)]
&+D_{2}(t)^{\top}\tilde{z}^{0}_{2}(t))+\bar{D}_{2}(t)^{\top}\tilde{z}^{1}_{2}(t)).
\end{align*}
Moreover, we also confirm that the process $\hat{W}_{i}(\cdot) := W_{i}(\cdot) - \tilde{W}_{i}(\cdot) \in \mM_{\mathbb{F}}(\mathbb{R}^{n+n+nd+nd})$ fulfills the MF-AFBSDDEswCN  with $(\grave{\Psi}_{i}, \grave{\ddot{\Psi}}_{i}, \grave{\Phi}_{i}, \grave{\aleph}_{i})$. Consequently, the estimate \eqref{11} applied to $\hat{W}_{i}(\cdot)$ results in \eqref{12}.

\end{proof}

\subsection{ Method of continuity}  We  now introduce the method of continuity and present a useful lemma that will serve to prove the existence and uniqueness result. Following Assumption 1, we define a new set of coefficients: $(\Psi^{0}_{i},\ddot{\Psi}^{0}_{i}, \Phi^{0}_{i}, \aleph^{0}_{i}), i=1,2$  with

\begin{align}\label{21}
&\Psi^{0}_{1}(y_{1},y_{2}):=\sH\bar{\chi}_{11}(\frac{\sH^{\top}y_{1}-\tau \sH^{\top}y_{2}}{1-\tau^{2}})+\sH\bar{\chi}_{12}(\frac{\sH^{\top}y_{2}-\tau \sH^{\top}y_{1}}{1-\tau^{2}}),\no\\
&\Psi^{0}_{2}(y_{1},y_{2}):=\sH\bar{\chi}_{21}(\frac{\sH^{\top}y_{2}-\tau \sH^{\top}y_{1}}{1-\tau^{2}})+\sH\bar{\chi}_{22}(\frac{\sH^{\top}y_{1}-\tau \sH^{\top}y_{2}}{1-\tau^{2}}),\no\\
&\ddot{\Psi}^{0}_{1}(t,y_{1},y_{2}):=\ddot{\chi}_{1}(a_{1}M_{1}(t+\theta)^{\top}y_{1}+c_{1}M_{2}(t+\theta)^{\top}y_{2}),\no\\
&\ddot{\Psi}^{0}_{2}(t,y_{1},y_{2}):=\ddot{\chi}_{2}(a_{2}M_{1}(t+\theta)^{\top}y_{1}+c_{2}M_{2}(t+\theta)^{\top}y_{2}),\no\\
&  \Phi^{0}_{i}(x_{1},x_{2}):=0, f^{0}_{i}(t,\mathfrak{U}_{i}):=0,\no\\
&b^{0}_{i}(t,\mathfrak{U}_{i}):=H_{i}(t)\chi_{i}(t, H_{i}(t)^{\top}y_{i}+\tau\bar{H}_{i}(t)^{\top}y'_{i}+\tau\ddot{H}_{i}(t)^{\top}y''_{i}+\tau\tilde{H}_{i}(t+\theta)^{\top}y'''_{i}\no\\
&\quad\quad\quad\quad\quad\quad\quad\quad\quad\quad\quad\quad\quad\quad\quad\quad\quad\quad\quad\quad\quad\quad\quad\quad+D_{i}(t)^{\top}z^{0}_{i}+\bar{D}_{i}(t)^{\top}z^{1}_{i}),\no\\
&\sigma^{0,0}_{i}(t,\mathfrak{U}_{i}):=D_{i}(t)\chi_{i}(t, H_{i}(t)^{\top}y_{i}+\tau\bar{H}_{i}(t)^{\top}y'_{i}+\tau\ddot{H}_{i}(t)^{\top}y''_{i}+\tau\tilde{H}_{i}(t+\theta)^{\top}y'''_{i}\no\\
&\quad\quad\quad\quad\quad\quad\quad\quad\quad\quad\quad\quad\quad\quad\quad\quad\quad\quad\quad\quad\quad\quad\quad\quad+D_{i}(t)^{\top}z_{i}+\bar{D}_{i}(t)^{\top}z^{\prime}_{i}),\no\\
 &\sigma^{1,0}_{i}(t,\mathfrak{U}_{i}):=\bar{D}_{i}(t)\chi_{i}(t, H_{i}(t)^{\top}y_{i}+\tau\bar{H}_{i}(t)^{\top}y'_{i}+\tau\ddot{H}_{i}(t)^{\top}y''_{i}+\tau\tilde{H}_{i}(t+\theta)^{\top}y'''_{i}\no\\
 &\quad\quad\quad\quad\quad\quad\quad\quad\quad\quad\quad\quad\quad\quad\quad\quad\quad\quad\quad\quad\quad\quad\quad\quad+D_{i}(t)^{\top}z_{i}+\bar{D}_{i}(t)^{\top}z^{\prime}_{i}),i=1,2,
\end{align}
for any
 \begin{align*}
 &\mathfrak{U}_{i}:=(x_{1}^{\prime\top},y_{1}^{\prime\top}, x_{2}^{\prime\top}, y_{2}^{\prime\top},x_{1}^{\prime\prime\top},y_{1}^{\prime\prime\top}, x_{2}^{\prime\prime\top}, y_{2}^{\prime\prime\top},x_{1}^{\prime\prime\prime\top},y_{1}^{\prime\prime\prime\top}, x_{2}^{\prime\prime\prime\top}, y_{2}^{\prime\prime\prime\top}, x_{i}^{\top},y_{i}^{\top},z_{i}^{\top})^{\top}\\
 &\in \mR^{n+n+n+n+n+n+n+n+n+n+n+n+n+n+nd+nd},i=1,2.
\end{align*}
 Within this framework, we use the notations \\$\aleph^{0}_{i}:=((f^{0}_{i})^{\top},(b^{0}_{i})^{\top}, (\sigma^{0,0}_{i})^{\top}, (\sigma^{1,0}_{i})^{\top})^{\top} $ and $ \sigma^{k,0}_{i}:=((\sigma^{k,0}_{i1})^{\top},(\sigma^{k,0}_{i2})^{\top},\cdots, (\sigma^{k,0}_{id})^{\top})^{\top},k=0,1$ in this context.  A simple verification confirms that the coefficients  $(\Psi^{0}_{i}, \ddot{\Psi}^{0}_{i},\Phi^{0}_{i}, \aleph^{0}_{i})$  satisfy Assumption 1 with  the same parameters(If needed, make suitable adjustments, ensuring they do not alter the overall proof).

For any $\xi_{i}\in \mR^{n}, \ddot{\xi}_{i}(\cdot)\in L^{2}(-\theta,0; \mR^{n}),  \zeta_{i} \in L^{2}_{\sF_{T}}(\Omega, \mR^{n}),$ and any $\gamma_{i}(\cdot):=(\phi_{i}(\cdot)^{\top}, \psi_{i}(\cdot)^{\top},\\\beta_{i}^{0}(\cdot)^{\top},\beta_{i}^{1}(\cdot)^{\top})^{\top}\in \sM_{\mF}(\mR^{n+n+nd+nd})$ and $\beta_{i}^{k}:=((\beta^{k}_{i1})^{\top},(\beta^{k}_{i2})^{\top},\cdots, (\beta^{k}_{id})^{\top})^{\top},i=1,2, k=0,1,$ we now proceed to introduce a family of MF-AFBSDDEswCN defined by the parameter $\alpha\in [0,1],$

\begin{align}\label{23}
\begin{cases}
&\dif x_{1}^{\alpha}(t)=[b_{1}^{\alpha}(t,\mathfrak{U}_{1}^{\alpha}(t))+\psi_{1}(t)]\dif t
+\sum^{d}_{j=1}[\sigma_{1j}^{0,\alpha}(t,\mathfrak{U}^{\alpha}_{1}(t))+\beta^{0}_{1j}(t)]\dif B^{0}_{j}(t)\\
&\quad\quad\quad\quad\quad\quad\quad\quad\quad\quad\quad\quad+\sum^{d}_{j=1}[\sigma_{1j}^{1,\alpha}(t,\mathfrak{U}^{\alpha}_{1}(t))+\beta^{1}_{1j}(t)]\dif B^{1}_{j}(t), t\in [0,T],\\
&\dif y^{\alpha}_{1}(t)^{\top}(t) =\{[f^{\alpha}_{1}(t,\mathfrak{U}_{1}^{\alpha}(t))+\phi_{1}(t)]\dif t+\sum^{d}_{j=1} z^{0,\alpha}_{1j}(t)\dif B^{0}_{j}(t)+\sum^{d}_{j=1} z^{1,\alpha}_{1j}(t)\dif B^{1}_{j}(t),\\
& x_{1}^{\alpha}(0)=\Psi_{1}^{\alpha}(y_{1}^{\alpha}(0),y_{2}^{\alpha}(0))+ \xi_{1},\\ &x_{1}^{\alpha}(t)=\ddot{\Psi}_{1}(t,\mE[y^{\alpha}_{1}(t+\theta)],\mE[y^{\alpha}_{2}(t+\theta)])+\ddot{\xi}_{1}(t), t\in [-\theta, 0),\\
  &y^{\alpha}_{1}(T)=\Phi_{1}^{\alpha}(x^{\alpha}_{1}(T),x^{\alpha}_{2}(T) )+\zeta_{1}, y^{\alpha}_{1}(t)=\eta_{1}(t), t\in (T, T+\theta],
\end{cases}
\end{align}

 \begin{align}\label{24}
\begin{cases}
&\dif x_{2}^{\alpha}(t)=[b_{2}^{\alpha}(t,\mathfrak{U}_{2}^{\alpha}(t))+\psi_{2}(t)]\dif t
 +\sum^{d}_{j=1}[\sigma_{1j}^{0,\alpha}(t,\mathfrak{U}^{\alpha}_{2}(t))+\beta^{0}_{2j}(t)]\dif B^{0}_{j}(t)\\
&\quad\quad\quad\quad\quad\quad\quad\quad\quad\quad\quad\quad+\sum^{d}_{j=1}[\sigma_{2j}^{1,\alpha}(t,\mathfrak{U}^{\alpha}_{2}(t))+\beta^{1}_{2j}(t)]\dif B^{1}_{j}(t),t\in [0,T],\\
&\dif y_{2 }^{\alpha}(t) =[f_{2}^{\alpha}(t,\mathfrak{U}_{2}^{\alpha}(t))+\phi_{2}(t)]\dif t+\sum^{d}_{j=1} z^{0,\alpha}_{2j}(t)\dif B^{0}_{j}(t)+\sum^{d}_{j=1} z^{1,\alpha}_{2j}(t)\dif B^{1}_{j}(t),\\
& x_{2 }^{\alpha}(0)=\Psi_{2}^{\alpha}(y_{1}^{\alpha}(0),y_{2}^{\alpha}(0))+ \xi_{2}, \\ &x_{2}^{\alpha}(t)=\ddot{\Psi}_{2}(t,\mE[y^{\alpha}_{1}(t+\theta)],\mE[y^{\alpha}_{2}(t+\theta)])+\ddot{\xi}_{2}(t), t\in [-\theta, 0),\\
 & y_{2}^{\alpha}(T)=\Phi_{2}^{\alpha}(x_{1}^{\alpha}(T),x_{2}^{\alpha}(T) )+\zeta_{2}, y^{\alpha}_{2}(t)=\eta_{2}(t), t\in (T, T+\theta],
\end{cases}
\end{align}
 where $(\Psi_{i}^{\alpha}, \ddot{\Psi}_{i}^{\alpha},\Phi_{i}^{\alpha},\aleph_{i}^{\alpha}):=\alpha(\Psi_{i}, \ddot{\Psi}_{i}, \Phi_{i}, \aleph_{i})+(1-\alpha)(\Psi_{i}^{0}, \ddot{\Psi}_{i}^{0}, \Phi_{i}^{0}, \aleph_{i}^{0}), i=1,2,$
  \begin{align*}
 \mathfrak{U}^{\alpha}_{i}(t):&=\big(\big. (\mE[x^{\alpha}_{1}(t)])^{\top}, (\mE[y^{\alpha}_{1}(t)])^{\top}, (\mE[x^{\alpha}_{2}(t)])^{\top}, (\mE[y^{\alpha}_{2}(t)])^{\top}, \\
&\quad\quad\quad\quad\quad(\mE^{\sF^{0}_{t}}[x^{\alpha}_{1}(t)])^{\top}, (\mE^{\sF^{0}_{t}}[y^{\alpha}_{1}(t)])^{\top}, (\mE^{\sF^{0}_{t}}[x^{\alpha}_{2}(t)])^{\top}, (\mE^{\sF^{0}_{t}}[y^{\alpha}_{2}(t)])^{\top},\\
&\quad\quad \quad(x^{\alpha}_{1}(t-\theta))^{\top}, (\mE^{\sF_{t}}[y^{\alpha}_{1}(t+\theta)])^{\top}, (x^{\alpha}_{2}(t-\theta))^{\top}, (\mE^{\sF_{t}}[y^{\alpha}_{2}(t+\theta)])^{\top}, \\
  &\quad\quad \quad\quad
  x^{\alpha}_{i}(t)^{\top}, y^{\alpha}_{i}(t)^{\top}, z^{0,\alpha}_{i}(t)^{\top},z^{1,\alpha}_{i}(t)^{\top}\big)\big.^{\top}, i=1,2.
  \end{align*}
 It is evident that the coefficients $(\Psi_{i}^{\alpha}+\xi_{i},\ddot{\Psi}_{i}^{\alpha}+\ddot{\xi}_{i}, \Phi_{i}^{\alpha}+\zeta_{i}, \aleph_{i}^{\alpha}+ \gamma_{i}), i=1,2$ of Eq.\eqref{23} and Eq.\eqref{24} fulfill $\mathrm{Assumption\, 1}.$
For the sake of simplicity,  Eq.\eqref{23} and Eq.\eqref{24} are collectively referred to as System $(\pi_{1}).$

Two extreme cases are observed.  In the case where $\alpha = 1$ and $(\xi_{i}, \ddot{\xi}_{i}(\cdot),\zeta_{i}, \gamma_{i}(\cdot)), i=1,2$   vanish,  System $(\pi_{1})$  reduces to System $(\pi)$, which we intend to study.  In the instance where $\alpha = 0,$  System $(\pi_{1})$ degenerates into a decoupled form, which can be solved using the established results for SDEs (see \cite{zhang2017}) and BSDEs (see \cite{PardouxPeng1990}).

\bl\label{t4}
Under Assumption 1 for the coefficients   $(\Psi_{i}, \ddot{\Psi}_{i}, \Phi_{i},  \aleph_{i}),i=1,2,$  there is an absolute constant $\delta_{0} > 0$ such that if  for some $\alpha_{0}\in [0,1),$ System $(\pi_{1})$ is uniquely solvable in $\mM_{\mF}(\mR^{n+n+nd+nd})\times \mM_{\mF}(\mR^{n+n+nd+nd}) $ for any $(\xi_{i}, \ddot{\xi}_{i}(\cdot),\zeta_{i}, \gamma_{i}(\cdot)) \in \mR^{n}\times  L^{2}(-\theta,0; \mR^{n})\times L^{2}_{\sF_{T}}(\mR^{n})\times \sM_{\mF}(\mR^{n+n+nd+nd}), i=1,2,$  then replacing $\alpha_0$ by any $\alpha \in (\alpha_0, (\alpha_0 + \delta_0) \wedge
1]$, the same conclusion remains true.

\el

\begin{proof}
Let $\delta_0 > 0$ be a value specified below, and let $\delta \in (0, \delta_0 \wedge (1 - \alpha_0)]$. Set $\alpha = \alpha_0 + \delta$ and $(\xi_{i}, \ddot{\xi}_{i}(\cdot), \zeta_{i}, \gamma_{i}(\cdot)) \in \mathbb{R}^n \times  L^{2}(-\theta,0; \mR^{n}) \times L^2_{\mathcal{F}_T}(\mathbb{R}^n) \times \sM_{\mathbb{F}}(\mathbb{R}^{n+n+nd+nd})$. For any\\
 $W_{i}(\cdot) = (x_{i}(\cdot)^{\top}, y_{i}(\cdot)^{\top}, z^{0}_{i}(\cdot)^{\top}, z^{1}_{i}(\cdot)^{\top})^{\top} \in \mM_{\mathbb{F}}(\mathbb{R}^{n+n+nd+nd})$,  set
$$\mathfrak{U}(\cdot):=(\mathfrak{U}_{1}(\cdot)^{\top},\mathfrak{U}_{2}(\cdot)^{\top})^{\top},$$
where
\begin{align*}
\mathfrak{U}_{i}(t):&=\big(\big. (\mE[x_{1}(t)])^{\top}, (\mE[y_{1}(t)])^{\top}, (\mE[x_{2}(t)])^{\top}, (\mE[y_{2}(t)])^{\top}, \\
 &\quad\quad\quad\quad\quad(\mE^{\sF^{0}_{t}}[x_{1}(t)])^{\top}, (\mE^{\sF^{0}_{t}}[y_{1}(t)])^{\top}, (\mE^{\sF^{0}_{t}}[x_{2}(t)])^{\top}, (\mE^{\sF^{0}_{t}}[y_{2}(t)])^{\top},\\
&\quad\quad \quad(x_{1}(t-\theta))^{\top}, (\mE^{\sF_{t}}[y_{1}(t+\theta)])^{\top}, (x_{2}(t-\theta))^{\top}, (\mE^{\sF_{t}}[y_{2}(t+\theta)])^{\top}, \\
  &\quad\quad \quad\quad
  x_{i}(t)^{\top}, y_{i}(t)^{\top}, z^{0}_{i}(t)^{\top}, z^{1}_{i}(t)^{\top}\big)\big.^{\top}, i=1,2.
\end{align*}
We consider the following MF-AFBSDDEswCN:
\begin{align}\label{25}
\begin{cases}
&\dif \bar{x}_{1}(t)=[b_{1}^{\alpha_{0}}(t,\mathds{V}_{1}(t))+\bar{\psi}_{1}(t)]\dif t
+\sum^{d}_{j=1}[\sigma_{1j}^{0,\alpha_{0}}(t,\mathds{V}_{1}(t))+\bar{\beta}^{0}_{1j}(t)]\dif B^{0}_{j}(t),\\
&\quad\quad\quad\quad\quad\quad\quad\quad\quad\quad\quad\quad+\sum^{d}_{j=1}[\sigma_{1j}^{1,\alpha_{0}}(t,\mathds{V}_{1}(t))+\bar{\beta}^{1}_{1j}(t)]\dif B^{1}_{j}(t),
 t\in [0,T],\\
&\dif \bar{y}_{1}(t) =\{[f_{1}^{\alpha_{0}}(t,\mathds{V}_{1}(t))+\bar{\phi}_{1}(t)]\dif t+\sum^{d}_{j=1} \bar{z}^{0}_{1j}(t)\dif B^{0}_{j}(t)+\sum^{d}_{j=1} \bar{z}^{1}_{1j}(t)\dif B^{1}_{j}(t),\\
& \bar{x}_{1}(0)=\Psi_{1}^{\alpha_{0}}(\bar{y}_{1}(0),\bar{y}_{2}(0))+ \bar{\xi}_{1}, \\ &\bar{x}_{1}(t)=\ddot{\Psi}_{1}(t,\mE[\bar{y}_{1}(t+\theta)],\mE[\bar{y}_{2}(t+\theta)])+\ddot{\bar{\xi}}_{1}(t), t\in [-\theta, 0), \\
 & \bar{y}_{1}(T)=\Phi_{1}^{\alpha}(\bar{x}_{1}(T),\bar{x}_{2}(T) )+\bar{\zeta}_{1}, \bar{y}_{1}(t)=\eta_{1}(t), t\in (T, T+\theta],
\end{cases}
\end{align}

 \begin{align}\label{26}
\begin{cases}
&\dif \bar{x}_{2}(t)=[b_{2}^{\alpha_{0}}(t,\mathds{V}_{2}(t))+\bar{\psi}_{2}(t)]\dif t
 +\sum^{d}_{j=1}[\sigma_{2j}^{0,\alpha_{0}}(t,\mathds{V}_{2}(t))+\bar{\beta}^{0}_{2j}(t)]\dif B^{0}_{j}(t),\\
&\quad\quad\quad\quad\quad\quad\quad\quad\quad\quad\quad\quad+\sum^{d}_{j=1}[\sigma_{2j}^{1,\alpha_{0}}(t,\mathds{V}_{2}(t))+\bar{\beta}^{1}_{2j}(t)]\dif B^{1}_{j}(t),
 t\in [0,T],\\
&\dif \bar{y}_{2 }(t) =[f_{2}^{\alpha_{0}}(t,\mathds{V}_{2}(t))+\bar{\phi}_{2}(t)]\dif t+\sum^{d}_{j=1} \bar{z}^{0}_{2j}(t)\dif B^{0}_{j}(t)+\sum^{d}_{j=1} \bar{z}^{1}_{2j}(t)\dif B^{1}_{j}(t),\\
& \bar{x}_{2 }(0)=\Psi_{2}^{\alpha_{0}}(\bar{y}_{1}(0),\bar{y}_{2}(0))+ \bar{\xi}_{2}, \\ &\bar{x}_{2}(t)=\ddot{\Psi}_{2}(t,\mE[\bar{y}_{1}(t+\theta)],\mE[\bar{y}_{2}(t+\theta)])+\ddot{\bar{\xi}}_{2}(t), t\in [-\theta, 0),\\
  & \bar{y}_{2}(T)=\Phi_{2}^{\alpha_{0}}(\bar{x}_{1}(T),\bar{x}_{2}(T) )+\bar{\zeta}_{2}, \bar{y}_{2}(t)=\eta_{2}(t), t\in (T, T+\theta],
\end{cases}
\end{align}
where
\begin{align*}
&\bar{\xi}_{i}:=\xi_{i}+\delta[\Psi_{i}(y_{1}(0),y_{2}(0))-\Psi^{0}_{i}(y_{1}(0),y_{2}(0))],\\
&\ddot{\bar{\xi}}_{i}(t):=\ddot{\xi}_{i}(t)+\delta[\ddot{\Psi}_{i}(t,\mE[y_{1}(t+\theta)],\mE[y_{2}(t+\theta)])-\ddot{\Psi}^{0}_{i}(t,\mE[y_{1}(t+\theta)],\mE[y_{2}(t+\theta)])],\\
&\bar{\zeta}_{i}:=\zeta_{i}+\delta[\Phi_{i}(x_{1}(T),x_{2}(T))-\Phi^{0}_{i}(x_{1}(T),x_{2}(T))],\\
&\bar{\gamma}_{i}(t):=\gamma_{i}(t)+\delta[\aleph_{i}(t,\mathfrak{U}_{i}(t))-\aleph^{0}_{i}(t,\mathfrak{U}_{i}(t))].\\
\end{align*}
\begin{align*}
\mathds{V}_{1}(\cdot)&:=\big(\big. (\mE[\bar{x}_{1}(t)])^{\top}, (\mE[\bar{y}_{1}(t)])^{\top}, (\mE[\bar{x}_{2}(t)])^{\top}, (\mE[\bar{y}_{2}(t)])^{\top}, \\
 &\quad\quad\quad\quad\quad(\mE^{\sF^{0}_{t}}[\bar{x}_{1}(t)])^{\top}, (\mE^{\sF^{0}_{t}}[\bar{y}_{1}(t)])^{\top}, (\mE^{\sF^{0}_{t}}[\bar{x}_{2}(t)])^{\top}, (\mE^{\sF^{0}_{t}}[\bar{y}_{2}(t)])^{\top},\\
&\quad\quad \quad(\bar{x}_{1}(t-\theta))^{\top}, (\mE^{\sF_{t}}[\bar{y}_{1}(t+\theta)])^{\top}, (\bar{x}_{2}(t-\theta))^{\top}, (\mE^{\sF_{t}}[\bar{y}_{2}(t+\theta)])^{\top}, \\
  &\quad\quad \quad\quad
  \bar{x}_{1}(t)^{\top}, \bar{y}_{1}(t)^{\top}, \bar{z}^{0}_{1}(t)^{\top}, \bar{z}^{1}_{1}(t)^{\top}\big)\big.^{\top},
\end{align*}
\begin{align*}
\mathds{V}_{2}(\cdot)&:=\big(\big. (\mE[\bar{x}_{1}(t)])^{\top}, (\mE[\bar{y}_{1}(t)])^{\top}, (\mE[\bar{x}_{2}(t)])^{\top}, (\mE[\bar{y}_{2}(t)])^{\top}, \\
 &\quad\quad\quad\quad\quad(\mE^{\sF^{0}_{t}}[\bar{x}_{1}(t)])^{\top}, (\mE^{\sF^{0}_{t}}[\bar{y}_{1}(t)])^{\top}, (\mE^{\sF^{0}_{t}}[\bar{x}_{2}(t)])^{\top}, (\mE^{\sF^{0}_{t}}[\bar{y}_{2}(t)])^{\top},\\
&\quad\quad \quad(\bar{x}_{1}(t-\theta))^{\top}, (\mE^{\sF_{t}}[\bar{y}_{1}(t+\theta)])^{\top}, (\bar{x}_{2}(t-\theta))^{\top}, (\mE^{\sF_{t}}[\bar{y}_{2}(t+\theta)])^{\top}, \\
  &\quad\quad \quad\quad
  \bar{x}_{2}(t)^{\top}, \bar{y}_{2}(t)^{\top}, \bar{z}^{0}_{2}(t)^{\top}, \bar{z}^{1}_{2}(t)^{\top}\big)\big.^{\top}.
\end{align*}
One can easily check that\\ $(\bar{\xi}_{i}, \ddot{\bar{\xi}}_{i}(\cdot),\bar{\zeta}_{i}, \bar{\gamma}_{i}(\cdot)) \in \mR^{n}\times L^{2}(-\theta,0;\mR^{n})\times L^{2}_{\sF_{T}}(\mR^{n})\times \sM_{\mF}(\mR^{n+n+nd+nd}), i=1,2.$  From the assumptions in lemma, it holds that Eqs.\eqref{25} and \eqref{26} admit a unique solution $\mathds{V}(\cdot):=(\mathds{V}_{1}(\cdot)^{\top},\mathds{V}_{2}(\cdot)^{\top})^{\top}.$
Thus, given the arbitrariness of $\mathfrak{U}(\cdot):=(\mathfrak{U}_{1}(\cdot)^{\top},\mathfrak{U}_{2}(\cdot)^{\top})^{\top}$, we proceeded to define a mapping $G:$ $$\mM_{\mF}(\mR^{n+n+nd+nd})\times \mM_{\mF}(\mR^{n+n+nd+nd}) \rightarrow \mM_{\mF}(\mR^{n+n+nd+nd})\times \mM_{\mF}(\mR^{n+n+nd+nd}), $$
$$\mathfrak{U}(\cdot) \rightarrow \mathds{V}(\cdot).  $$  If we can prove that $G$ is a   contractive mapping  when $\delta$ is
small enough, we can easily  get the result in the lemma.  For given $\mathfrak{U}(\cdot), \tilde{\mathfrak{U}}(\cdot) \in \mM_{\mF}(\mR^{n+n+nd+nd})\times \mM_{\mF}(\mR^{n+n+nd+nd}), $  set $\mathds{V}(\cdot):=G(\mathfrak{U}(\cdot)), \tilde{\mathds{V}}(\cdot):=G(\tilde{\mathfrak{U}}(\cdot)), \hat{l}:=l-\tilde{l}, l:=\mathfrak{U}, \mathds{V}.  $
From estimate \eqref{12} in Lemma \ref{t1}, we then obtain
\begin{align*}
&\|\hat{\mathds{V}}_{1}(\cdot)\|_{\mM_{\mF}(\mathbb{R}^{n+n+nd+nd})}+\|\hat{\mathds{V}}_{2}(\cdot)\|_{\mM_{\mF}(\mathbb{R}^{n+n+nd+nd})} = \mathbb{E}\bigg[\bigg.\blacklozenge_{\hat{\mathds{V}}_{1}}\bigg]\bigg.+\mathbb{E}\bigg[\bigg.\blacklozenge_{\hat{\mathds{V}}_{2}}\bigg]\bigg.\\
&\leq \delta^2 C\bigg\{\bigg.|\check{\Psi}_{1} - \check{\Psi}_{1}^0|^2 +\sup_{-\theta\leq t\leq 0}|\check{\Psi}_{1}(t) - \check{\Psi}_{1}^0(t)|^2 + \mathbb{E}\bigg[\bigg.|\check{\Phi}_{1}|^2 + \bigg(\bigg.\int_0^T |\check{f}_{1}(t)|\dif t\bigg)\bigg.^2 \\
&\quad+\bigg(\bigg.\int_0^T |\check{b}_{1}(t) - \check{b}^0_{1}(t)|dt\bigg)\bigg.^2 + \int_0^T |\check{\sigma}^{0}_{1}(t) - \check{\sigma}^{0,0}_{1}(t)|^2 \dif t + \int_0^T |\check{\sigma}^{1}_{1}(t) - \check{\sigma}^{1,0}_{1}(t)|^2 \dif t\bigg]\bigg.\\
&\quad +|\check{\Psi}_{2} - \check{\Psi}^0_{2}|^2 +\sup_{-\theta\leq t\leq 0}|\check{\Psi}_{2}(t) - \check{\Psi}_{2}^0(t)|^2+ \mathbb{E}\bigg[\bigg.|\check{\Phi}_{2}|^2 + \bigg(\bigg.\int_0^T |\check{f}_{2}(t)|dt\bigg)\bigg.^2 \\
&\quad+\bigg(\bigg.\int_0^T |\check{b}_{2}(t) - \check{b}^0_{2}(t)|\dif t\bigg)\bigg.^2 + \int_0^T |\check{\sigma}^{0}_{2}(t) - \check{\sigma}^{0,0}_{2}(t)|^2 \dif t + \int_0^T |\check{\sigma}^{1}_{2}(t) - \check{\sigma}^{1,0}_{2}(t)|^2 \dif t\bigg]\bigg.\bigg\}\bigg.,
\end{align*}

\begin{align*}
    &\check{\Psi}_{i}:= \Psi_{i}(y_{1}(0), y_{2}(0)) - \Psi_{i}(\tilde{y}_{1}(0), \tilde{y}_{2}(0)), \quad \check{\Psi}_{i}^{0} := \Psi^{0}_{i}(y_{1}(0), y_{2}(0)) - \Psi^{0}_{i}(\tilde{y}_{1}(0), \tilde{y}_{2}(0)) \\
    &\check{\Phi}_{i}:= \Phi_{i}(x_{1}(T), x_{2}(T)) - \Phi_{i}(\tilde{x}_{1}(T), \tilde{x}_{2}(T)),\\
    & \check{\ddot{\Psi}}_{i}^{0}(t) := \Psi^{0}_{i}(t,
    \mE[y_{1}(t+\theta)],  \mE[y_{2}(t+\theta)]) - \Psi^{0}_{i}(t,
    \mE[\tilde{y}_{1}(t+\theta)],  \mE[\tilde{y}_{2}(t+\theta)]) \\
    & \check{\ddot{\Psi}}_{i}(t) := \Psi_{i}(t,
    \mE[y_{1}(t+\theta)],  \mE[y_{2}(t+\theta)]) - \Psi_{i}(t,
    \mE[\tilde{y}_{1}(t+\theta)],  \mE[\tilde{y}_{2}(t+\theta)]) \\
    &\check{l}_{i}(t):= l_{i}(t, \mathfrak{U}_{i}(t)) - l_{i}(t, \tilde{\mathfrak{U}}_{i}(t)) \quad \text{with} \quad l_{i} = f_{i}, b_{i}, b^0_{i}, \sigma^{0,0}_{i}, \sigma^{1,0}_{i}, \sigma^{0}_{i}, \sigma^{1}_{i}, i=1,2.
\end{align*}
Basing on the Lipschitz condition, it holds that

$$\| \hat{\mathds{V}}(\cdot) \|_{\mM_{\mathbb{F}}(\mathbb{R}^{n+n+nd+nd})\times \mM_{\mathbb{F}}(\mathbb{R}^{n+n+nd+nd})}^2 \leq \delta^2 C \| \hat{\mathfrak{U}}(\cdot) \|_{\mM_{\mathbb{F}}(\mathbb{R}^{n+n+nd+nd})\times \mM_{\mathbb{F}}(\mathbb{R}^{n+n+nd+nd})}^2.$$
Choosing $\delta_{0}:=\frac{1}{2\sqrt{C}},$ then for any $\delta\in (0,\delta_{0}\wedge (1-\alpha_{0})],$ the above inequality shows that $G$ is contractive.
It is evident that the unique fixed point corresponds precisely to the unique solution of System $(\pi_{1})$  when $\alpha=\alpha_{0}+\delta$ and $\xi_{i},\ddot{\xi}_{i}(\cdot), \zeta_{i},\gamma_{i}(\cdot), i=1,2.$ The proof is complete.

\end{proof}

The well-posed of System $(\pi)$ be listed as follows.
\bt\label{tt}
Provided that Assumption 1 holds for the coefficients $(\Psi_{i},\ddot{\Psi}_{i},\Phi_{i},\aleph_{i}), i=1,2 ,$ System $(\pi)$ admits a unique solution.
\et
\begin{proof}
As we previously observed, System $(\pi_{1})$ possesses a unique solution when $\alpha_{0}=0.$ Through the repeated application of  Lemma \ref{t4}, this unique solvability is subsequently extended from $\alpha = 0$ to $\alpha > 0.$ Considering that the fixed step size is $\delta_{0}>0$ , a limited quantity of these extensions suffices to confirm the unique solvability of  System $(\pi)$ for $\alpha=1$,  thus concluding the proof.

\end{proof}

\section{Application to stochastic linear-convex problems}
Here, we examine $Problem (LC),$ a problem formulated in Section 1 based on the solvability of extended  MF-FBSDEswCN.

\subsection{Convex Criterion Functional}
As observed under Assumption 3, Lemma \ref{A1} indicates that the mappings  $\nabla g_{1i}(\cdot): \mathbb{R}^{m} \rightarrow \mathbb{R}^{m}, $ $\nabla\tilde{g}_{3i}(t,\cdot): \mathbb{R}^{n} \rightarrow \mathbb{R}^{n},$ $\nabla \bar{g}_{4i}(t,\cdot): \mathbb{R}^{k} \rightarrow \mathbb{R}^{k}$ and $\nabla g_{4i}(t,\cdot): \mathbb{R}^{k} \rightarrow \mathbb{R}^{k},i=1,2$ are bijective. Let $(\nabla g_{1i})^{-1}(\cdot),$ $(\nabla \tilde{g}_{3i})^{-1}(t, \cdot)$ $(\nabla \bar{g}_{4i})^{-1}(t, \cdot)$ and $(\nabla g_{4i})^{-1}(t, \cdot), i=1,2$  represent the corresponding inverse mappings, respectively. Moreover, under the Assumption 3, we are also aware that
\begin{align}\label{61}
\mathbb{E}&\bigg\{\bigg.g_{21}(x_{1}(T)+x_{2}(T) ) +g_{22}( x_{1}(T)+x_{2}(T) )
  + \int^{T}_{0}g_{31}(t,x_{1}(t) )\dif t +\int^{T}_{0}g_{32}(t,x_{2}(t))\dif t\no\\
 &+ \int^{T}_{0}\tilde{g}_{31}(t-\theta,x_{1}(t-\theta)] )\dif t +\int^{T}_{0}\tilde{g}_{32}(t-\theta,x_{2}(t-\theta))\dif t\no\\
 & + \int^{T}_{0}\ddot{g}_{31}(t,\mE^{\sF^{0}_{t}}[x_{1}(t)] )\dif t +\int^{T}_{0}\ddot{g}_{32}(t,\mE^{\sF^{0}_{t}}[x_{2}(t)])\dif t\no\\
 &+ \int^{T}_{0}\bar{g}_{31}(t,\mE[x_{1}(t)] )\dif t +\int^{T}_{0}\bar{g}_{32}(t,\mE[x_{2}(t)])\dif t\no\\
 &  + \int^{T}_{T-\theta}\langle \tau \tilde{H}_{1}(t+\theta)^{\top}\mE^{\sF_{t}}[\eta_{1}(t+\theta)],  u_{1}(t)\rangle\dif t\no\\
&+ \int^{T}_{T-\theta}\langle \tau \tilde{H}_{2}(t+\theta)^{\top}\mE^{\sF_{t}}[\eta_{2}(t+\theta)],  u_{2}(t)\rangle\dif t\no\\
&+\int^{T}_{T-\theta}\langle \tilde{F}_{2}(t+\theta)^{\top}\mE^{\sF_{t}}[\eta_{2}(t+\theta)]+\tilde{F}_{1}(t+\theta)^{\top}\mE^{\sF_{t}}[\eta_{1}(t+\theta)], x_{1}(t) \rangle\dif t\no\\
 &+\int^{T}_{T-\theta}\langle \tilde{F}_{2}(t+\theta)^{\top}\mE^{\sF_{t}}[\eta_{1}(t+\theta)], x_{2}(t) \rangle\dif t \bigg\}\bigg. < \infty.
\end{align}
 for any $x_{1}(\cdot), x_{2}(\cdot) \in S^{2}_{\mF}(-\theta, T;\mathbb{R}^n), u_{1}(\cdot), u_{2}(\cdot)\in L_{\mF}^2(-\theta, T; \mathbb{R}^k)$. However, the final integral in the criterion functional \eqref{5} maybe diverge to $\infty$. To facilitate later analysis, we  present
\begin{align}\label{62}
\sU_{1} := \left\{u_{1}(\cdot) \in L_{\mF}^2(-\theta,T; \mathbb{R}^k) \ \bigg| \ \mathbb{E} \int_0^T g_{41}(t, u_{1}(t)) dt+\int^{0}_{-\theta}\bar{g}_{41}(t,u_{1} (t) )\dif t  < \infty\right\}.
\end{align}

\begin{align}\label{62+}
\sU_{2} := \left\{u_{2}(\cdot) \in L_{\mF}^2(-\theta,T; \mathbb{R}^k) \ \bigg| \ \mathbb{E} \int_0^T g_{42}(t, u_{2}(t)) dt+\int^{0}_{-\theta}\bar{g}_{42}(t,  u_{2}(t))\dif t< \infty \right\}.
\end{align}
We observe that for $u_{1}(\cdot), u_{2}(\cdot)$ belonging to $L_{\mF}^2(-\theta,T; \mathbb{R}^k)$,$|J(\xi_{1},\xi_{2}, u_{1}(\cdot),u_{2}(\cdot) )| < \infty$  if and only if $(u_{1}, u_{2}) \in \sU_{1}\times \sU_{2} $.  It's obvious that the set  $\sU_{1}, \sU_{2}$ as specified in  \eqref{62}, \eqref{62+} are nonempty. We define $L_{\mF}^2(-\theta,0; \mathbb{R}^n)\times L_{\mF}^2(-\theta,T; \mathbb{R}^n) \times \sU_{1}\times \sU_{2}$ as the admissible control set. When $(\nu_{1}, \nu_{2},u_{1}, u_{2}) \in L_{\mF}^2(-\theta,0; \mathbb{R}^n)\times L_{\mF}^2(-\theta,T; \mathbb{R}^n) \times \sU_{1}\times \sU_{2}$,  we refer to  them  as quadruple admissible controls. Additionally, $x_{i}(\cdot) \equiv x_{i}(\cdot; \nu_{i}, u_{i}),i=1,2$ and $(\nu_{1}, u_{1}, x_{1},\nu_{2}, u_{2}, x_{2})$ are called the corresponding admissible state and an admissible sextet, respectively.

In contrast to the majority of stochastic optimal control problems, Problem (LC) presented in this article features two initial processes  $\nu_{1},
\nu_{2}$ and two process controls $u_{1},u_{2}$.   Nevertheless, it remains a form of Bolza problem.

\subsection{Stochastic Hamiltonian system}The optimal control quartet of Problem (LC) will be characterized using an extended MF-FBSDEswCN, also referred to as a stochastic Hamiltonian system in control theory.

\begin{equation}\label{22++}
J(\nu_{1},\nu_{2}, u_{1}, u_{2}) - J(\nu^*_{1},\nu^*_{2}, u^*_{1},u^*_{2}) = \sum_{i=1}^{7} \Gamma_i
\end{equation}

where
\begin{align*}
\Gamma_1 &:= g_{11}(\xi_{1}+\tau\xi_{2} ) - g_{11}(\xi^{*}_{1}+\tau\xi^{*}_{2})+g_{12}(\xi_{2}+\tau\xi_{1} ) - g_{12}(\xi^{*}_{2}+\tau\xi^{*}_{1}), \\
\Gamma_2 &:= \sum^{2}_{i=1}\{\mathbb{E} \left[ g_{2i}(x_{1}(T)+ x_{2}(T)) - g_{2i}(x^*_{1}(T)+x^*_{2}(T) ) \right]\}, \\
\Gamma_3 &:= \sum^{2}_{i=1}\bigg\{\bigg.\mathbb{E} \int_{0}^{T} \left[ g_{3i}(x_{i}(t) ) - g_{3i}(x^*_{i}(t) ) \right] dt\bigg\}\bigg. \\
&\quad+\sum^{2}_{i=1}\bigg\{\bigg.\mathbb{E} \int_{0}^{T} \left[ \tilde{g}_{3i}(t-\theta,x_{i}(t-\theta) ) - \tilde{g}_{3i}(t-\theta,x^*_{i}(t-\theta) ) \right] dt\bigg\}\bigg.\\
&\quad+\sum^{2}_{i=1}\bigg\{\bigg.\mathbb{E} \int_{0}^{T} \left[ \bar{g}_{3i}(t,\mE[x_{i}(t)] ) - \bar{g}_{3i}(t,\mE[x^*_{i}(t)] )\right] dt\\
 &\quad+\sum^{2}_{i=1}\bigg\{\bigg.\mathbb{E} \int_{0}^{T} \left[ \ddot{g}_{3i}(t,\mE^{\sF_{t}^{0}}[x_{i}(t)] ) - \ddot{g}_{3i}(t,\mE^{\sF_{t}^{0}}[x^*_{i}(t)] )\right] dt\bigg\}\bigg.,\\
\Gamma_4 &:=  \sum^{2}_{i=1}\bigg\{\bigg.\mathbb{E} \int_{0}^{T} \left[ g_{4i}(t, u_{i}(t)) - g_{4i}(t, u^*_{i}(t)) \right] dt\bigg\}\bigg.,\\
\Gamma_5&:=\sum^{2}_{i=1}\bigg\{\bigg.\mathbb{E}\int^{T}_{T-\theta}\langle \tau \tilde{H}_{i}(t+\theta)^{\top}\mE^{\sF_{t}}[\eta_{i}(t+\theta)],  u_{i}(t)- u^{*}_{i}(t)\rangle\dif t\bigg\}\bigg.,\\
\Gamma_6 &:=  \sum^{2}_{i=1}\bigg\{\bigg.\mathbb{E} \int_{-\theta}^{0} \left[ \bar{g}_{4i}(t, u_{i}(t)) - \bar{g}_{4i}(t, u^*_{i}(t)) \right] dt\bigg\}\bigg.,\\
\Gamma_7 &:=\mathbb{E}\int^{T}_{T-\theta}\langle \tilde{F}_{2}(t+\theta)^{\top}\mE^{\sF_{t}}[\eta_{2}(t+\theta)]+\tilde{F}_{1}(t+\theta)^{\top}\mE^{\sF_{t}}[\eta_{1}(t+\theta)], x_{1}(t)-x^{*}_{1}(t) \rangle\dif t\\
&\quad+\mathbb{E}\int^{T}_{T-\theta}\langle \tilde{F}_{2}(t+\theta)^{\top}\mE^{\sF_{t}}[\eta_{1}(t+\theta)], x_{2}(t)-x^{*}_{2}(t) \rangle\dif t.
\end{align*}
Due to the convexity of $g_{2i},$  $g_{3i},$  $\bar{g}_{3i},$ $\ddot{g}_{3i},$ $\tilde{g}_{3i},$ and the uniform convexity of $g_{1i},$ $\bar{g}_{4i}$, and $g_{4i}$, Lemma \ref{100} (2) in the Appendix works to yield

\begin{align}\label{27}
&\Gamma_1 \geq \langle \nabla g_{11}(\xi^*_{1}+\tau\xi^*_{2}), \xi_{1} - \xi^*_{1} \rangle + \langle \tau\nabla g_{11}(\xi^*_{1}+\tau\xi^*_{2}), \xi_{2} - \xi^*_{2} \rangle\no\\
&\quad\quad+\langle \nabla g_{12}(\tau\xi^*_{1}+\xi^*_{2}), \xi_{2} - \xi^*_{2} \rangle + \langle \tau\nabla g_{12}(\tau\xi^*_{1}+\xi^*_{2}), \xi_{1} - \xi^*_{1} \rangle\no\\
&\quad\quad+\frac{\delta}{2} |\xi_{1}+\tau\xi_{2} - \xi^*_{1}-\tau\xi^*_{2}|^2+\frac{\delta}{2} |\tau\xi_{1}+\xi_{2} - \tau\xi^*_{1}-\xi^*_{2}|^2,\no\\
&\Gamma_2\geq \sum^{2}_{i=1}\bigg\{\bigg.\mathbb{E} [ \langle \nabla g_{2i}(x^*_{1}(T)+x^*_{2}(T) ), (x_{1}(T)+x_{2}(T)) - (x^*_{1} (T)+(x^*_{2} (T)\rangle]\bigg\}\bigg.,\no \\
 &\Gamma_3
 \geq \sum^{2}_{i=1}\mathbb{E} \int_{0}^{T} \left[ \langle \nabla g_{3i}( t,x^*_{i}(t)), x_{i}(t) - x^*_{i}(t) \rangle \right] dt \no\\
&\quad\quad+\sum^{2}_{i=1}\mathbb{E} \int_{0}^{T} \left[ \langle \nabla \ddot{g}_{3i}( t,\mE^{\sF_{t}^{0}}[x^*_{i}(t)]), \mE^{\sF_{t}^{0}}[x_{i}(t) - x^*_{i}(t)] \rangle \right] dt\no\\
&\quad\quad+\sum^{2}_{i=1}\mathbb{E} \int_{0}^{T} \left[ \langle \nabla \tilde{g}_{3i}( t-\theta, x^*_{i}(t-\theta)), x_{i}(t-\theta) - x^*_{i}(t-\theta)] \rangle \right] dt\no\\
&\quad\quad+\sum^{2}_{i=1}\mathbb{E} \int_{0}^{T} \left[ \langle \nabla \bar{g}_{3i}( t,\mE[x^*_{i}(t)]), \mE[x_{i}(t) - x^*_{i}(t)] \rangle \right] dt\no\\
&\Gamma_4\geq \sum^{2}_{i=1}\mathbb{E} \int_{0}^{T} \left[ \langle \nabla g_{4i}(t, u^*_{i}(t)), u_{i}(t) - u^*_{i}(t) \rangle + \frac{\delta}{2} |u_{i}(t) - u^*_{i}(t)|^2 \right] dt\no\\
&\Gamma_6\geq \sum^{2}_{i=1}\mathbb{E} \int_{-\theta}^{0} \left[ \langle \nabla \bar{g}_{4i}(t, u^*_{i}(t)), u_{i}(t) - u^*_{i}(t) \rangle + \frac{\delta}{2} |u_{i}(t) - u^*_{i}(t)|^2 \right] dt.
\end{align}
Certainly, \eqref{27}  rephrases \eqref{22++} as an inequality. To render this inequality more manageable, we
utilize a duality perspective and introduce a BSDE. This BSDE, specified on  $[0, T],$ adopts the
subsequent structure:
\begin{align}\label{28}
\begin{cases}
&\dif y_{1}(t)=-[\nabla g_{31}(t,x^{*}_{1}(t))+\nabla \tilde{g}_{31}(t, x^{*}_{1}(t))1_{0\leq t \leq T-\theta}+\nabla \bar{g}_{31}(t,\mE[x^{*}_{1}(t)])\\
&\quad\quad\quad\quad\quad\quad+\nabla \ddot{g}_{31}(t,\mE^{\sF_{t}^{0}}[x^{*}_{1}(t)])+F_{1}(t)^{\top}y_{1}(t)+K_{1}(t)^{\top}z^{0}_{1}(t)+\bar{K}_{1}(t)^{\top}z^{1}_{1}(t)\\
&\quad\quad\quad\quad\quad\quad
+\bar{F}_{1}(t)^{\top}\mE[y_{1}(t)]+\bar{F}_{2}(t)^{\top}\mE[y_{2}(t)]+\ddot{F}_{1}(t)^{\top}\mE^{\sF_{t}^{0}}[y_{1}(t)]+\ddot{F}_{2}(t)^{\top}\mE^{\sF_{t}^{0}}[y_{2}(t)]\\
&\quad\quad\quad\quad\quad\quad\quad\quad\quad\quad\quad\quad+\tilde{F}_{1}(t+\theta)^{\top}\mE^{\sF_{t}}[y_{1}(t+\theta)]+\tilde{F}_{2}(t+\theta)^{\top}\mE^{\sF_{t}}[y_{2}(t+\theta)] ]\dif t\\
&\quad\quad\quad\quad\quad\quad\quad\quad\quad\quad\quad\quad\quad+\sum^{d}_{j=1}z^{0}_{1j}(t)\dif B^{0}_{j}(t)+\sum^{d}_{j=1}z^{1}_{1j}(t)\dif B^{1}_{j}(t), t\in [0,T],\\
& y_{1}(T)=\nabla g_{21}(x^{*}_{1}(T)+x^{*}_{2}(T))+\nabla g_{22}(x^{*}_{1}(T)+x^{*}_{2}(T)),
\end{cases}
\end{align}
and

\begin{align}\label{29}
\begin{cases}
&\dif y_{2}(t)=-[\nabla g_{32}(t,x^{*}_{2}(t))+\nabla \tilde{g}_{32}(t,x^{*}_{2}(t))1_{0\leq t \leq T-\theta}+\nabla \bar{g}_{32}(t,\mE[x^{*}_{2}(t)])\\
&\quad\quad\quad\quad\quad\quad+\nabla \ddot{g}_{32}(t,\mE^{\sF_{t}^{0}}[x^{*}_{2}(t)])+F_{2}(t)^{\top}y_{2}(t)
+K_{2}(t)^{\top}z^{0}_{2}(t)+\bar{K}_{2}(t)^{\top}z^{1}_{2}(t)\\
&\quad\quad\quad\quad\quad\quad+\bar{F}_{2}(t)^{\top}\mE[y_{1}(t)]+\ddot{F}_{2}(t)^{\top}\mE^{\sF_{t}^{0}}[y_{1}(t)]+\tilde{F}_{2}(t+\theta)^{\top}\mE^{\sF_{t}}[y_{1}(t+\theta)] ]\dif t\\
&\quad\quad\quad\quad\quad\quad\quad\quad\quad\quad\quad\quad+\sum^{d}_{j=1}z^{0}_{2j}(t)\dif B^{0}_{j}(t)+\sum^{d}_{j=1}z^{1}_{2j}(t)\dif B^{1}_{j}(t), t\in [0,T],\\
&  y_{2}(T)=\nabla g_{21}(x^{*}_{1}(T)+x^{*}_{2}(T))+\nabla g_{22}(x^{*}_{1}(T)+x^{*}_{2}(T)).
\end{cases}
\end{align}
Eq.\eqref{28} has a unique solution $(y_{1},z_{1}^{0}, ,z_{1}^{1})\in L^{2}_{\mF}(0,T+\theta;\mR^{n})\times L^{2}_{\mF}(0,T;\mR^{nd})\times L^{2}_{\mF}(0,T;\mR^{nd})$ and Eq.\eqref{29} has a unique solution $(y_{2},z_{2}^{0},,z_{2}^{1})\in L^{2}_{\mF}(0,T+\theta; \mR^{n})\times L^{2}_{\mF}(0,T;\mR^{nd})\times L^{2}_{\mF}(0,T;\mR^{nd}).$
Applying It\^{o}'s formula to $\langle y_{1}(t), x_{1}(t)-x^{*}_{1}(t) \rangle, $   we have
\begin{align}\label{30}
&\mE\bigg\{\bigg. \langle \nabla g_{21}(x^{*}_{1}(T)+x^{*}_{2}(T)),   x_{1}(T)-x^{*}_{1}(T)\rangle+\langle \nabla g_{22}(x^{*}_{1}(T)+x^{*}_{2}(T)),   x_{1}(T)-x^{*}_{1}(T)\rangle  \no\\
  &\quad +\int^{T}_{0}\langle \nabla g_{31}(x^{*}_{1}(t)),   x_{1}(t)-x^{*}_{1}(t)\rangle \dif t+\int^{T}_{0}\langle \nabla \bar{g}_{31}(\mE[x^{*}_{1}(t)]),   x_{1}(t)-x^{*}_{1}(t)\rangle \dif t\no\\
  &\quad +\int^{T}_{0}\langle \nabla \ddot{g}_{31}(\mE^{\sF^{0}_{t}}[x^{*}_{1}(t)]),   x_{1}(t)-x^{*}_{1}(t)\rangle \dif t\no\\
  &\quad+\int^{T}_{0}\langle 1_{0\leq t\leq T-\theta}\nabla \tilde{g}_{31}(x^{*}_{1}(t)),   x_{1}(t)-x^{*}_{1}(t)\rangle \dif t   \bigg\}\bigg.\no\\
&=\langle \sH^{\top}y_{1}(0), \xi_{1}-\xi^{*}_{1}\rangle-\mE\int^{T}_{0}\langle F_{1}(t)^{\top}y_{1}(t), x_{1}(t)-x^{*}_{1}(t) \rangle\dif t\no\\
&\quad-\mE\int^{T}_{0}\langle \bar{F}_{1}(t)^{\top}\mE[y_{1}(t)], x_{1}(t)-x^{*}_{1}(t) \rangle\dif t-\mE\int^{T}_{0}\langle \bar{F}_{2}(t)^{\top}\mE[y_{2}(t)], x_{1}(t)-x^{*}_{1}(t) \rangle\dif t\no\\
&\quad-\mE\int^{T}_{0}\langle \ddot{F}_{1}(t)^{\top}\mE^{\sF_{t}^{0}}[y_{1}(t)], x_{1}(t)-x^{*}_{1}(t) \rangle\dif t-\mE\int^{T}_{0}\langle \ddot{F}_{2}(t)^{\top}\mE^{\sF_{t}^{0}}[y_{2}(t)], x_{1}(t)-x^{*}_{1}(t) \rangle\dif t\no\\
&\quad-\mE\int^{T}_{0}\langle \tilde{F}_{1}(t+\theta)^{\top}\mE^{\sF_{t}}[y_{1}(t+\theta)], x_{1}(t)-x^{*}_{1}(t) \rangle\dif t\no\\
&\quad-\mE\int^{T}_{0}\langle \tilde{F}_{2}(t+\theta)^{\top}\mE^{\sF_{t}}[y_{2}(t+\theta)], x_{1}(t)-x^{*}_{1}(t) \rangle\dif t\no\\
&\quad+\mE\int^{T}_{0}\langle F_{1}(t)^{\top}y_{1}(t), x_{1}(t)-x^{*}_{1}(t) \rangle\dif t+\mE\int^{T}_{0}\langle \bar{F}_{2}(t)^{\top}\mE[y_{1}(t)], x_{2}(t)-x^{*}_{2}(t) \rangle\dif t\no\\
&\quad+\mE\int^{T}_{0}\langle \ddot{F}_{2}(t)^{\top}\mE^{\sF_{t}^{0}}[y_{1}(t)], x_{2}(t)-x^{*}_{2}(t) \rangle\dif t+\mE\int^{T}_{0}\langle \tilde{F}_{2}(t)^{\top}y_{1}(t), x_{2}(t-\theta)-x^{*}_{2}(t-\theta) \rangle\dif t\no\\
&\quad+\mE\int^{T}_{0}\langle \bar{F}_{1}(t)^{\top}\mE[y_{1}(t)], x_{1}(t)-x^{*}_{1}(t) \rangle\dif t
+\mE\int^{T}_{0}\langle \ddot{F}_{1}(t)^{\top}\mE^{\sF_{t}^{0}}[y_{1}(t)], x_{1}(t)-x^{*}_{1}(t) \rangle\dif t\no\\
&\quad+\mE\int^{T}_{0}\langle \tilde{F}_{1}(t)^{\top}y_{1}(t), x_{1}(t-\theta)-x^{*}_{1}(t-\theta) \rangle\dif t\no \\
&\quad+\mE\int^{T}_{0}\langle H_{1}(t)^{\top}y_{1}(t)+\tau \bar{H}_{1}(t)^{\top}\mE[y_{1}(t)]+\tau \ddot{H}_{1}(t)^{\top}\mE^{\sF^{0}_{t}}[y_{1}(t)]\no\\
&\quad\quad\quad\quad\quad\quad\quad\quad\quad\quad\quad\quad\quad\quad\quad\quad+ D_{1}(t)^{\top}z^{0}_{1}(t)+ \bar{D}_{1}(t)^{\top}z^{1}_{1}(t), u_{1}(t)-u^{*}_{1}(t)   \rangle\dif t\no\\
&\quad+\mE\int^{T}_{0}\langle y_{1}(t),   \tau \tilde{H}_{1}(t)[u_{1}(t-\theta)-u^{*}_{1}(t-\theta)]\rangle\dif t.
\end{align}
Similarly, applying It\^{o}'s formula to $\langle y_{2}(t), x_{2}(t)-x^{*}_{2}(t) \rangle$  it yields

\begin{align}\label{31}
&\mE\bigg\{\bigg. \langle \nabla g_{21}(x^{*}_{1}(T)+x^{*}_{2}(T)),   x_{2}(T)-x^{*}_{2}(T)\rangle+\langle \nabla g_{22}(x^{*}_{1}(T)+x^{*}_{2}(T)),   x_{2}(T)-x^{*}_{2}(T)\rangle  \no\\
  &\quad+\int^{T}_{0}\langle \nabla g_{32}(x^{*}_{2}(t)),   x_{2}(t)-x^{*}_{2}(t)\rangle \dif t+\int^{T}_{0}\langle \nabla \bar{g}_{32}(\mE[x^{*}_{2}(t)]),   x_{2}(t)-x^{*}_{2}(t)\rangle \dif t \no\\
  &\quad +\int^{T}_{0}\langle \nabla \ddot{g}_{32}(\mE^{\sF^{0}_{t}}[x^{*}_{2}(t)]),   x_{2}(t)-x^{*}_{2}(t)\rangle \dif t\no\\
  &\quad+\int^{T}_{0}\langle 1_{0\leq t\leq T-\theta}\nabla \tilde{g}_{32}(x^{*}_{2}(t)),   x_{2}(t)-x^{*}_{2}(t)\rangle \dif t \bigg\}\bigg.\no\\
&=\langle \sH^{\top}y_{2}(0), \xi_{2}-\xi^{*}_{2}\rangle-\mE\int^{T}_{0}\langle F_{2}(t)^{\top}y_{2}(t), x_{2}(t)-x^{*}_{2}(t) \rangle\dif t\no\\
&\quad-\mE\int^{T}_{0}\langle \bar{F}_{2}(t)^{\top}\mE[y_{1}(t)], x_{2}(t)-x^{*}_{2}(t) \rangle\dif t \rangle\dif t\no\\
&\quad-\mE\int^{T}_{0}\langle \ddot{F}_{2}(t)^{\top}\mE^{\sF^{0}_{t}}[y_{1}(t)], x_{2}(t)-x^{*}_{2}(t) \rangle\dif t \rangle\dif t\no\\
&\quad-\mE\int^{T}_{0}\langle \tilde{F}_{2}(t+\theta)^{\top}\mE^{\sF_{t}}[y_{1}(t+\theta)], x_{2}(t)-x^{*}_{2}(t) \rangle\dif t \rangle\dif t\no\\
&\quad+\mE\int^{T}_{0}\langle F_{2}(t)^{\top}y_{2}(t), x_{2}(t)-x^{*}_{2}(t) \rangle\dif t+\mE\int^{T}_{0}\langle \bar{F}_{2}(t)^{\top}\mE[y_{2}(t)], x_{1}(t)-x^{*}_{1}(t) \rangle\dif t\no\\
&\quad +\mE\int^{T}_{0}\langle \ddot{F}_{2}(t)^{\top}\mE^{\sF^{0}_{t}}[y_{2}(t)], x_{1}(t)-x^{*}_{1}(t) \rangle\dif t\no\\
&\quad+\mE\int^{T}_{0}\langle \tilde{F}_{2}(t)^{\top}y_{2}(t), x_{1}(t-\theta)-x^{*}_{1}(t-\theta) \rangle\dif t\no\\
&\quad+\mE\int^{T}_{0}\langle H_{2}(t)^{\top}y_{2}(t)+\tau \bar{H}_{2}(t)^{\top}\mE[y_{2}(t)]+\tau \ddot{H}_{2}(t)^{\top}\mE^{\sF^{0}_{t}}[y_{2}(t)]\no\\
&\quad\quad\quad\quad\quad\quad\quad\quad\quad\quad\quad\quad\quad\quad\quad\quad+ D_{2}(t)^{\top}z^{0}_{2}(t)+ \bar{D}_{2}(t)^{\top}z^{1}_{2}(t), u_{2}(t)-u^{*}_{2}(t)   \rangle\dif t\no\\
&\quad+\mE\int^{T}_{0}\langle y_{2}(t),   \tau \tilde{H}_{2}(t)[u_{2}(t-\theta)-u^{*}_{2}(t-\theta)]\rangle\dif t.
\end{align}
Since
\begin{align*}
&\mE\int^{T}_{0}\langle \tilde{F}_{2}(t)^{\top}y_{2}(t), x_{1}(t-\theta)-x^{*}_{1}(t-\theta) \rangle\dif t=\mE\int^{T-\theta}_{-\theta}\langle \tilde{F}_{2}(t+\theta)^{\top}y_{2}(t+\theta), x_{1}(t)-x^{*}_{1}(t) \rangle\dif t\\
&=\mE\int^{T-\theta}_{0}\langle \tilde{F}_{2}(t+\theta)^{\top}y_{2}(t+\theta), x_{1}(t)-x^{*}_{1}(t) \rangle\dif t\\
&\quad \quad\quad\quad \quad\quad\quad \quad\quad\quad \quad\quad+\mE\int^{0}_{-\theta}\langle \tilde{F}_{2}(t+\theta)^{\top}y_{2}(t+\theta), x_{1}(t)-x^{*}_{1}(t) \rangle\dif t\\
&=\mE\int^{T}_{0}\langle \tilde{F}_{2}(t+\theta)^{\top}y_{2}(t+\theta), x_{1}(t)-x^{*}_{1}(t) \rangle\dif t\\
&\quad \quad\quad\quad \quad\quad\quad \quad\quad\quad \quad\quad+\mE\int^{0}_{-\theta}\langle \tilde{F}_{2}(t+\theta)^{\top}y_{2}(t+\theta), x_{1}(t)-x^{*}_{1}(t) \rangle\dif t\\
&\quad \quad\quad\quad \quad\quad\quad \quad\quad\quad \quad\quad-\mE\int^{T}_{T-\theta}\langle \tilde{F}_{2}(t+\theta)^{\top}\eta_{2}(t+\theta), x_{1}(t)-x^{*}_{1}(t) \rangle\dif t\\
&=\mE\int^{T}_{0}\langle \tilde{F}_{2}(t+\theta)^{\top}\mE^{\sF_{t}}[y_{2}(t+\theta)], x_{1}(t)-x^{*}_{1}(t) \rangle\dif t\\
&\quad \quad\quad\quad \quad\quad\quad \quad\quad\quad \quad\quad+\mE\int^{0}_{-\theta}\langle \tilde{F}_{2}(t+\theta)^{\top}\mE[y_{2}(t+\theta)], x_{1}(t)-x^{*}_{1}(t) \rangle\dif t\\
& \quad \quad\quad\quad \quad\quad\quad \quad\quad\quad \quad\quad-\mE\int^{T}_{T-\theta}\langle \tilde{F}_{2}(t+\theta)^{\top}\mE^{\sF_{t}}[\eta_{2}(t+\theta)], x_{1}(t)-x^{*}_{1}(t) \rangle\dif t,\\
&\mE\int^{T}_{0}\langle \tilde{F}_{1}(t)^{\top}y_{1}(t), x_{1}(t-\theta)-x^{*}_{1}(t-\theta) \rangle\dif t\\
&=\mE\int^{T}_{0}\langle \tilde{F}_{1}(t+\theta)^{\top}y_{1}(t+\theta), x_{1}(t)-x^{*}_{1}(t) \rangle\dif t\\
&\quad \quad\quad\quad \quad\quad\quad \quad\quad\quad \quad\quad+\mE\int^{0}_{-\theta}\langle \tilde{F}_{1}(t+\theta)^{\top}y_{1}(t+\theta), x_{1}(t)-x^{*}_{1}(t) \rangle\dif t\\
&\quad \quad\quad\quad \quad\quad\quad \quad\quad\quad \quad\quad-\mE\int^{T}_{T-\theta}\langle \tilde{F}_{1}(t+\theta)^{\top}y_{1}(t+\theta), x_{1}(t)-x^{*}_{1}(t) \rangle\dif t\\
&=\mE\int^{T}_{0}\langle \tilde{F}_{1}(t+\theta)^{\top}\mE^{\sF_{t}}[y_{1}(t+\theta)], x_{1}(t)-x^{*}_{1}(t) \rangle\dif t\\
&\quad \quad\quad\quad \quad\quad\quad \quad\quad\quad \quad\quad+\mE\int^{0}_{-\theta}\langle\tilde{F}_{1}(t+\theta)^{\top}\mE[y_{1}(t+\theta)], x_{1}(t)-x^{*}_{1}(t) \rangle\dif t\\
&\quad \quad\quad\quad \quad\quad\quad \quad\quad\quad \quad\quad-\mE\int^{T}_{T-\theta}\langle \tilde{F}_{1}(t+\theta)^{\top}\mE^{\sF_{t}}[\eta_{1}(t+\theta)], x_{1}(t)-x^{*}_{1}(t) \rangle\dif t,\\
&\mE\int^{T}_{0}\langle \tilde{F}_{2}(t)^{\top}y_{1}(t), x_{2}(t-\theta)-x^{*}_{2}(t-\theta) \rangle\dif t\\
&=\mE\int^{T}_{0}\langle \tilde{F}_{2}(t+\theta)^{\top}y_{1}(t+\theta), x_{2}(t)-x^{*}_{2}(t) \rangle\dif t\\
&\quad \quad\quad\quad \quad\quad\quad \quad\quad\quad \quad\quad+\mE\int^{0}_{-\theta}\langle \tilde{F}_{2}(t+\theta)^{\top}y_{1}(t+\theta), x_{2}(t)-x^{*}_{2}(t) \rangle\dif t\\
&\quad \quad\quad\quad \quad\quad\quad \quad\quad\quad \quad\quad-\mE\int^{T}_{T-\theta}\langle \tilde{F}_{2}(t+\theta)^{\top}y_{1}(t+\theta), x_{2}(t)-x^{*}_{2}(t) \rangle\dif t\\
&=\mE\int^{T}_{0}\langle \tilde{F}_{2}(t+\theta)^{\top}\mE^{\sF_{t}}[y_{1}(t+\theta)], x_{2}(t)-x^{*}_{2}(t) \rangle\dif t\\
&\quad \quad\quad\quad \quad\quad\quad \quad\quad\quad \quad\quad+\mE\int^{0}_{-\theta}\langle \tilde{F}_{2}(t+\theta)^{\top}\mE[y_{1}(t+\theta)], x_{2}(t)-x^{*}_{2}(t) \rangle\dif t\\
&\quad \quad\quad\quad \quad\quad\quad \quad\quad\quad \quad\quad-\mE\int^{T}_{T-\theta}\langle \tilde{F}_{2}(t+\theta)^{\top}\mE^{\sF_{t}}[\eta_{1}(t+\theta)], x_{2}(t)-x^{*}_{2}(t) \rangle\dif t,
\end{align*}
and
\begin{align*}
&\mE\int^{T}_{0}\langle y_{i}(t),   \tau \tilde{H}_{i}(t)[u_{i}(t-\theta)-u^{*}_{i}(t-\theta)]\rangle\dif t
=\mE\int^{T}_{0}\langle \tau \tilde{H}_{i}(t)^{\top}y_{i}(t),  u_{i}(t-\theta)-u^{*}_{i}(t-\theta)\rangle\dif t\\
&=\mE\int^{T-\theta}_{0}\langle \tau \tilde{H}_{i}(t+\theta)^{\top}y_{i}(t+\theta),  u_{i}(t)-u^{*}_{i}(t)\rangle\dif t\\
&\quad\quad\quad\quad\quad\quad\quad\quad\quad\quad\quad+\mE\int^{0}_{-\theta}\langle \tau \tilde{H}_{i}(t+\theta)^{\top}y_{i}(t+\theta),  u_{i}(t)-u^{*}_{i}(t)\rangle\dif t\\
&=\mE\int^{T}_{0}\langle \tau \tilde{H}_{i}(t+\theta)^{\top}\mE^{\sF_{t}}[y_{i}(t+\theta)],  u_{i}(t)-u^{*}_{i}(t)\rangle\dif t\\
&\quad\quad\quad\quad\quad\quad\quad\quad\quad\quad\quad +\mE\int^{0}_{-\theta}\langle \tau \tilde{H}_{i}(t+\theta)^{\top}\mE[y_{i}(t+\theta)],  u_{i}(t)-u^{*}_{i}(t)\rangle\dif t\\
&\quad\quad\quad\quad\quad\quad\quad\quad\quad\quad\quad-\mE\int^{T}_{T-\theta}\langle \tau \tilde{H}_{i}(t+\theta)^{\top}\mE^{\sF_{t}}[y_{i}(t+\theta)],  u_{i}(t)-u^{*}_{i}(t)\rangle\dif t,\\
&\int^{T}_{0}\langle 1_{0\leq t\leq T-\theta}\langle\nabla \tilde{g}_{3i}(x^{*}_{i}(t)),   x_{i}(t)-x^{*}_{i}(t)\rangle \dif t
=\int^{T-\theta}_{0}\langle\nabla \tilde{g}_{3i}(x^{*}_{i}(t)),   x_{i}(t)-x^{*}_{i}(t)\rangle \dif t\\
&=\int^{T-\theta}_{-\theta}\langle\nabla \tilde{g}_{3i}(x^{*}_{i}(t)),   x_{i}(t)-x^{*}_{i}(t)\rangle \dif t
-\int^{0}_{-\theta}\langle\nabla \tilde{g}_{3i}(x^{*}_{i}(t)),   x_{i}(t)-x^{*}_{i}(t)\rangle \dif t\\
&=\int^{T}_{0}\langle\nabla \tilde{g}_{3i}(x^{*}_{i}(t-\theta)),   x_{i}(t-\theta)-x^{*}_{i}(t-\theta)\rangle \dif t
 -\int^{0}_{-\theta}\langle\nabla \tilde{g}_{3i}(x^{*}_{i}(t)),   x_{i}(t)-x^{*}_{i}(t)\rangle \dif t,i=1,2,
\end{align*}
using the above calculations and summing up equations \eqref{30}-\eqref{31}  yield
\begin{align}\label{34}
&\mE\bigg\{\bigg. \langle\nabla g_{21}(x^{*}_{1}(T)+x^{*}_{2}(T)),   (x_{1}(T)+x_{2}(T))-(x^{*}_{1}(T)+x^{*}_{2}(T))\rangle\no\\
& \quad\quad +\langle\nabla g_{22}(x^{*}_{1}(T)+x^{*}_{2}(T)),     (x_{1}(T)+x_{2}(T))-(x^{*}_{1}(T)+x^{*}_{2}(T))\rangle\no\\
&\quad\quad+\int^{T}_{0}\langle \nabla g_{31}(x^{*}_{1}(t)),   x_{1}(t)-x^{*}_{1}(t)\rangle \dif t
+\int^{T}_{0}\langle \nabla g_{32}(x^{*}_{2}(t)),   x_{2}(t)-x^{*}_{2}(t)\rangle \dif t\no\\
 &\quad+\int^{T}_{0}\langle \nabla \bar{g}_{31}(\mE[x^{*}_{1}(t)]),   x_{1}(t)-x^{*}_{1}(t)\rangle \dif t+\int^{T}_{0}\langle \nabla \bar{g}_{32}(\mE[x^{*}_{2}(t)]),   x_{2}(t)-x^{*}_{2}(t)\rangle \dif t\no\\
&\quad +\int^{T}_{0}\langle \nabla \ddot{g}_{31}(\mE^{\sF^{0}_{t}}[x^{*}_{1}(t)]),   x_{1}(t)-x^{*}_{1}(t)\rangle \dif t+\int^{T}_{0}\langle \nabla \ddot{g}_{32}(\mE^{\sF^{0}_{t}}[x^{*}_{2}(t)]),   x_{2}(t)-x^{*}_{2}(t)\rangle \dif t\no\\
 &\quad+\int^{T}_{0}\langle \nabla \tilde{g}_{31}(x^{*}_{1}(t-\theta)),   x_{1}(t-\theta)-x^{*}_{1}(t-\theta)\rangle \dif t\no\\
 &\quad+\int^{T}_{0}\langle \nabla \tilde{g}_{32}(x^{*}_{2}(t-\theta)),   x_{2}(t-\theta)-x^{*}_{2}(t-\theta)\rangle \dif t  \bigg\}\bigg.\no\\
  &\quad+\mE\int^{T}_{T-\theta}\langle \tilde{F}_{2}(t+\theta)^{\top}\mE^{\sF_{t}}[\eta_{2}(t+\theta)], x_{1}(t)-x^{*}_{1}(t) \rangle\dif t\no\\
 &\quad+\mE\int^{T}_{T-\theta}\langle \tilde{F}_{1}(t+\theta)^{\top}\mE^{\sF_{t}}[\eta_{1}(t+\theta)], x_{1}(t)-x^{*}_{1}(t) \rangle\dif t\no\\
  &\quad+\mE\int^{T}_{T-\theta}\langle \tilde{F}_{2}(t+\theta)^{\top}\mE^{\sF_{t}}[\eta_{1}(t+\theta)], x_{2}(t)-x^{*}_{2}(t) \rangle\dif t\no\\
  & \quad + \mE\int^{T}_{T-\theta}\langle \tau \tilde{H}_{1}(t+\theta)^{\top}\mE^{\sF_{t}}[\eta_{1}(t+\theta)],  u_{1}(t)-u^{*}_{1}(t)\rangle\dif t\no\\
&\quad+ \mE\int^{T}_{T-\theta}\langle \tau \tilde{H}_{2}(t+\theta)^{\top}\mE^{\sF_{t}}[\eta_{2}(t+\theta)],  u_{2}(t)-u^{*}_{2}(t)\rangle\dif t\no\\
  &=\langle \sH^{\top}y_{1}(0), \xi_{1}-\xi^{*}_{1}\rangle
+\langle \sH^{\top}y_{2}(0), \xi_{2}-\xi^{*}_{2}\rangle\no\\
&+\int^{0}_{-\theta}\langle\nabla \tilde{g}_{31}(x^{*}_{1}(t))+\tilde{F}_{2}(t+\theta)^{\top}\mE[y_{2}(t+\theta)]+\tilde{F}_{1}(t+\theta)^{\top}\mE^{\sF_{t}}[y_{1}(t+\theta)],   x_{1}(t)-x^{*}_{1}(t)\rangle \dif t\no\\
&+\int^{0}_{-\theta}\langle\nabla \tilde{g}_{32}(x^{*}_{2}(t))+\tilde{F}_{2}(t+\theta)^{\top}\mE[y_{1}(t+\theta)],   x_{2}(t)-x^{*}_{2}(t)\rangle \dif t\no\\
&+\int^{0}_{-\theta}\langle\tau \tilde{H}_{1}(t+\theta)^{\top}\mE[y_{1}(t+\theta)],   u_{1}(t)-u^{*}_{1}(t)\rangle \dif t\no\\
&+\int^{0}_{-\theta}\langle\tau \tilde{H}_{2}(t+\theta)^{\top}\mE[y_{2}(t+\theta)],   u_{2}(t)-u^{*}_{2}(t)\rangle \dif t\no\\
&+\mE\int^{T}_{0}\langle H_{1}^{\top}(t)y_{1}(t)+\tau \bar{H}_{1}^{\top}(t)\mE[y_{1}(t)]+\tau \ddot{H}_{1}^{\top}(t)\mE^{\sF^{0}_{t}}[y_{1}(t)]+\tau \tilde{H}_{1}^{\top}(t+\theta)\mE^{\sF_{t}}[y_{1}(t+\theta)]\no\\
&\quad\quad\quad\quad\quad\quad\quad\quad\quad\quad\quad\quad\quad\quad\quad\quad\quad+ D_{1}^{\top}(t)z^{0}_{1}(t)+\bar{D}_{1}^{\top}(t)z^{1}_{1}(t), u_{1}(t)-u^{*}_{1}(t)   \rangle\dif t\no\\
& +\mE\int^{T}_{0}\langle    H_{2}(t)^{\top}y_{2}(t)+ \tau\bar{H}_{2}(t)^{\top}\mE[y_{2}(t)]+\tau \ddot{H}_{2}^{\top}(t)\mE^{\sF^{0}_{t}}[y_{2}(t)]+\tau \tilde{H}_{2}(t+\theta)^{\top}\mE^{\sF_{t}}[y_{2}(t+\theta)]\no\\
&\quad\quad\quad\quad\quad\quad\quad\quad\quad\quad\quad\quad\quad\quad\quad\quad\quad+D_{2}(t)^{\top}z^{0}_{2}(t)+\bar{D}_{2}(t)^{\top}z^{1}_{2}(t),   u_{2}(t)-u^{*}_{2}(t) \rangle\dif t.
\end{align}
Combining \eqref{22++}-\eqref{34} leads to
\begin{align}\label{35}
&J(\nu_{1},\nu_{2}, u_{1}, u_{2}) - J(\nu^*_{1},\nu^*_{2}, u^*_{1},u^*_{2})\no\\
&\geq\langle \nabla g_{11}(\xi_{1}^{*}+\tau\xi_{2}^{*})+\tau \nabla g_{12}(\tau\xi_{1}^{*}+\xi_{2}^{*})+\sH^{\top}y_{1}(0), \xi_{1}-\xi^{*}_{1}\rangle\no\\
&\quad +\langle\tau\nabla g_{11}(\xi_{1}^{*}+\tau\xi_{2}^{*})+\nabla g_{12}(\tau\xi_{1}^{*}+\xi_{2}^{*})+\sH^{\top}y_{2}(0), \xi_{2}-\xi^{*}_{2}\rangle\no\\
&\quad+\int^{0}_{-\theta}\langle\nabla \tilde{g}_{31}(x^{*}_{1}(t))+\tilde{F}_{2}(t+\theta)^{\top}\mE[y_{2}(t+\theta)]\no\\
&\quad \quad\quad\quad\quad \quad\quad+\tilde{F}_{1}(t+\theta)^{\top}\mE[y_{1}(t+\theta)],   x_{1}(t)-x^{*}_{1}(t)\rangle \dif t\no\\
&\quad+\int^{0}_{-\theta}\langle\nabla \tilde{g}_{32}(x^{*}_{2}(t))+\tilde{F}_{2}(t+\theta)^{\top}\mE[y_{1}(t+\theta)],   x_{2}(t)-x^{*}_{2}(t)\rangle \dif t\no\\
&\quad +\mE\int^{T}_{0}\langle \nabla g_{41}(u^{*}_{1}(t))+ H_{1}^{\top}(t)y_{1}(t)+ \tau\bar{H}_{1}(t)^{\top}\mE[y_{1}(t)]+ \tau\ddot{H}_{1}^{\top}(t)\mE^{\sF^{0}_{t}}[y_{1}(t)]\no\\
&\quad \quad\quad\quad\quad + \tau\tilde{H}_{1}(t+\theta)^{\top}\mE^{\sF_{t}}[y_{1}(t+\theta)]+D_{1}(t)^{\top}z^{0}_{1}(t)+\bar{D}_{1}(t)^{\top}z^{0}_{1}(t), u_{1}(t)-u^{*}_{1}(t)   \rangle\dif t\no\\
&\quad +\mE\int^{T}_{0}\langle \nabla g_{42}(u^{*}_{2}(t))+H_{2}(t)^{\top}y_{2}(t)+ \tau\bar{H}_{2}(t)^{\top}\mE[y_{2}(t)]+ \tau\ddot{H}_{2}(t)^{\top}\mE^{\sF^{0}_{t}}[y_{2}(t)]\no\\
& \quad \quad\quad\quad\quad + \tau\tilde{H}_{2}^{\top}(t+\theta)\mE^{\sF_{t}}[y_{2}(t+\theta)]+D^{\top}_{2}(t)z^{0}_{2}(t)+\bar{D}^{\top}_{2}(t)z^{1}_{2}(t),   u_{2}(t)-u^{*}_{2}(t) \rangle\dif t\no\\
&+\int^{0}_{-\theta}\langle \bar{g}_{41}(t,u_{1}(t))+\tau \tilde{H}_{1}(t+\theta)^{\top}\mE[y_{1}(t+\theta)],   u_{1}(t)-u^{*}_{1}(t)\rangle \dif t\no\\
&+\int^{0}_{-\theta}\langle \bar{g}_{42}(t,u_{2}(t))+\tau \tilde{H}_{2}(t+\theta)^{\top}\mE[y_{2}(t+\theta)],   u_{2}(t)-u^{*}_{2}(t)\rangle \dif t\no\\
&\quad+\frac{\delta}{2}\bigg\{\bigg. |\xi_{1}+\tau\xi_{2} - \xi^*_{1}-\tau\xi^*_{2}|^2+ |\tau\xi_{1}+\xi_{2} - \tau\xi^*_{1}-\xi^*_{2}|^2\no\\
& \quad+\mE\int^{T}_{0} |u_{1}(t)-u^{*}_{1}(t)|^{2}\dif t +\mE\int^{T}_{0} |u_{2}(t)-u^{*}_{2}(t)|^{2}\dif t    \bigg\}\bigg..
\end{align}
Denote
\begin{align*}
&\Gamma^{*}_7:=\nabla g_{11}(\xi_{1}^{*}+\tau\xi_{2}^{*})+\tau \nabla g_{12}(\tau\xi_{1}^{*}+\xi_{2}^{*})+\sH^{\top}y_{1}(0),\\
&\Gamma^{*}_8:=\tau\nabla g_{11}(\xi_{1}^{*}+\tau\xi_{2}^{*})+\nabla g_{12}(\tau\xi_{1}^{*}+\xi_{2}^{*})+\sH^{\top}y_{2}(0),\\
&\Gamma^{*}_9:= \nabla g_{41}(u^{*}_{1}(t))+ H_{1}(t)^{\top}y_{1}(t)+ \tau\bar{H}_{1}(t)^{\top}\mE[y_{1}(t)]+ \tau\ddot{H}_{1}(t)^{\top}\mE^{\sF^{0}_{t}}[y_{1}(t)]\\
& \quad \quad \quad\quad\quad \quad \quad\quad\quad \quad \quad+ \tau\tilde{H}_{1}(t+\theta)^{\top}\mE[y_{1}(t+\theta)]+D_{1}(t)^{\top}z^{0}_{1}(t)+\bar{D}_{1}(t)^{\top}z^{1}_{1}(t),\\
&\Gamma^{*}_{10}:=\nabla g_{42}(u^{*}_{2}(t))+H_{2}(t)^{\top}y_{2}(t)+ \tau\bar{H}_{2}^{\top}(t)\mE[y_{2}(t)]+ \tau\ddot{H}_{2}(t)^{\top}\mE^{\sF^{0}_{t}}[y_{2}(t)]\\
& \quad \quad \quad\quad\quad \quad \quad\quad\quad \quad \quad+ \tau\tilde{H}_{2}(t+\theta)^{\top}\mE[y_{2}(t+\theta)]+D_{2}(t)^{\top}z^{0}_{2}(t)+\bar{D}_{2}(t)^{\top}z^{1}_{2}(t), t\in [0,T],\\
\end{align*}
and
\begin{align*}
&\bar{\Gamma}^{*}_7:=\nabla \tilde{g}_{31}(x^{*}_{1}(t))+\tilde{F}_{2}(t+\theta)^{\top}\mE[y_{2}(t+\theta)]+\tilde{F}_{1}(t+\theta)^{\top}\mE[y_{1}(t+\theta)],\\
&\bar{\Gamma}^{*}_8:=\nabla \tilde{g}_{32}(x^{*}_{2}(t))+\tilde{F}_{2}(t+\theta)^{\top}\mE[y_{1}(t+\theta)],\\
&\bar{\Gamma}^{*}_9:=\bar{g}_{41}(t,u_{1}(t))+\tau \tilde{H}_{1}(t+\theta)^{\top}\mE[y_{1}(t+\theta)],\\
&\bar{\Gamma}^{*}_{10}:=\bar{g}_{42}(t,u_{1}(t))+\tau \tilde{H}_{1}(t+\theta)^{\top}\mE[y_{2}(t+\theta)], t\in [-\theta,0).
\end{align*}
This section focuses on identifying optimal control pairs for Problem (LC). Considering inequality \eqref{35}, the optimality of $(\nu_{1},\nu_{2},u_{1},u_{2})$ appears equivalent to the following condition
$$(\Gamma^{*}_7, \Gamma^{*}_8, \Gamma^{*}_9, \Gamma^{*}_{10}) = (0, 0,0,0), (\bar{\Gamma}^{*}_7, \bar{\Gamma}^{*}_8, \bar{\Gamma}^{*}_9, \bar{\Gamma}^{*}_{10}) = (0, 0,0,0).$$
Then,  we have
\begin{align}\label{36}
\begin{cases}
\xi^{*}_{1}=\frac{1}{1-\tau^{2}}\bigg\{\bigg.(\nabla g_{11})^{-1}(-\frac{\sH^{\top}y_{1}(0)-\tau \sH^{\top}y_{2}(0)}{1-\tau^{2}})-\tau(\nabla g_{12})^{-1}(-\frac{\sH^{\top}y_{2}(0)-\tau \sH^{\top}y_{1}(0)}{1-\tau^{2}})\bigg\}\bigg.,\\
\xi^{*}_{2}=\frac{1}{1-\tau^{2}}\bigg\{\bigg.(\nabla g_{12})^{-1}(-\frac{\sH^{\top}y_{2}(0)-\tau \sH^{\top}y_{1}(0)}{1-\tau^{2}})-\tau(\nabla g_{11})^{-1}(-\frac{\sH^{\top}y_{1}(0)-\tau \sH^{\top}y_{2}(0)}{1-\tau^{2}})\bigg\}\bigg.,\\
u^{*}_{1}(t)=(\nabla g_{41})^{-1}(-( H_{1}(t)^{\top}y_{1}(t)+ \tau\bar{H}_{1}(t)^{\top}\mE[y_{1}(t)]+\tau\ddot{H}_{1}(t)^{\top}\mE^{\sF^{0}_{t}}[y_{1}(t)]\\
\quad\quad\quad\quad\quad\quad\quad+ \tau\tilde{H}_{1}(t+\theta)^{\top}\mE^{\sF_{t}}[y_{1}(t+\theta)]+D_{1}(t)^{\top}z^{0}_{1}(t)+\bar{D}_{1}(t)^{\top}z^{1}_{1}(t))),\\
u^{*}_{2}(t)=(\nabla g_{42})^{-1}(-( H_{2}(t)^{\top}y_{2}(t)+ \tau\bar{H}_{2}(t)^{\top}\mE[y_{2}(t)]+\tau\ddot{H}_{1}(t)^{\top}\mE^{\sF^{0}_{t}}[y_{1}(t)]\\
\quad\quad\quad\quad\quad\quad\quad+ \tau\tilde{H}_{2}(t+\theta)^{\top}\mE^{\sF_{t}}[y_{2}(t+\theta)]+D_{2}(t)^{\top}z^{0}_{2}(t)+\bar{D}_{2}(t)^{\top}z^{1}_{2}(t))), t\in [0,T],
\end{cases}
\end{align}
and
\begin{align}\label{36++++-}
\begin{cases}
x^{*}_{1}(t)=(\nabla \tilde{g}_{31})^{-1}(-( \tilde{F}_{2}(t+\theta)^{\top}\mE[y_{2}(t+\theta)]+\tilde{F}_{1}(t+\theta)^{\top}\mE[y_{1}(t+\theta)])),\\
x^{*}_{2}(t)=(\nabla \tilde{g}_{32})^{-1}(-( \tilde{F}_{2}(t+\theta)^{\top}\mE[y_{2}(t+\theta)])),\\
u^{*}_{1}(t)=(\nabla \bar{g}_{41})^{-1}(-( \tau \tilde{H}_{1}(t+\theta)^{\top}\mE[y_{1}(t+\theta)])),\\
u^{*}_{2}(t)=(\nabla \bar{g}_{42})^{-1}(-( \tau \tilde{H}_{2}(t+\theta)^{\top}\mE[y_{2}(t+\theta)])), t\in [-\theta,0),
\end{cases}
\end{align}
This conjecture will be formally demonstrated in the following main result in this section.

At the conclusion of this section, we combine the SDEs for the state $x^{*}_{1}(\cdot), x^{*}_{2}(\cdot)$ [see \eqref{4}-\eqref{4d}], the BSDEs \eqref{28}, \eqref{29} and the expression \eqref{36} into a single system, referred to as a stochastic Hamiltonian system, as follows:
\begin{align}\label{37a}
&\dif x^{*}_{1}(t)=[F_{1}(t)x^{*}_{1}(t)\no\\
 &\quad\quad+H_{1}(t)(\nabla g_{41})^{-1}(-( H^{\top}_{1}(t)y_{1}(t)+ \tau\bar{H}_{1}(t)^{\top}\mE[y_{1}(t)]+ \tau\ddot{H}_{1}(t)^{\top}\mE^{\sF^{0}_{t}}[y_{1}(t)]\no\\
&\quad\quad\quad\quad\quad\quad+ \tau\tilde{H}_{1}(t+\theta)^{\top}\mE^{\sF_{t}}[y_{1}(t+\theta)]+D_{1}(t)^{\top}z^{0}_{1}(t)+\bar{D}_{1}(t)^{\top}z^{1}_{1}(t)))\no\\
&\quad\quad+ \tau\bar{H}_{1}(t)\mE[(\nabla g_{41})^{-1}(-( H^{\top}_{1}(t)y_{1}(t)+ \tau\bar{H}_{1}(t)^{\top}\mE[y_{1}(t)]+ \tau\ddot{H}_{1}(t)^{\top}\mE^{\sF^{0}_{t}}[y_{1}(t)]\no\\
&\quad\quad\quad\quad\quad\quad+ \tau\tilde{H}_{1}(t+\theta)^{\top}\mE^{\sF_{t}}[y_{1}(t+\theta)]+D_{1}(t)^{\top}z^{0}_{1}(t)+\bar{D}_{1}(t)^{\top}z^{1}_{1}(t)))]\no\\
&\quad\quad+ \tau\ddot{H}_{1}(t)\mE^{\sF^{0}_{t}}[(\nabla g_{41})^{-1}(-( H^{\top}_{1}(t)y_{1}(t)+ \tau\bar{H}_{1}(t)^{\top}\mE[y_{1}(t)]+ \tau\ddot{H}_{1}(t)^{\top}\mE^{\sF^{0}_{t}}[y_{1}(t)]\no\\
&\quad\quad\quad\quad\quad\quad+ \tau\tilde{H}_{1}(t+\theta)^{\top}\mE^{\sF_{t}}[y_{1}(t+\theta)]+D_{1}(t)^{\top}z^{0}_{1}(t)+\bar{D}_{1}(t)^{\top}z^{1}_{1}(t)))]\no\\
&\quad\quad+ 1_{\{\theta \leq t \leq T\}}\tau\tilde{H}_{1}(t)(\nabla g_{41})^{-1}(-( H^{\top}_{1}(t-\theta)y_{1}(t-\theta)+ \tau\bar{H}_{1}(t-\theta)^{\top}\mE[y_{1}(t-\theta)]\no\\
&\quad\quad\quad+ \tau\tilde{H}_{1}(t)^{\top}\mE^{\sF_{t-\theta}}[y_{1}(t)]+\tau\ddot{H}_{1}(t-\theta)^{\top}\mE^{\sF^{0}_{t-\theta}}[y_{1}(t-\theta)]\no\\
&\quad\quad\quad+D_{1}(t-\theta)^{\top}z^{0}_{1}(t-\theta)+\bar{D}_{1}(t-\theta)^{\top}z^{1}_{1}(t-\theta)))\no\\
&\quad\quad+ 1_{\{0 \leq t < \theta\}}\tau\tilde{H}_{1}(t)(\nabla \bar{g}_{41})^{-1}(-( \tau \tilde{H}_{1}(t)^{\top}\mE[y_{1}(t)]))\no\\
& \quad\quad +\bar{F}_{2}(t)\mE[ x^{*}_{2}(t)]+\bar{F}_{1}(t)\mE[ x^{*}_{1}(t)]+\ddot{F}_{2}(t)\mE^{
\sF_{t}^{0}}[ x^{*}_{2}(t)]+\ddot{F}_{1}(t)\mE^{
\sF_{t}^{0}}[ x^{*}_{1}(t)]\no\\
&\quad\quad\quad\quad+\tilde{F}_{2}(t) x^{*}_{2}(t-\theta)+\tilde{F}_{1}(t) x^{*}_{1}(t-\theta)+\upsilon_{1}(t)]\dif t\no\\
&\quad\quad +\sum^{d}_{j=1}[K_{1j}(t)x^{*}_{1}(t) \no\\
&\quad\quad\quad\quad +D_{1j}(t)(\nabla g_{41})^{-1}(-( H^{\top}_{1}(t)y_{1}(t)+ \tau\bar{H}_{1}(t)^{\top}\mE[y_{1}(t)]+ \tau\ddot{H}_{1}(t)^{\top}\mE^{\sF^{0}_{t}}[y_{1}(t)]\no\\
&\quad\quad\quad+ \tau\tilde{H}_{1}(t+\theta)^{\top}\mE^{\sF_{t}}[y_{1}(t+\theta)]+D_{1}(t)^{\top}z^{0}_{1}(t)+\bar{D}_{1}(t)^{\top}z^{1}_{1}(t)))+\iota_{1j}(t)]\dif B^{0}_{j}(t), \no\\
&\quad\quad+\sum^{d}_{j=1}[\bar{K}_{1j}(t)x^{*}_{1}(t)\no\\
&\quad\quad\quad\quad +\bar{D}_{1j}(t)(\nabla g_{41})^{-1}(-( H^{\top}_{1}(t)y_{1}(t)+ \tau\bar{H}_{1}(t)^{\top}\mE[y_{1}(t)]+ \tau\ddot{H}_{1}(t)^{\top}\mE^{\sF^{0}_{t}}[y_{1}(t)]\no\\
&\quad\quad\quad+ \tau\tilde{H}_{1}(t+\theta)^{\top}\mE^{\sF_{t}}[y_{1}(t+\theta)]+D_{1}(t)^{\top}z^{0}_{1}(t)+\bar{D}_{1}(t)^{\top}z^{1}_{1}(t)))+\bar{\iota}_{1j}(t)]\dif B^{1}_{j}(t),
\end{align}
and
\begin{align}\label{37b}
&\dif y_{1}(t)=-[\nabla g_{31}(x^{*}_{1}(t))+\nabla \bar{g}_{31}(\mE[x^{*}_{1}(t)])+\nabla \ddot{g}_{31}(\mE^{\sF^{0}_{t}}[x^{*}_{1}(t)])+1_{\{0\leq t\leq T-\theta\}}\nabla \tilde{g}_{31}(x^{*}_{1}(t))\no\\
& \quad \quad +F_{1}(t)^{\top}y_{1}(t)+K_{1}(t)^{\top}z^{0}_{1}(t)+\bar{K}_{1}(t)^{\top}z^{1}_{1}(t)+\bar{F}_{1}(t)^{\top}\mE[y_{1}(t)]+\bar{F}_{2}(t)^{\top}\mE[y_{2}(t)] \no\\
&\quad\quad+\ddot{F}_{1}(t)^{\top}\mE^{\sF^{0}_{t}}[y_{1}(t)]+\ddot{F}_{2}(t)^{\top}\mE^{\sF^{0}_{t}}[y_{2}(t)] \no\\ &\quad\quad+\tilde{F}_{1}(t+\theta)^{\top}\mE^{\sF_{t}}[y_{1}(t+\theta)]+\tilde{F}_{2}(t+\theta)^{\top}\mE^{\sF_{t}}[y_{2}(t+\theta)] ]\dif t\no\\
&\quad\quad+\sum^{d}_{i=1}z^{0}_{1j}(t)\dif B^{0}_{j}(t)+\sum^{d}_{i=1}z^{1}_{1j}(t)\dif B^{1}_{j}(t).
\end{align}

\begin{align}\label{37c}
\dif x^{*}_{2}(t)&=[F_{2}(t)x^{*}_{2}(t)\no\\
&\quad+H_{2}(t)(\nabla g_{42})^{-1}(-( H_{2}(t)^{\top}y_{2}(t)+ \tau\bar{H}_{2}(t)^{\top}\mE[y_{2}(t)]+\tau\ddot{H}_{2}(t)^{\top}\mE^{\sF^{0}_{t}}[y_{2}(t)]\no\\
&\quad+ \tau\tilde{H}_{2}(t+\theta)^{\top}\mE^{\sF_{t}}[y_{2}(t+\theta)]+D_{2}(t)^{\top}z^{0}_{2}(t)+\bar{D}_{2}(t)^{\top}z^{1}_{2}(t)))\no\\
&\quad+ \tau \bar{H}_{2}(t)\mE[(\nabla g_{42})^{-1}(-( H^{\top}_{2}(t)y_{2}(t)+ \tau\bar{H}_{2}(t)^{\top}\mE[y_{2}(t)]+\tau\ddot{H}_{2}(t)^{\top}\mE^{\sF^{0}_{t}}[y_{2}(t)]\no\\
&\quad+ \tau\tilde{H}_{2}(t+\theta)^{\top}\mE^{\sF_{t}}[y_{2}(t+\theta)]+D_{2}(t)^{\top}z^{0}_{2}(t)+\bar{D}_{2}(t)^{\top}z^{1}_{2}(t)))]\no\\
&\quad+ \tau \ddot{H}_{2}(t)\mE[(\nabla g_{42})^{-1}(-( H^{\top}_{2}(t)y_{2}(t)+ \tau\bar{H}_{2}(t)^{\top}\mE[y_{2}(t)]+\tau\ddot{H}_{2}(t)^{\top}\mE^{\sF^{0}_{t}}[y_{2}(t)]\no\\
&\quad+ \tau\tilde{H}_{2}(t+\theta)^{\top}\mE^{\sF_{t}}[y_{2}(t+\theta)]+D_{2}(t)^{\top}z^{0}_{2}(t)+\bar{D}_{2}(t)^{\top}z^{1}_{2}(t)))]\no\\
&\quad+ 1_{\{\theta \leq t \leq T\}}\tau\tilde{H}_{2}(t)(\nabla g_{42})^{-1}(-( H_{2}(t-\theta)^{\top}y_{2}(t-\theta)+ \tau\bar{H}_{2}(t-\theta)^{\top}\mE[y_{2}(t-\theta)]\no\\
&\quad+ \tau\tilde{H}_{2}(t)^{\top}\mE^{\sF_{t-\theta}}[y_{2}(t)]+ \tau\ddot{H}_{2}(t-\theta)^{\top}\mE^{\sF^{0}_{t-\theta}}[y_{2}(t-\theta)]\no\\
&\quad+D_{2}(t-\theta)^{\top}z^{0}_{2}(t-\theta)+\bar{D}_{2}(t-\theta)^{\top}z^{1}_{2}(t-\theta)))\no\\
&\quad+ 1_{\{0 \leq t < \theta\}}\tau\tilde{H}_{2}(t)(\nabla \bar{g}_{42})^{-1}(-( \tau \tilde{H}_{2}(t)^{\top}\mE[y_{2}(t)]))\no\\
&\quad +\bar{F}_{2}(t)\mE[ x^{*}_{1}(t) ]+\ddot{F}_{2}(t)\mE^{
\sF_{t}^{0}}[ x^{*}_{1}(t)] +\tilde{F}_{2}(t) x^{*}_{1}(t-\theta) +\upsilon_{2}(t)]\dif t\no\\
&\quad+\sum^{d}_{j=1}[K_{2j}(t)x^{*}_{1}(t) +D_{2j}(t)(\nabla g_{42})^{-1}(-( H^{\top}_{2}(t)y_{2}(t)+ \tau\bar{H}_{2}(t)^{\top}\mE[y_{2}(t)]\no\\
&\quad+ \tau\ddot{H}_{2}(t)^{\top}\mE^{\sF^{0}_{t}}[y_{2}(t)]+ \tau\tilde{H}_{2}(t+\theta)^{\top}\mE^{\sF_{t}}[y_{2}(t+\theta)]\no\\
&\quad+D_{2}(t)^{\top}z^{0}_{2}(t)+\bar{D}_{2}(t)^{\top}z^{1}_{2}(t)))+\iota_{2j}(t)]\dif B^{0}_{j}(t), \no\\
&\quad+\sum^{d}_{j=1}[\bar{K}_{2j}(t)x^{*}_{2}(t) +\bar{D}_{2j}(t)(\nabla g_{42})^{-1}(-( H^{\top}_{2}(t)y_{2}(t)+ \tau\bar{H}_{2}(t)^{\top}\mE[y_{2}(t)]\no\\
&\quad+ \tau\ddot{H}_{2}(t)^{\top}\mE^{\sF^{0}_{t}}[y_{2}(t)]+ \tau\tilde{H}_{2}(t+\theta)^{\top}\mE^{\sF_{t}}[y_{2}(t+\theta)]\no\\
&\quad+D_{2}(t)^{\top}z^{0}_{2}(t)+\bar{D}_{2}(t)^{\top}z^{1}_{2}(t)))+\bar{\iota}_{2j}(t)]\dif B^{1}_{j}(t),
\end{align}
and
\begin{align}\label{37d}
&\dif y_{2}(t)=-[\nabla g_{32}(x^{*}_{2}(t))+\nabla \bar{g}_{32}(\mE[x^{*}_{2}(t)])+\nabla \ddot{g}_{32}(\mE^{\sF^{0}_{t}}[x^{*}_{2}(t)])\no\\
& \quad\quad\quad\quad  +1_{0\leq t\leq T-\theta}\nabla \tilde{g}_{32}(x^{*}_{2}(t))+F_{2}(t)^{\top}y_{2}(t)+K_{2}(t)^{\top}z^{0}_{2}(t)+\bar{K}_{2}(t)^{\top}z^{1}_{2}(t)\no\\
&\quad\quad\quad\quad+\bar{F}_{2}(t)^{\top}\mE[y_{1}(t)]+\ddot{F}_{2}(t)\mE^{
\sF_{t}^{0}}[ y^{*}_{1}(t)]+\tilde{F}_{2}(t+\theta)^{\top}\mE^{\sF_{t}}[y_{1}(t+\theta)] ]\dif t\no\\
&\quad\quad\quad\quad\quad\quad\quad\quad+\sum^{d}_{j=1}z^{0}_{2j}(t)\dif B^{0}_{j}(t)+\sum^{d}_{j=1}z^{1}_{2j}(t)\dif B^{1}_{j}(t), t\in [0,T].
\end{align}
with the initial values and terminal values
\begin{align*}
\begin{cases}
& x^{*}_{1}(0)=\frac{\sH}{1-\tau^{2}}\bigg\{\bigg.(\nabla g_{11})^{-1}(-\frac{\sH^{\top}y_{1}(0)-\tau \sH^{\top}y_{2}(0)}{1-\tau^{2}})-\tau\nabla (g_{12})^{-1}(-\frac{\sH^{\top}y_{2}(0)-\tau \sH^{\top}y_{1}(0)}{1-\tau^{2}})\bigg\}\bigg.+x_{0},\\
 &x^{*}_{2}(0)=\frac{\sH}{1-\tau^{2}}\bigg\{\bigg.(\nabla g_{12})^{-1}(-\frac{\sH^{\top}y_{2}(0)-\tau \sH^{\top}y_{1}(0)}{1-\tau^{2}})-\tau\nabla (g_{11})^{-1}(-\frac{\sH^{\top}y_{1}(0)-\tau \sH^{\top}y_{2}(0)}{1-\tau^{2}})\bigg\}\bigg.+x_{0},\\
&x^{*}_{1}(t)=(\nabla \tilde{g}_{31})^{-1}(-( \tilde{F}_{2}(t+\theta)^{\top}\mE[y_{2}(t+\theta)]+\tilde{F}_{1}(t+\theta)^{\top}\mE[y_{1}(t+\theta)])),\\
&x^{*}_{2}(t)=(\nabla \tilde{g}_{32})^{-1}(-( \tilde{F}_{2}(t+\theta)^{\top}\mE[y_{2}(t+\theta)])), t\in [-\theta, 0),\\
& y_{1}(T)=\nabla g_{21}(x^{*}_{1}(T)+x^{*}_{2}(T))+\nabla g_{22}(x^{*}_{1}(T)+x^{*}_{2}(T)),\\
 & y_{2}(T)=\nabla g_{21}(x^{*}_{1}(T)+x^{*}_{2}(T))+\nabla g_{22}(x^{*}_{1}(T)+x^{*}_{2}(T))\\
& y_{1}(t)=\eta_{1}(t),y_{2}(t)=\eta_{2}(t), t\in (T, T+\theta].
\end{cases}
\end{align*}
We denote   \eqref{37a}-\eqref{37d} by   $(\pi_{2}).$
Now, we present our main result as follows:
\bt\label{205}
Let Assumptions 2 and 3 be satisfied. Then, the Hamiltonian system $(\pi_{2})$ possesses a unique solution
$$W_{1}(\cdot):=(x_{1}^{*}(\cdot)^{\top},y_{1}(\cdot)^{\top}, z^{0}_{1}(\cdot)^{\top}, z^{1}_{1}(\cdot)^{\top} )^{\top},  W_{2}(\cdot):=(x^{*}_{2}(\cdot)^{\top},y_{2}(\cdot)^{\top},z^{0}_{2}(\cdot)^{\top},z^{1}_{2}(\cdot)^{\top} )^{\top}.$$
Set
\begin{align*}
\mathfrak{U}_{i}(t):&=\big(\big. (\mE[x_{1}(t)])^{\top}, (\mE[y_{1}(t)])^{\top}, (\mE[x_{2}(t)])^{\top}, (\mE[y_{2}(t)])^{\top}, \\
 &\quad\quad\quad\quad\quad(\mE^{\sF^{0}_{t}}[x_{1}(t)])^{\top}, (\mE^{\sF^{0}_{t}}[y_{1}(t)])^{\top}, (\mE^{\sF^{0}_{t}}[x_{2}(t)])^{\top}, (\mE^{\sF^{0}_{t}}[y_{2}(t)])^{\top},\\
&\quad\quad \quad(x_{1}(t-\theta))^{\top}, (\mE^{\sF_{t}}[y_{1}(t+\theta)])^{\top}, (x_{2}(t-\theta))^{\top}, (\mE^{\sF_{t}}[y_{2}(t+\theta)])^{\top}, \\
  &\quad\quad \quad\quad
  x_{i}(t)^{\top}, y_{i}(t)^{\top}, z^{0}_{i}(t)^{\top}, z^{1}_{i}(t)^{\top}\big)\big.^{\top}, i=1,2.
\end{align*}
Using this solution, the quartet $(\nu^{*}_{1},\nu^{*}_{2},u^{*}_{1},u^{*}_{2})$ defined by \eqref{36} and \eqref{36++++-}
constitutes the unique optimal controls for Problem (LC)
\et
\begin{proof}
First, we aim to demonstrate the unique solvability of  System $(\pi_{2})$.
The Hamiltonian system $(\pi_{2})$ takes the form of the System $(\pi)$. To utilize Theorem  \ref{tt} for establishing the unique solvability of $(\pi_{2})$, it is sufficient to confirm Assumption 1. In the subsequent parts, we will present a detailed verification solely of the monotonicity condition for the coefficients  $\aleph_{i}(\cdot),i=1,2$ [cf. Assumption 1(iii)3)], while omitting other specifics.  Suppose that  $\bar{W}_{1}(\cdot):=(\bar{x}_{1}(\cdot)^{\top},\bar{y}_{1}(\cdot)^{\top}, \bar{z}^{0}_{1}(\cdot)^{\top}, \bar{z}^{1}_{1}(\cdot)^{\top} )^{\top},  \bar{W}_{2}(\cdot):=(\bar{x}_{2}(\cdot)^{\top},\bar{y}_{2}(\cdot)^{\top},\bar{z}^{0}_{2}(\cdot)^{\top},\bar{z}^{1}_{2}(\cdot)^{\top} )^{\top}$ is a solution of the Hamiltonian system $(\pi_{2})$, $\tilde{W}_{1}(\cdot):=(\tilde{x}_{1}(\cdot)^{\top},\tilde{y}_{1}(\cdot)^{\top}, \tilde{z}^{0}_{1}(\cdot)^{\top}, \tilde{z}^{1}_{1}(\cdot)^{\top} )^{\top},$  $\tilde{W}_{2}(\cdot):=(\tilde{x}_{2}(\cdot)^{\top},\tilde{y}_{2}(\cdot)^{\top},\tilde{z}^{0}_{2}(\cdot)^{\top},\tilde{z}^{1}_{2}(\cdot)^{\top}  )^{\top}$ is another solution of the Hamiltonian system $(\pi_{2})$ and set
\begin{align*}
\bar{\mathfrak{U}}_{i}(t)&:=((\mE[\bar{x}_{1}(t)])^{\top}, (\mE[\bar{y}_{1}(t)])^{\top}), (\mE[\bar{x}_{2}(t)])^{\top}, (\mE[\bar{y}_{2}(t)])^{\top}, \\
 &\quad\quad\quad\quad\quad(\mE^{\sF^{0}_{t}}[x_{1}(t)])^{\top}, (\mE^{\sF^{0}_{t}}[y_{1}(t)])^{\top}, (\mE^{\sF^{0}_{t}}[x_{2}(t)])^{\top}, (\mE^{\sF^{0}_{t}}[y_{2}(t)])^{\top},\\
&\quad\quad \quad(\bar{x}_{1}(t-\theta))^{\top}, (\mE^{\sF_{t}}[\bar{y}_{1}(t+\theta)])^{\top}, (\bar{x}_{2}(t-\theta))^{\top}, (\mE^{\sF_{t}}[\bar{y}_{2}(t+\theta)])^{\top}, \\
 &\quad\quad\quad\quad\quad\quad\quad\quad\quad\quad\quad\quad \bar{x}_{i}(t)^{\top}, \bar{y}_{i}(t)^{\top}, \bar{z}^{0}_{i}(t)^{\top}, \bar{z}^{1}_{i}(t)^{\top})^{\top}, i=1,2.
\end{align*}

\begin{align*}
\tilde{\mathfrak{U}}_{i}(t)&:=((\mE[\tilde{x}_{1}(t)])^{\top}, (\mE[\tilde{y}_{1}(t)])^{\top}, (\mE[\tilde{x}_{2}(t)])^{\top},(\mE[\tilde{y}_{2}(t)])^{\top}, \\
&\quad\quad\quad\quad\quad(\mE^{\sF^{0}_{t}}[x_{1}(t)])^{\top}, (\mE^{\sF^{0}_{t}}[y_{1}(t)])^{\top}, (\mE^{\sF^{0}_{t}}[x_{2}(t)])^{\top}, (\mE^{\sF^{0}_{t}}[y_{2}(t)])^{\top},\\
&\quad\quad \quad(\tilde{x}_{1}(t-\theta))^{\top}, (\mE^{\sF_{t}}[\tilde{y}_{1}(t+\theta)])^{\top}, (\tilde{x}_{2}(t-\theta))^{\top}, (\mE^{\sF_{t}}[\tilde{y}_{2}(t+\theta)])^{\top}, \\
& \quad\quad\quad\quad\quad\quad\quad\quad\quad\quad\quad\quad \tilde{x}_{i}(t)^{\top}, \tilde{y}_{i}(t)^{\top}, \tilde{z}^{0}_{i}(t)^{\top}, \tilde{z}^{1}_{i}(t)^{\top})^{\top}, i=1,2 .
\end{align*}
The following result follows from a straightforward calculation and denote $\hat{l}:=\tilde{l}-\bar{l}.$
\begin{align}\label{302}
&\mE\int^{T}_{0}[\langle\aleph_{1}(\tilde{\mathfrak{U}}_{1}(t)) -\aleph_{1}(\bar{\mathfrak{U}}_{1}(t)),\tilde{W}_{1}(t)-\bar{W}_{1}(t)\rangle+\langle\aleph_{2}(\tilde{\mathfrak{U}}_{2}(t))-\aleph_{2}(\bar{\mathfrak{U}}_{2}(t)),
\tilde{W}_{2}(t)-\bar{W}_{2}(t)\rangle]\dif t\no\\
&=-\mE\int^{T}_{0}\langle \nabla g_{31}(\tilde{x}_{1}(t))-\nabla g_{31}(\bar{x}_{1}(t)),   \hat{x}_{1}(t)\rangle\dif t\no\\
&-\int^{T}_{0}\mE\langle \nabla g_{32}(\tilde{x}_{2}(t))-\nabla g_{32}(\bar{x}_{2}(t)),   \hat{x}_{2}(t)\rangle\dif t\no\\
&\quad-\mE\int^{T}_{0}\langle \nabla \bar{g}_{31}(\mE[\tilde{x}_{1}(t)])-\nabla g_{31}(\mE[\bar{x}_{1}(t)]),   \mE[\hat{x}_{1}(t)]\rangle\dif t\no\\
&-\mE\int^{T}_{0}\langle \nabla \bar{g}_{32}(\mE[\tilde{x}_{2}(t)])-\nabla \bar{g}_{32}(\mE[\bar{x}_{2}(t)]),   \mE[\hat{x}_{2}(t)]\rangle\dif t\no\\
&\quad-\mE\int^{T}_{0}\langle \nabla \ddot{g}_{31}(\mE^{\sF^{0}_{t}}[\tilde{x}_{1}(t)])-\nabla g_{31}(\mE^{\sF^{0}_{t}}[\bar{x}_{1}(t)]),   \mE^{\sF^{0}_{t}}[\hat{x}_{1}(t)]\rangle\dif t\no\\
&-\mE\int^{T}_{0}\langle \nabla \ddot{g}_{32}(\mE^{\sF^{0}_{t}}[\tilde{x}_{2}(t)])-\nabla \bar{g}_{32}(\mE^{\sF^{0}_{t}}[\bar{x}_{2}(t)]),   \mE^{\sF^{0}_{t}}[\hat{x}_{2}(t)]\rangle\dif t\no\\
&\quad-\mE\int^{T}_{0}1_{0\leq t\leq T-\theta}\langle \nabla \tilde{g}_{31}(\tilde{x}_{1}(t))-\nabla \tilde{g}_{31}(\bar{x}_{1}(t)),   \hat{x}_{1}(t)\rangle\dif t\no\\
&-\mE\int^{T}_{0}1_{0\leq t\leq T-\theta}\langle \nabla \tilde{g}_{32}(\tilde{x}_{2}(t))-\nabla \tilde{g}_{32}(\bar{x}_{2}(t)),   \hat{x}_{2}(t)\rangle\dif t\no\\
&+\mE\int^{T}_{0}\langle\hat{\alpha}_{1}(t), H^{\top}_{1}(t)\hat{y}_{1}(t)+\tau\bar{H}^{\top}_{1}(t)\mE[\hat{y}_{1}(t)]+\tau\ddot{H}^{\top}_{1}(t)\mE^{\sF_{t}^{0}}[\hat{y}_{1}(t)]\no\\
&\quad \quad\quad\quad \quad\quad\quad \quad\quad\quad \quad+\tau\tilde{H}^{\top}_{1}(t+\theta)\mE^{\sF_{t}}[\hat{y}_{1}(t+\theta)] +D^{\top}_{1}(t)\hat{z}^{0}_{1}(t)+\bar{D}^{\top}_{1}(t)\hat{z}^{1}_{1}(t)\rangle\dif t\no\\
&+\mE\int^{T}_{0}\langle\hat{\alpha}_{2}(t), H^{\top}_{2}(t)\hat{y}_{2}(t)+\tau\bar{H}^{\top}_{2}(t)\mE[\hat{y}_{2}(t)]+\tau\ddot{H}^{\top}_{2}(t)\mE^{\sF_{t}^{0}}[\hat{y}_{2}(t)]\no\\
&\quad \quad\quad\quad \quad\quad\quad \quad\quad\quad \quad+\tau\tilde{H}^{\top}_{2}(t+\theta)\mE^{\sF_{t}}[\hat{y}_{2}(t+\theta)] +D^{\top}_{2}(t)\hat{z}^{0}_{2}(t)+\bar{D}^{\top}_{2}(t)\hat{z}^{1}_{2}(t)\rangle\dif t\no\\
&+\int^{0}_{-\theta}\langle\hat{\beta}_{1}(t),    \tau \tilde{H}_{1}(t+\theta)^{\top}\mE[\hat{y}_{1}(t+\theta)]\rangle \dif t
+\int^{0}_{-\theta}\langle\hat{\beta}_{2}(t),    \tau \tilde{H}_{2}(t+\theta)^{\top}\mE[\hat{y}_{2}(t+\theta)]\rangle \dif t
\end{align}
where
\begin{align*}
&\tilde{\alpha}_{1}(t):= (\nabla g_{41})^{-1}(-( H^{\top}_{1}(t)\tilde{y}_{1}(t)+\tau\bar{H}^{\top}_{1}(t)\mE[\tilde{y}_{1}(t)]+\tau\ddot{H}^{\top}_{1}(t)\mE^{\sF_{t}^{0}}[\tilde{y}_{1}(t)]\\
&\quad\quad\quad\quad\quad\quad\quad\quad+\tau\tilde{H}^{\top}_{1}(t+\theta)\mE^{\sF_{t}}[\tilde{y}_{1}(t+\theta)]
+D^{\top}_{1}(t)\tilde{z}^{0}_{1}(t)+\bar{D}^{\top}_{1}(t)\tilde{z}^{1}_{1}(t))), t\in [0,T],\\
&\tilde{\alpha}_{2}(t):= (\nabla g_{41})^{-1}(-( H^{\top}_{2}(t)\tilde{y}_{2}(t)+\tau\bar{H}^{\top}_{2}(t)\mE[\tilde{y}_{2}(t)]+\tau\ddot{H}^{\top}_{2}(t)\mE^{\sF_{t}^{0}}[\tilde{y}_{2}(t)]\\
&\quad\quad\quad\quad\quad\quad\quad\quad+\tau\tilde{H}^{\top}_{2}(t+\theta)\mE^{\sF_{t}}[\tilde{y}_{2}(t+\theta)]
+D^{\top}_{2}(t)\tilde{z}^{0}_{2}(t)+\bar{D}^{\top}_{2}(t)\tilde{z}^{1}_{2}(t))), t\in [0,T],\\
&\bar{\alpha}_{1}(t):= (\nabla g_{41})^{-1}(-( H^{\top}_{1}(t)\bar{y}_{1}(t)+\tau\bar{H}^{\top}_{1}(t)\mE[\bar{y}_{1}(t)]+\tau\ddot{H}^{\top}_{1}(t)\mE^{\sF_{t}^{0}}[\bar{y}_{1}(t)]\\
&\quad\quad\quad\quad\quad\quad\quad\quad
+\tau\bar{H}^{\top}_{1}(t+\theta)\mE^{\sF_{t}}[\bar{y}_{1}(t+\theta)]+D^{\top}_{1}(t)\bar{z}^{0}_{1}(t)+\bar{D}^{\top}_{1}(t)\bar{z}^{1}_{1}(t))),t\in [0,T],\\
&\bar{\alpha}_{2}(t):= (\nabla g_{41})^{-1}(-( H^{\top}_{2}(t)\bar{y}_{2}(t)+\tau\bar{H}^{\top}_{2}(t)\mE[\bar{y}_{2}(t)]+\tau\ddot{H}^{\top}_{2}(t)\mE^{\sF_{t}^{0}}[\bar{y}_{2}(t)]\\
&\quad\quad\quad\quad\quad\quad\quad\quad+\tau\bar{H}^{\top}_{2}(t+\theta)\mE^{\sF_{t}}[\bar{y}_{2}(t+\theta)]
+D^{\top}_{2}(t)\bar{z}^{0}_{2}(t)+\bar{D}^{\top}_{2}(t)\bar{z}^{1}_{2}(t))),t\in [0,T],\\
&\tilde{\beta}_{1}(t):= (\nabla \bar{g}_{41})^{-1}(-( \tau \tilde{H}_{1}(t+\theta)^{\top}\mE[\tilde{y}_{1}(t+\theta)])),t\in [-\theta,0),\\
&\tilde{\beta}_{2}(t):=(\nabla \bar{g}_{42})^{-1}(-( \tau \tilde{H}_{2}(t+\theta)^{\top}\mE[\tilde{y}_{2}(t+\theta)])),t\in [-\theta,0),\\
&\bar{\beta}_{1}(t):=(\nabla \bar{g}_{41})^{-1}(-( \tau \tilde{H}_{1}(t+\theta)^{\top}\mE[\bar{y}_{1}(t+\theta)])),t\in [-\theta,0),\\
&\bar{\beta}_{2}(t):= (\nabla \bar{g}_{42})^{-1}(-( \tau \tilde{H}_{2}(t+\theta)^{\top}\mE[\bar{y}_{2}(t+\theta)])),t\in [-\theta,0).
 \end{align*}
Given the convexity of $g_{31}, g_{32}, \bar{g}_{31}, \bar{g}_{32},\tilde{g}_{31}, \tilde{g}_{32}$ and uniform convexity of $g_{41}, g_{42},$ $\bar{g}_{41}, \bar{g}_{42},$ it follows from Lemma \ref{100} that
\begin{align}\label{301}
&-\mE\langle \nabla g_{31}(\tilde{x}_{1}(t))-\nabla g_{31}(\bar{x}_{1}(t)),   \hat{x}_{1}(t)\rangle
-\mE\langle \nabla g_{32}(\tilde{x}_{2}(t))-\nabla g_{32}(\bar{x}_{2}(t)),   \hat{x}_{2}(t)\rangle\no\\
&-\mE\langle \nabla \bar{g}_{31}(\mE[\tilde{x}_{1}(t)])
-\nabla g_{31}(\mE[\bar{x}_{1}(t)]),   \mE[\hat{x}_{1}(t)]\rangle\no\\
&-\mE\langle \nabla \bar{g}_{32}(\mE[\tilde{x}_{2}(t)])
-\nabla \bar{g}_{32}(\mE[\bar{x}_{2}(t)]),   \mE[\hat{x}_{2}(t)]\rangle\no\\
&-\mE\langle \nabla \ddot{g}_{31}(\mE^{\sF_{t}^{0}}[\tilde{x}_{1}(t)])
-\nabla \ddot{g}_{31}(\mE^{\sF_{t}^{0}}[\bar{x}_{1}(t)]),   \mE^{\sF_{t}^{0}}[\hat{x}_{1}(t)]\rangle\no\\
&-\mE\langle \nabla \ddot{g}_{32}(\mE^{\sF_{t}^{0}}[\tilde{x}_{2}(t)])
-\nabla \ddot{g}_{32}(\mE^{\sF_{t}^{0}}[\bar{x}_{2}(t)]),   \mE^{\sF_{t}^{0}}[\hat{x}_{2}(t)]\rangle\no\\
&-\mE\langle \nabla \tilde{g}_{31}(\tilde{x}_{1}(t))
-\nabla \tilde{g}_{31}(\bar{x}_{1}(t)),  \hat{x}_{1}(t)\rangle\no\\
&-\mE\langle \nabla \tilde{g}_{32}(\tilde{x}_{2}(t))
-\nabla \tilde{g}_{32}(\bar{x}_{2}(t)),  \hat{x}_{2}(t)\rangle \leq 0,
\end{align}
and
\begin{align}\label{301+}
&\mE\langle  \hat{\alpha}_{1}(t),   H_{1}(t)^{\top}\hat{y}_{1}(t)+ \tau\ddot{H}_{1}(t)^{\top}\mE^{\sF^{0}_{t}}[\hat{y}_{1}(t)]+ \tau\bar{H}_{1}(t)^{\top}\mE[\hat{y}_{1}(t)]\no\\
&\quad\quad\quad\quad\quad\quad\quad\quad+ \tau\tilde{H}_{1}(t+\theta)^{\top}\mE^{\sF_{t}}[\hat{y}_{1}(t+\theta)]+D_{1}(t)^{\top}\hat{z}^{0}_{1}(t)+\bar{D}_{1}(t)^{\top}\hat{z}^{1}_{1}(t)\rangle\no\\
&\quad +\mE\langle  \hat{\alpha}_{2}(t),   H_{2}(t)^{\top}\hat{y}_{2}(t)+ \tau\ddot{H}_{2}(t)^{\top}\mE^{\sF^{0}_{t}}[\hat{y}_{21}(t)]+ \tau\bar{H}_{2}(t)^{\top}\mE[\hat{y}_{2}(t)]\no\\
&\quad\quad\quad\quad\quad\quad\quad\quad+ \tau\tilde{H}_{2}(t+\theta)^{\top}\mE^{\sF_{t}}[\hat{y}_{2}(t+\theta)]+D_{2}(t)^{\top}\hat{z}^{0}_{2}(t)+\bar{D}_{2}(t)^{\top}\hat{z}^{1}_{2}(t)\rangle\no\\
&\leq -\mE\langle \hat{\alpha}_{2}(t), \nabla g_{42}(\tilde{\alpha}_{2}(t))- \nabla g_{42}(\bar{\alpha}_{2}(t))   \rangle -\langle \hat{\alpha}_{2}(t), \nabla g_{42}(\tilde{\alpha}_{2}(t))- \nabla g_{41}(\bar{\alpha}_{2}(t))   \rangle\no\\
&\leq -\delta\mE[|\hat{\alpha}_{1}(t)|^{2}] -\delta\mE[|\hat{\alpha}_{2}(t)|^{2}].
\end{align}
Similarly, we can prove that
\begin{align}\label{301+x}
\langle\hat{\beta}_{1}(t),    \tau \tilde{H}_{1}(t+\theta)^{\top}\mE[\hat{y}_{1}(t+\theta)]\rangle
+\langle\hat{\beta}_{2}(t),    \tau \tilde{H}_{2}(t+\theta)^{\top}\mE[\hat{y}_{2}(t+\theta)]\rangle \leq 0.
\end{align}
Then, We can derive  the corresponding monotonicity conditions for $\aleph_{1}, \aleph_{2}$ is satisfied by giving the following  definitions:
\begin{align}
\begin{cases}
&\bar{\chi}_{11}(v):=\frac{1}{1-\tau^{2}}\bigg.(\nabla g_{11})^{-1}(-v), \bar{\chi}_{12}(v):=-\frac{\tau}{1-\tau^{2}}\nabla (g_{12})^{-1}(-v),\\
&\bar{\chi}_{21}(v):=\frac{1}{1-\tau^{2}}\bigg.(\nabla g_{12})^{-1}(-v), \bar{\chi}_{22}(v):=-\frac{\tau}{1-\tau^{2}}\nabla (g_{11})^{-1}(-v),\\
&  \ddot{\chi}_{1}(t,y_{1},y_{2}):=(\nabla \tilde{g}_{31})^{-1}(t,-(F_{1}(t+\theta)^{\top}y_{1}+F_{2}(t+\theta)^{\top}y_{2})),t\in [-\theta, 0),\\
&  \ddot{\chi}_{2}(t,y_{1},y_{2}):=(\nabla \tilde{g}_{32})^{-1}(t,-(F_{2}(t+\theta)^{\top}y_{2})),t\in [-\theta, 0),\\
& \chi_{2}(t,u):=(\nabla g_{42})^{-1}(t,-u), t\in [0,T],\\
&  \chi_{1}(t,u):=(\nabla g_{41})^{-1}(t,-u), \chi_{2}(t,u):=(\nabla g_{42})^{-1}(t,-u), t\in [0,T].
\end{cases}
\end{align}

Now, we intend to prove that $(\nu^{*}_{1},\nu^{*}_{2},u^{*}_{1},u^{*}_{2})$  is  the unique quartet of optimal controls. Let $(\nu_{1},\nu_{2},u_{1},u_{2})$
be another arbitrary admissible  quartet. By $(\nu^{*}_{1},\nu^{*}_{2},u^{*}_{1},u^{*}_{2})$ in \eqref{36},   we know that $(\Gamma^{*}_7, \Gamma^{*}_8, \Gamma^{*}_9, \Gamma^{*}_{10}) = (0, 0,0,0)$  and $(\bar{\Gamma}^{*}_7, \bar{\Gamma}^{*}_8, \bar{\Gamma}^{*}_9, \bar{\Gamma}^{*}_{10}) = (0, 0,0,0).$  Then,   $\eqref{35}$  is reduced to

\begin{align}\label{41}
&J(\nu_{1},\nu_{2}, u_{1}, u_{2}) - J(\nu^{*}_{1},\nu^{*}_{2},u_{1}^{*},u_{2}^{*})\no\\
&\geq \frac{\delta}{2}\bigg\{\bigg. |\xi_{1}+\tau\xi_{2} - \xi^*_{1}-\tau\xi^*_{2}|^2+ |\tau\xi_{1}+\xi_{2} - \tau\xi^*_{1}-\xi^*_{2}|^2\no\\
& +\mE\int^{T}_{0} |u_{1}(t)-u^{*}_{1}(t)|^{2}\dif t +\mE\int^{T}_{0} |u_{2}(t)-u^{*}_{2}(t)|^{2}\dif t    \bigg\}\bigg.> 0.
\end{align}

Thus, $(\nu^{*}_{1},\nu^{*}_{2},u^{*}_{1},u^{*}_{2})$  is the unique quartet of optimal controls for Problem (LC). The proof is thereby established.
\end{proof}

\begin{exa}  Let $n = d = m = k = T = 1, \eta_{i}(t)=0, t\in (T,T+\theta], i=1,2$. Consider the following controlled system:
\begin{equation*}
\begin{cases}
&dx_{1}(t) = \bigg\{\bigg. \frac{1}{1000}\mE[x_{2}(t)]+\frac{1}{1000}\mE^{\sF^{0}_{t}}[x_{2}(t)]+\frac{1}{1000}1_{\{0\leq t\leq \theta\}}x_{2}(t-\theta)\bigg\}\bigg. \dif t\\
&\quad\quad\quad\quad\quad\quad+B^{1}(t) \dif B^{0}(t)+[u_{1}(t) + \sin(B^{1}(t))] \dif B^{1}(t),\\
&dx_{2}(t) = \bigg\{\bigg.\frac{1}{1000}\mE[x_{1}(t)]+\frac{1}{1000}\mE^{\sF^{0}_{t}}[x_{1}(t)]+\frac{1}{1000}1_{\{0\leq t\leq \theta\}}x_{1}(t-\theta)\bigg\}\bigg. \dif t\\
&\quad\quad\quad\quad\quad\quad+B^{1}(t) \dif B^{0}(t)+[u_{2}(t) + \sin(B^{1}(t))] \dif B^{1}(t), t\in [0,T]\\
&x_{1}(0) = \xi_{1},x_{2}(0) = \xi_{2},\\
& x_{1}(t)=\kappa_{1}(t), x_{2}(t)=\kappa_{2}(t), t\in [-\theta, 0).
\end{cases}
\end{equation*}
and the following criterion functional:
\begin{align*}
J(\nu_{1},\nu_{2}, u_{1}(\cdot), u_{2}(\cdot)) &= g(\xi_{1}+\frac{1}{1000}\xi_{2})+g(\xi_{2}+\frac{1}{1000}\xi_{1}) + \mathbb{E} \bigg\{\bigg. \frac{1}{2} |x_{1}(1)+x_{2}(1)|^2 \\
&\quad +\frac{1}{2} |x_{1}(1)+x_{2}(1)|^2+ \int_0^1 g(u_{1}(s)) \dif s+ \int_0^1 g(u_{2}(s)) \dif s \bigg\}\bigg..
\end{align*}
where $g(\cdot)$ is given by \begin{equation}
g_{1}(u) = g_{4}(t,u) = g(u) := \begin{cases} e^{u}-u-1, & u \geq 0 \\ e^{-u}+u-1, & u < 0. \end{cases}
\end{equation}

It is clear that the function  $g(\cdot)$ described above is continuously differentiable.
More specifically,

\begin{equation}\label{201}
\nabla g(u) = \frac{dg}{du}(u) = \begin{cases} e^{u}-1, & u \geq 0 \\ 1-e^{-u}, & u < 0. \end{cases}
\end{equation}
Since $(d^2 g)/(du^2)(u) \geq 1$ for any $u \in \mathbb{R},$
we have
$\langle \nabla g(u) - \nabla g(\bar{u}), u - \bar{u} \rangle \geq |u - \bar{u}|^2$ for  any $u,\bar{u}\in \mR.$
Thus, the Lemma \ref{100} in Appendix  shows that $f(\cdot)$ is uniformly convex. Furthermore, since $f(0)=\nabla f(0)=0,$  $Assumption\, 3$ $(3) (4)$ hold.

From \eqref{201}, we have

$$ (\nabla g)^{-1}(u) = \begin{cases}
\ln(1+u), & u \geq 0 \\
-\ln(1-u), & u < 0.
\end{cases} $$
For  $i=1,2,j=1,2,$ set
\begin{align}
\bar{\chi}_{ij}(u) = \chi_{i}(t, u) = (\nabla f)^{-1}(-u)
=
\begin{cases}
-\ln(1+u), & u \geq 0, \\
\ln(1-u), & u < 0 ,
\end{cases}
\end{align}
Obviously, we know that
\begin{align}\label{203}
\frac{\partial \chi_{i}}{\partial u}(t, u) = -\frac{1}{1 + |u|} \in [-1, 0), i=1,2.
\end{align}
On the one hand, we claim that the previously defined $\chi_{i}(\cdot,\cdot)$ satisfies $Assumption\, 1(\mathrm{iii}).$ In fact,, \eqref{203} implies the Lipschitz continuity, meaning  the third inequality of $\mathrm{Assumption\, 1\,(iii)\,1)}$ is valid. furthermore, for any $-\infty < \bar{u} < u < \infty$, the Lagrange's mean value theorem  guarantees the existence of a value  $\tilde{u} \in (\bar{u}, u)$ such that
\begin{equation*}
    \chi_{i}(t, u) - \chi_{i}(t, \bar{u}) = \frac{\partial \chi_{i}}{\partial u}(t, \tilde{u})(u - \bar{u}) \geq -(u - \bar{u})
\end{equation*}
i.e., $u - \bar{u} \geq -(\chi_{i}(t, u) - \chi_{i}(t, \bar{u}))$. Therefore,
\begin{equation*}
    (\chi(t, u) - \chi(t, \bar{u}))(u - \bar{u}) \leq -|\chi(t, u) - \chi(t, \bar{u})|^2.
\end{equation*}
 The fifth inequality of $\mathrm{Assumption\, 1\,(iii)\,1)}$ is hereby confirmed.

In addition, owing to \eqref{203} and Lagrange's mean value theorem, there is no constant $c > 0$ for which
\begin{equation*}
    |\chi(t, u) -\chi(t, \bar{u})|^2 \geq c|u - \bar{u}|^2 \quad \text{for any } u, \bar{u} \in \mathbb{R}.
\end{equation*}
This indicates that we cannot replace $h$ with $-u$, i.e., the linear version of domination-monotonicity conditions in \cite{yu2022} does not hold.

Clearly, Assumptions 2 and 3 hold true in this case.

Theorem \ref{205} implies that Problem (LC) with the special setting above admits a unique quartet  of optimal controls$(\tau=\frac{1}{1000}).$
\begin{align}\label{101}
\begin{cases}
\xi^{*}_{1}=\frac{1}{1-\tau^{2}}\bigg\{\bigg.(\nabla g)^{-1}(-\frac{y_{1}(0)-\tau y_{2}(0)}{1-\tau^{2}})-\tau\nabla (g)^{-1}(-\frac{y_{2}(0)-\tau y_{1}(0)}{1-\tau^{2}})\bigg\}\bigg.,\\
\xi^{*}_{2}=\frac{1}{1-\tau^{2}}\bigg\{\bigg.(\nabla g)^{-1}(-\frac{y_{2}(0)-\tau y_{1}(0)}{1-\tau^{2}})-\tau(\nabla g)^{-1}(-\frac{y_{1}(0)-\tau y_{2}(0)}{1-\tau^{2}})\bigg\}\bigg.,\\
\kappa_{1}(t)=0, t\in [-\theta,0),\\
\kappa_{1}(t)=0,t\in [-\theta,0),\\
u^{*}_{1}(t)=(\nabla g)^{-1}(-z^{1}_{1}(t)),t\in [0,T],\\
u^{*}_{2}(t)=(\nabla g)^{-1}(- z^{1}_{2}(t)),t\in [0,T],\\
u^{*}_{1}(t)=0,t\in [-\theta,0),\\
u^{*}_{2}(t)=0,t\in [-\theta,0),
\end{cases}
\end{align}
where $(x^*_{i}(\cdot)^{\top}, y_{i}(\cdot)^{\top}, z^{0}_{i}(\cdot)^{\top}, z^{1}_{i}(\cdot)^{\top})^{\top} \in \mM_{\mF}(\mathbb{R}^{1+1+1+1}),i=1,2$ is the unique solution to the following stochastic Hamiltonian system:
\begin{equation}\label{102}
\begin{cases}
&dx^{*}_{1}(t) = \bigg\{\bigg.  \frac{1}{1000}\mE[x_{2}(t)]+\frac{1}{1000}\mE^{\sF^{0}_{t}}[x_{2}(t)]+\frac{1}{1000}1_{\{0\leq t\leq \theta\}}x_{2}(t-\theta)\bigg\}\bigg.\dif t \\
&\quad\quad\quad\quad\quad\quad\quad\quad\quad\quad\quad+ B^{1}(t) \dif B^{0}(t)
+[(\nabla g)^{-1}(- z^{1}_{1}(t)) + \sin(B^{1}(t))] \dif B^{1}(t),\\
&dy_{1}(t) =\bigg\{\bigg.-\frac{1}{1000}\mE[y_{2}(t)]-\frac{1}{1000}\mE^{\sF^{0}_{t}}[y_{2}(t)]-\frac{1}{1000}1_{\{-\theta\leq t\leq 0\}}\mE[y_{2}(t+\theta)]\bigg\}\bigg.\dif t\\
& \quad\quad\quad\quad\quad\quad\quad\quad\quad\quad\quad\quad\quad\quad\quad\quad\quad+ z^{0}_{1}(t) \dif B^{0}(t)+ z^{1}_{1}(t) \dif B^{1}(t),\\
&dx^{*}_{2}(t) =\bigg\{\bigg.  \frac{1}{1000}\mE[x_{1}(t)]+\frac{1}{1000}\mE^{\sF^{0}_{t}}[x_{1}(t)]+\frac{1}{1000}1_{\{0\leq t\leq \theta\}}x_{1}(t-\theta) \bigg\}\bigg. \dif t\\
&\quad\quad\quad\quad\quad\quad\quad\quad\quad\quad\quad+ B^{1}(t) \dif B^{0}(t)
+[(\nabla g)^{-1}(- z^{1}_{2}(t)) + \sin(B^{1}(t))] \dif B^{1}(t),\\
& \quad\quad\quad\quad\quad\quad\quad\quad\quad\quad\quad\quad\quad\quad\quad\quad\quad+ z^{0}_{1}(t) \dif B^{0}(t)+ z^{1}_{1}(t) \dif B^{1}(t),\\
&dy_{2}(t) =\bigg\{\bigg.-\frac{1}{1000}\mE[y_{1}(t)]-\frac{1}{1000}\mE^{\sF^{0}_{t}}[y_{1}(t)]-\frac{1}{1000}1_{\{-\theta\leq t\leq 0\}}\mE[y_{1}(t+\theta)]\bigg\}\bigg.\dif t\\
& \quad\quad\quad\quad\quad\quad\quad\quad\quad\quad\quad\quad\quad\quad \quad\quad\quad + z^{0}_{2}(t) \dif B^{0}(t)+ z^{1}_{2}(t) \dif B^{1}(t),\\
&x^{*}_{1}(0) = \xi^{*}_{1},x^{*}_{2}(0) = \xi^{*}_{2},\\
& x^{*}_{1}(t)=\sH\xi^{*}_{1}+\kappa^{*}_{1}(t), x^{*}_{2}(t)=\sH\xi^{*}_{2}+\kappa^{*}_{2}(t), t\in [-\theta, 0),\\
& y_{1}(1) =\frac{1}{2}( x^{*}_{1}(1)+x^{*}_{2}(1)), y_{2}(1) =\frac{1}{2}( x^{*}_{1}(1)+x^{*}_{2}(1)),\\
& y_{1}(t)=0,y_{2}(t)=0, t\in (T, T+\theta].
\end{cases}
\end{equation}
Set
\begin{align*}
\nu^{*}_{i}(t)=
\begin{cases}
&\sH\xi^{*}_{i}+x_{0}, t=0,\\
&\kappa^{*}_{i}(t), t\in [-\theta, 0).
\end{cases}
\end{align*}
As a matter of fact, we can verify that the unique  solution to\eqref{102} takes on the subsequent form:
\begin{equation}\label{103}
\begin{cases}
\nu^{*}_{1}(\cdot)=\nu^{*}_{2}(\cdot),u^{*}_{1}(\cdot)=u^{*}_{2}(\cdot) \\
z^{0}_{1}(t)=z^{0}_{2}(t) =  B^{1}(t),\\
z^{1}_{1}(t)=z^{1}_{2}(t) = u^*_{i}(t) + \sin(B^{1}(t)),\\
\mE[x^*_{1}(t)]=\mE[y_{1}(t)]=\mE[x^*_{2}(t)]=\mE[y_{2}(t)]=0,\\
\mE^{\sF^{0}_{t}}[x^*_{1}(t)]=\mE^{\sF^{0}_{t}}[y_{1}(t)]=\mE^{\sF^{0}_{t}}[x^*_{2}(t)]=\mE^{\sF^{0}_{t}}[y_{2}(t)]=0, t\in [0,T]\\
\end{cases}
\end{equation}
for any $(\omega, t) \in \Omega \times [0, 1]$, where the unique optimal pair $(\nu^*_{1},\nu^*_{1}, u^*_{1}, u^*_{2})$ satisfies that $\nu^*_{i}(t) \equiv 0,i=1,2,$ and
\begin{equation}
\begin{cases}
e^{u^*_{i}(t)} - 1 + u^*_{i}(t) + B^{1}(t)+ \sin(B^{1}(t)) = 0, \quad u^*_{i}(t) \geq 0\\
-e^{-u^*_{i}(t)} + 1 + u^*_{i}(t) + B^{1}(t)+ \sin(B^{1}(t))  = 0, \quad u^*_{i}(t) < 0, t\in [0,T],\\
u^*_{i}(t)=0, t\in [-\theta,0),
i=1,2.
\end{cases}
\end{equation}

\end{exa}

\section{Application to stochastic linear-quadratic problems with input constraints}
 Stochastic linear-quadratic (LQ) problems with input constraints are of paramount importance in both control theory and real-world applications because they address the fundamental need to make optimal decisions under uncertainty while respecting physical or operational limits.  the LQ framework provides a tractable yet powerful setting for analyzing systems with linear dynamics and quadratic costs, which naturally penalize deviations from desired states and excessive control effort. When randomness is introduced (e.g., process noise, uncertain parameters), the problem becomes a stochastic optimal control problem-critical for systems where future conditions are imperfectly known.  Prior to presenting the main results of this section, we present some illustrations. In this section, we turn our attention to the (LQ-IC) problem introduced in Section 1 under Assumptions 2 and 4. Consequently, it follows that
\begin{align*}
&|\mJ(\nu_{1},\nu_{2}, u_{1},u_{2})|< \infty,\,\\
&\mbox{for any}\, (\nu_{1},\nu_{2}, u_{1},u_{2})\in L^{2}_{\mF}(-\theta, 0; \mR^{n})\times L^{2}_{\mF}(-\theta, 0; \mR^{n})
\times L^{2}_{\mF}(-\theta,T;\mR^{k})
\times L^{2}_{\mF}(-\theta,T;\mR^{k}).
\end{align*}
In contrast to Problem  (LC), the controls (inputs) in Problem (LQ-IC) are restricted to closed convex sets. These constraints on $u_i(\cdot), i=1,2$ can exhibit time-dependency and randomness. Our approach can still be applied even when the terms containing  $x_{i}(T),i=1,2$ and $x_{i}(\cdot),i=1,2$ are replaced by other convex functions with appropriate properties. Nevertheless, we employ the quadratic criterion functional \eqref{6} to present our main idea with clarity and conciseness.  Here, the constraint set $\mU(\cdot)$ is influenced  on both time and randomness, mirroring how control constraints vary with time and circumstances.  We propose adopting Assumption 5 as a workable solution. Despite this, it is not an overly strict condition and is applicable in numerous scenarios.
To be precise, Assumption 5 obviously holds when $\mU$ is independent of $(\omega, t)$ and Assumption 5(i) is satisfied.
Moreover, when $\mU(\omega, t)$ is uniformly bounded for almost all $(\omega, t) \in \Omega \times [-\theta, T]$, Assumption 5(iii) is superfluous.  In fact, under Assumption 5(ii),

$$
\Lambda = \{ (\omega, t, u) \in \Omega \times [-\theta, T] \times \mathbb{R}^k | 1_{\mU(\omega, t)}(u) = 1 \}
$$
is $\sP \otimes \mathcal{B}(\mathbb{R}^k)$-measurable.  Consequently, according to the measurable selection theorem, there is a $\mathcal{P}$-measurable process $a(\cdot)$ such that $a(\omega, t) \in \mU(\omega, t)$.

To define the optimal control pair for Problem (LQ-IC), we introduce the Hamiltonian system below:

\begin{align}\label{43m}
\dif x^{*}_{1}(t)&=[F_{1}(t)x^{*}_{1}(t)
+H_{1}(t)\Pi(t)(-I^{-1}_{1}(t)( H^{\top}_{1}(t)y_{1}(t)+ \tau\bar{H}_{1}(t)^{\top}\mE[y_{1}(t)]\no\\
&\quad+\tau\ddot{H}_{1}(t)^{\top}\mE^{\sF^{0}_{t}}[y_{1}(t)]+ \tau\tilde{H}_{1}(t+\theta)^{\top}\mE^{\sF_{t}}[y_{1}(t+\theta)]+D_{1}(t)^{\top}z^{0}_{1}(t)+\bar{D}_{1}(t)^{\top}z^{1}_{1}(t)))\no\\
&\quad+ \tau\bar{H}_{1}(t)\mE[\Pi(t)(-I^{-1}_{1}(t)( H^{\top}_{1}(t)y_{1}(t)+ \tau\bar{H}_{1}(t)^{\top}\mE[y_{1}(t)]+\tau\ddot{H}_{1}(t)^{\top}\mE^{\sF^{0}_{t}}[y_{1}(t)]\no\\
&\quad+ \tau\tilde{H}_{1}(t+\theta)^{\top}\mE^{\sF_{t}}[y_{1}(t+\theta)]+D_{1}(t)^{\top}z^{0}_{1}(t)+\bar{D}_{1}(t)^{\top}z^{1}_{1}(t)))]\no\\
&\quad+\tau\ddot{H}_{1}(t)^{\top}\mE^{\sF^{0}_{t}}[\Pi(t)(-I^{-1}_{1}(t)( H^{\top}_{1}(t)y_{1}(t)+ \tau\bar{H}_{1}(t)^{\top}\mE[y_{1}(t)]+\tau\ddot{H}_{1}(t)^{\top}\mE^{\sF^{0}_{t}}[y_{1}(t)]\no\\
&\quad+ \tau\tilde{H}_{1}(t+\theta)^{\top}\mE^{\sF_{t}}[y_{1}(t+\theta)]+D_{1}(t)^{\top}z^{0}_{1}(t)+\bar{D}_{1}(t)^{\top}z^{1}_{1}(t)))  ]\no\\
&\quad+ 1_{\{\theta \leq t \leq T\}}\tau\tilde{H}_{1}(t)\Pi(t-\theta)(-I^{-1}_{1}(t-\theta)( H^{\top}_{1}(t-\theta)y_{1}(t-\theta)\no\\
&\quad+ \tau\bar{H}_{1}(t-\theta)^{\top}\mE[y_{1}(t-\theta)]+\tau\ddot{H}_{1}(t-\theta)^{\top}\mE^{\sF^{0}_{t-\theta}}[y_{1}(t-\theta)]+ \tau\tilde{H}_{1}(t)^{\top}\mE^{\sF_{t-\theta}}[y_{1}(t)]\no\\
&\quad+D_{1}(t-\theta)^{\top}z^{0}_{1}(t-\theta)+\bar{D}_{1}(t-\theta)^{\top}z^{1}_{1}(t-\theta)))\no\\
&\quad+ 1_{\{0 \leq t < \theta\}}\tau\tilde{H}_{1}(t)\Pi(t-\theta)(-I^{-1}_{1}(t-\theta)( \tau \tilde{H}_{1}(t)^{\top}\mE[y_{1}(t)]))\no\\
& \quad +\bar{F}_{2}(t)\mE[ x^{*}_{2}(t)]+\tilde{F}_{2}(t) x^{*}_{2}(t-\theta)+\bar{F}_{1}(t)\mE[ x^{*}_{1}(t)]+\tilde{F}_{1}(t) x^{*}_{1}(t-\theta)\no\\
&\quad+\ddot{F}_{2}(t)\mE^{\sF^{0}_{t}}[ x^{*}_{2}(t)]+\ddot{F}_{1}(t)\mE^{\sF^{0}_{t}}[ x^{*}_{1}(t)]+\upsilon_{1}(t)]\dif t\no\\
&\quad +\sum^{d}_{j=1}[K_{1j}(t)x^{*}_{1}(t) +D_{1j}(t)\Pi(t)(-I^{-1}_{1}(t)( H^{\top}_{1}(t)y_{1}(t)+ \tau\bar{H}_{1}(t)^{\top}\mE[y_{1}(t)]\no\\
&\quad+\tau\ddot{H}_{1}(t)^{\top}\mE^{\sF^{0}_{t}}[y_{1}(t)]+ \tau\tilde{H}_{1}(t+\theta)^{\top}\mE^{\sF_{t}}[y_{1}(t+\theta)]\no\\
&\quad+D_{1}(t)^{\top}z^{0}_{1}(t)+\bar{D}_{1}(t)^{\top}z^{1}_{1}(t))) +\iota_{1j}(t)]\dif B^{0}_{j}(t) \no\\
&\quad +\sum^{d}_{j=1}[\bar{K}_{1j}(t)x^{*}_{1}(t) +\bar{D}_{1j}(t)\Pi(t)(-I^{-1}_{1}(t)( H^{\top}_{1}(t)y_{1}(t)+ \tau\bar{H}_{1}(t)^{\top}\mE[y_{1}(t)]\no\\
&\quad+ \tau\ddot{H}_{1}(t)^{\top}\mE^{\sF^{0}_{t}}[y_{1}(t)]+ \tau\tilde{H}_{1}(t+\theta)^{\top}\mE^{\sF_{t}}[y_{1}(t+\theta)]\no\\
&\quad+D_{1}(t)^{\top}z^{0}_{1}(t)+\bar{D}_{1}(t)^{\top}z^{1}_{1}(t)))+\bar{\iota}_{1j}(t)]\dif B^{1}_{j}(t),
\end{align}

\begin{align}\label{43n}
\dif y_{1}(t)&=-[R_{1}(t) x^*_{1}(t)+\bar{R}_{1}(t)\mE[x^*_{1}(t)]+\ddot{R}_{1}(t)\mE^{\sF_{t}^{0}}[x^*_{1}(t)]+1_{\{0\leq t\leq T-\theta\}}\tilde{R}_{1}(t) x^*_{1}(t)\no\\
& \quad +F_{1}(t)^{\top}y_{1}(t) +K_{1}(t)^{\top}z^{0}_{1}(t)+\bar{K}_{1}(t)^{\top}z^{1}_{1}(t)\no\\
&\quad+\bar{F}_{1}(t)^{\top}\mE[y_{1}(t)]+\bar{F}_{2}(t)^{\top}\mE[y_{2}(t)]+\ddot{F}_{1}(t)^{\top}\mE^{\sF_{t}^{0}}[y_{1}(t)]+\ddot{F}_{2}(t)^{\top}\mE^{\sF_{t}^{0}}[y_{2}(t)] \no\\
&\quad+\tilde{F}_{1}(t+\theta)^{\top}\mE^{\sF_{t}}[y_{1}(t+\theta)]+\tilde{F}_{2}(t+\theta)^{\top}\mE^{\sF_{t}}[y_{2}(t+\theta)] ]\dif t\no\\
&\quad+\sum^{d}_{i=1}z^{0}_{1j}(t)\dif B^{0}_{j}(t)+\sum^{d}_{i=1}z^{1}_{1j}(t)\dif B^{1}_{j}(t),
\end{align}

\begin{align}\label{43o}
\dif x^{*}_{2}(t)&=[F_{2}(t)x^{*}_{2}(t)+H_{2}(t)\Pi(t)(-I^{-1}_{2}(t)( H^{\top}_{2}(t)y_{2}(t)+ \tau\bar{H}_{2}(t)^{\top}\mE[y_{2}(t)]\no\\
&\quad+\tau\ddot{H}_{2}(t)^{\top}\mE^{\sF^{0}_{t}}[y_{2}(t)]+ \tau\tilde{H}_{2}(t+\theta)^{\top}\mE^{\sF_{t}}[y_{2}(t+\theta)]\no\\
&\quad+D_{2}(t)^{\top}z^{0}_{2}(t)+\bar{D}_{2}(t)^{\top}z^{1}_{2}(t)))\no\\
&\quad+ \tau \bar{H}_{2}(t)\mE[ \Pi(t)(-I^{-1}_{2}(t)( H^{\top}_{2}(t)y_{2}(t)+ \tau\bar{H}_{2}(t)^{\top}\mE[y_{2}(t)]\no\\
&\quad+\tau\ddot{H}_{2}(t)^{\top}\mE^{\sF^{0}_{t}}[y_{2}(t)]+ \tau\tilde{H}_{2}(t+\theta)^{\top}\mE^{\sF_{t}}[y_{2}(t+\theta)]\no\\
&\quad+D_{2}(t)^{\top}z^{0}_{2}(t)+\bar{D}_{2}(t)^{\top}z^{1}_{2}(t)))]\no\\
&\quad+ \tau \ddot{H}_{2}(t)\mE^{\sF^{0}_{t}}[ \Pi(t)(-I^{-1}_{2}(t)( H^{\top}_{2}(t)y_{2}(t)+ \tau\bar{H}_{2}(t)^{\top}\mE[y_{2}(t)]\no\\
&\quad+\tau\ddot{H}_{2}(t)^{\top}\mE^{\sF^{0}_{t}}[y_{2}(t)]+ \tau\tilde{H}_{2}(t+\theta)^{\top}\mE^{\sF_{t}}[y_{2}(t+\theta)]\no\\
&\quad+D_{2}(t)^{\top}z^{0}_{2}(t)+\bar{D}_{2}(t)^{\top}z^{1}_{2}(t)))]\no\\
&\quad+ 1_{\{\theta \leq t \leq T\}}\tau\tilde{H}_{2}(t)\Pi(t-\theta)(-I^{-1}_{2}(t-\theta)( H^{\top}_{2}(t-\theta)y_{2}(t-\theta)\no\\
&\quad+ \tau\bar{H}_{2}(t-\theta)^{\top}\mE[y_{2}(t-\theta)]+ \tau\tilde{H}_{2}(t)^{\top}\mE^{\sF_{t-\theta}}[y_{2}(t)]+ \tau\ddot{H}_{2}(t-\theta)^{\top}\mE^{\sF^{0}_{t-\theta}}[y_{2}(t-\theta)]\no\\ &\quad+D_{2}(t-\theta)^{\top}z^{0}_{2}(t-\theta)+\bar{D}_{2}(t-\theta)^{\top}z^{1}_{2}(t-\theta)))\no\\
&\quad+ 1_{\{0 \leq t < \theta\}}\tau\tilde{H}_{2}(t)\Pi(t-\theta)(-I^{-1}_{2}(t-\theta)( \tau \tilde{H}_{2}(t)^{\top}\mE[y_{2}(t)]))\no\\
&\quad +\bar{F}_{2}(t)\mE[ x^{*}_{1}(t) ]+\ddot{F}_{2}(t)\mE^{\sF^{0}_{t}}[ x^{*}_{1}(t) ] +\tilde{F}_{2}(t) x^{*}_{1}(t-\theta) +\upsilon_{2}(t)]\dif t\no\\
&\quad+\sum^{d}_{j=1}[K_{2j}(t)x^{*}_{1}(t) +D_{2j}(t)\Pi(t)(-I^{-1}_{2}(t)(  H^{\top}_{2}(t)y_{2}(t)+ \tau\bar{H}_{2}(t)^{\top}\mE[y_{2}(t)]\no\\
&\quad+\tau\ddot{H}_{2}(t)^{\top}\mE^{\sF^{0}_{t}}[y_{2}(t)]+ \tau\tilde{H}_{2}(t+\theta)^{\top}\mE^{\sF_{t}}[y_{2}(t+\theta)]\no\\
&\quad+D_{2}(t)^{\top}z^{0}_{2}(t)+\bar{D}_{2}(t)^{\top}z^{1}_{2}(t)))+\iota_{2j}(t)]\dif B^{0}_{j}(t) \no\\
&\quad+\sum^{d}_{j=1}[\bar{K}_{2j}(t)x^{*}_{2}(t) +\bar{D}_{2j}(t)\Pi(t)(-I^{-1}_{2}(t)( H^{\top}_{2}(t)y_{2}(t)+ \tau\bar{H}_{2}(t)^{\top}\mE[y_{2}(t)]\no\\
&\quad+\tau\ddot{H}_{2}(t)^{\top}\mE^{\sF^{0}_{t}}[y_{2}(t)]+ \tau\tilde{H}_{2}(t+\theta)^{\top}\mE^{\sF_{t}}[y_{2}(t+\theta)]\no\\
&\quad+D_{2}(t)^{\top}z^{0}_{2}(t)+\bar{D}_{2}(t)^{\top}z^{1}_{2}(t)))+\bar{\iota}_{2j}(t)]\dif B^{1}_{j}(t),
\end{align}

\begin{align}\label{43p}
\dif y_{2}(t)&=-[R_{2}(t) x^*_{2}(t)+\bar{R}_{2}(t)\mE[x^*_{2}(t)]+1_{\{0\leq t\leq T-\theta\}}\tilde{R}_{2}(t) x^*_{2}(t)+\ddot{R}_{2}(t)\mE^{\sF_{t}^{0}}[x^*_{2}(t)]\no\\
& \quad  +F_{2}(t)^{\top}y_{2}(t)+K_{2}(t)^{\top}z^{0}_{2}(t)+\bar{K}_{2}(t)^{\top}z^{1}_{2}(t)\no\\
&\quad+\bar{F}_{2}(t)^{\top}\mE[y_{1}(t)]+\ddot{F}_{2}(t)^{\top}\mE^{\sF_{t}^{0}}[y_{1}(t)]+\tilde{F}_{2}(t+\theta)^{\top}\mE^{\sF_{t}}[y_{1}(t+\theta)] ]\dif t\no\\
&\quad+\sum^{d}_{j=1}z^{0}_{2j}(t)\dif B^{0}_{j}(t)+\sum^{d}_{j=1}z^{1}_{2j}(t)\dif B^{1}_{j}(t), t\in [0,T],
\end{align}
with the initial values and terminal values
\begin{align*}
\begin{cases}
& x^{*}_{1}(0)=x^*_{1}(0) = \sH\bar{\Pi}_{0}(-A^{-1}_{1}\sH^{\top}y_{1}(0)) + x_0,\\
 &x^{*}_{2}(0)=x^*_{2}(0) = \sH\bar{\Pi}_{0}(-A^{-1}_{2}\sH^{\top}y_{2}(0)) + x_0,\\
&x^{*}_{1}(t)=\Pi_{0}(t)(-\tilde{R}_{1}^{-1}(t)( \tilde{F}_{2}(t+\theta)^{\top}\mE[y_{2}(t+\theta)]+\tilde{F}_{1}(t+\theta)^{\top}\mE[y_{1}(t+\theta)])), t\in [-\theta, 0),\\
&x^{*}_{2}(t)=\Pi_{0}(t)(-\tilde{R}_{1}^{-1}(t)(\tilde{F}_{2}(t+\theta)^{\top}\mE[y_{2}(t+\theta)])), t\in [-\theta, 0),\\
& y_{1}(T)=C_{1}(x^*_{1}(T)+x^*_{2}(T))+C_{2}(x^*_{1}(T)+x^*_{2}(T)),\\
 & y_{2}(T)= C_{1}(x^*_{1}(T)+x^*_{2}(T))+C_{2}(x^*_{1}(T)+x^*_{2}(T)),\\
& y_{1}(t)=\eta_{1}(t),y_{2}(t)=\eta_{2}(t), t\in (T, T+\theta].\\
\end{cases}
\end{align*}
We denote   \eqref{43m}-\eqref{43p} by   $(\pi_{3}).$
Here, $\bar{\Pi}_0(\cdot)$ represents the projections from $\mathbb{R}^m$   onto their corresponding closed convex subsets $\bar{\mU}_0$  under the specified norm
$$
|\cdot|_{A_{i}}:= \langle \cdot, \cdot \rangle_{A_{i}}^{1/2} := \langle A_{i} \cdot, \cdot \rangle^{1/2}, i=1,2,
$$
Here, $\Pi_0(t,\cdot), t\in [-\theta, 0]$  represent the projections from $\mathbb{R}^n$   onto their corresponding closed convex subsets $\mU_0(t), t\in [-\theta, 0]$ under the specified norm

$$
|\cdot|_{\tilde{R}_{i}(t)}:= \langle \cdot, \cdot \rangle_{\tilde{R}_{i}(t)}^{1/2} := \langle \tilde{R}_{i}(t) \cdot, \cdot \rangle^{1/2},  t\in [-\theta,0], i=1,2,
$$

 $\Pi(\omega, t, \cdot)$ represent the projections from  $\mathbb{R}^k$ onto their corresponding closed convex subset $\mU(\omega, t)$, under the specified norm
\begin{align*}
|\cdot|_{I_{i}(\omega, t)} := \langle \cdot, \cdot \rangle_{I_{i}(\omega, t)}^{1/2} := \langle I_{i}(\omega, t) \cdot, \cdot \rangle^{1/2},t
\in [-\theta, T], i=1,2,
\end{align*}
for almost every $(\omega, t) \in \Omega \times [-\theta, T]$.
Before giving our main result in this section,  we present the following lemma, which can be found in Lemma 3 of  \cite{LiuNiuWangyu2026}.
  \bl\label{d}
  Under Assumptions 4 and 5,  $\bar{\Pi}_0(\cdot)$  $\Pi_0(t,\cdot)$ and $\Pi(t,\cdot)$ are  Lipschitz continuous uniformly with respect to  $(\omega, t) $ a.e.$-\dif P\times \dif t.$  Moreover, $\bar{\Pi}_{0}, \Pi_0,\Pi$ is $\mathcal{B}(\mathbb{R}^m)$-measurable, $\sB([-\theta, 0])\times \mathcal{B}(\mathbb{R}^n)$-measurable  and  $\mathcal{P} \otimes \mathcal{B}(\mathbb{R}^k)$-measurable, respectively, Moreover,  the following inequalities hold

  \begin{align}\label{e7}
  \delta |\tilde{\nu}|^2 \leq |\tilde{\nu}|_{A_{i}}^2 \leq |A_{i}||\tilde{\nu}|^2 \text{ for any } \tilde{\nu} \in \mathbb{R}^m.
\end{align}

  \begin{align}\label{e8}
  |\bar{\Pi}_{0}(\tilde{\nu})-\bar{\Pi}_{0}(\bar{\tilde{\nu}})|^{2}\leq \frac{|A_{i}|}{\delta}|\tilde{\nu}-\bar{\tilde{\nu}}|^{2}\text{ for any } \nu, \bar{\nu }\in \mathbb{R}^m.
\end{align}

  \begin{align}\label{e}
  \delta |\nu|^2 \leq |\nu|_{\tilde{R}_{i}(t)}^2 \leq \|\tilde{R}_{i}(\cdot)\|_{L^{\infty}(-\theta,0;\mS^{n})} |\nu|^2 \text{ for any } \nu \in \mathbb{R}^n.
\end{align}

  \begin{align}\label{e2}
  |\Pi_{0}(t,\nu)-\Pi_{0}(t,\bar{\nu})|^{2}\leq \frac{\|\tilde{R}_{i}(\cdot)\|_{L^{\infty}(-\theta,0;\mS^{n})}}{\delta}|\nu-\bar{\nu}|^{2}\text{ for any } \nu, \bar{\nu }\in \mathbb{R}^n.
\end{align}

\begin{align}\label{e}
  \delta |u|^2 \leq |u|_{I_{i}(\omega,t)}^2 \leq \|I_{i}(\cdot)\|_{L^{\infty}_{\mF}(-\theta,T;\mS^{k})} |u|^2 \text{ for any } u \in \mathbb{R}^k.
\end{align}
  \begin{align}\label{e2}
  |\Pi(\omega,t,u)-\Pi(\omega,t,\bar{u})|^{2}\leq \frac{\|I_{i}(\cdot)\|_{L^{\infty}_{\mF}(-\theta,T;\mS^{k})}}{\delta}|u-\bar{u}|^{2} \text{ for any } u,\bar{u} \in \mathbb{R}^k.
\end{align}
 \el

Now, we give the main result in this section.

 \bt
 Suppose Assumptions 2, 4, and 5 are satisfied. Then, the Hamiltonian system $(\pi_{3})$ possesses a unique solution $ (W_{1}(\cdot),W_{2}(\cdot))$  where\\  $$W_{1}(\cdot) := (x^*_{1}(\cdot)^{\top}, y_{1}(\cdot)^{\top}, z^{0}_{1}(\cdot)^{\top}, z^{1}_{1}(\cdot)^{\top})^{\top} \in \mM_{\mathbb{F}}(\mathbb{R}^{n+n+nd+nd})$$, $$W_{2}(\cdot) := (x^*_{2}(\cdot)^{\top}, y_{2}(\cdot)^{\top}, z^{0}_{2}(\cdot)^{\top}, z^{1}_{2}(\cdot)^{\top})^{\top} \in \mM_{\mathbb{F}}(\mathbb{R}^{n+n+nd+nd}).$$  Moreover,

\begin{equation}\label{44+-+}
\begin{cases}
\xi^*_{1} =  \bar{\Pi}_{0}(-A_{1}^{-1}\sH^{\top}y_{1}(0)), \\
\xi^*_{2} = \bar{\Pi}_{0}(-A_{2}^{-1}\sH^{\top}y_{2}(0)), \\
u^{*}_{1}(t)=\Pi(t)(-I^{-1}_{1}(t)( H^{\top}_{1}(t)y_{1}(t)+ \tau\bar{H}_{1}^{\top}(t)\mE[y_{1}(t)]+ \tau\ddot{H}_{1}^{\top}(t)\mE^{\sF^{0}_{t}}[y_{1}(t)]\\
\quad\quad\quad\quad\quad\quad\quad\quad+ \tau\tilde{H}_{1}^{\top}(t+\theta)\mE^{\sF_{t}}[y_{1}(t+\theta)]+D^{\top}_{1}(t)z^{0}_{1}(t)+\bar{D}^{\top}_{1}(t)z^{1}_{1}(t))),\\
u^{*}_{2}(t)=\Pi(t)(-I^{-1}_{2}(t)( H^{\top}_{2}(t)y_{2}(t)+ \tau\bar{H}_{2}^{\top}(t)\mE[y_{2}(t)]+ \tau\ddot{H}_{2}^{\top}(t)\mE^{\sF^{0}_{t}}[y_{2}(t)]\\
\quad\quad\quad\quad\quad\quad\quad\quad+ \tau\tilde{H}_{2}^{\top}(t+\theta)\mE^{\sF_{t}}[y_{2}(t+\theta)]+D^{\top}_{2}(t)z^{0}_{2}(t)+\bar{D}^{\top}_{2}(t)z^{1}_{2}(t))), t\in [0,T],
\end{cases}
\end{equation}
and
\begin{align}\label{44++++-}
\begin{cases}
x^{*}_{1}(t)=\Pi_{0}(t)(-\tilde{R}^{-1}_{1}(t)( \tilde{F}_{2}(t+\theta)^{\top}\mE[y_{2}(t+\theta)]+\tilde{F}_{1}(t+\theta)^{\top}\mE[y_{1}(t+\theta)])),\\
x^{*}_{2}(t)=\Pi_{0}(t)(-\tilde{R}^{-1}_{2}(t)( \tilde{F}_{2}(t+\theta)^{\top}\mE[y_{2}(t+\theta)])),\\
u^{*}_{1}(t)=\Pi(t)(-I^{-1}_{1}(t)( \tau \tilde{H}_{1}(t+\theta)^{\top}\mE[y_{1}(t+\theta)])),\\
u^{*}_{2}(t)=\Pi(t)(-I^{-1}_{2}(t)( \tau \tilde{H}_{2}(t+\theta)^{\top}\mE[y_{2}(t+\theta)])), t\in [-\theta,0),
\end{cases}
\end{align}
is the quartet  of optimal controls for Problem (LQ-IC).
\et
\begin{proof}
$\mathbf{Step\, 1}:$ The objective herein is to establish the unique solvability of  $(\pi_{3})$.
To employ Theorem \ref{tt} or the demonstration of this unique solvability  of $(\pi_{3})$,  it is imperative to confirm Assumption 1. We present exclusively the detailed confirmations of $b_{i}(\cdot, 0) \in L_{\mF}^{2}(\mathbb{R}^n), i=1,2$ [refer to Assumption 1(i)] and the monotonicity property of $\aleph_{i}(\cdot),i=1,2$ [refer to Assumption 1(iii)3)]. Further specifics are not shown.

In the beginning, in this specific scenario, $$b_{1}(\cdot, 0) = H_{1}(\cdot)\Pi(\cdot, 0) + \upsilon_{1}(\cdot), b_{2}(\cdot, 0) = H_{2}(\cdot)\Pi(\cdot, 0) + \upsilon_{2}(\cdot).$$ Given that $H_{1}(\cdot), H_{2}(\cdot) \in L^{\infty}(0,T;\mathbb{R}^{n \times k})$ and $\upsilon_{1}(\cdot), \upsilon_{2}(\cdot) \in L_{\mF}^{2}(0,T;\mathbb{R}^n)$, the overall problem can be reduced to checking whether $\Pi(\cdot, 0) \in L_{\mF}^{2}(-\theta, T;\mathbb{R}^k)$.  Indeed, with the assistance of Assumption 5(iii) and the Lipschitz continuity of $\Pi(t, \cdot)$ (as stated in Lemma \ref{d}), hence, we deduce

\begin{align*}
|\Pi(t, 0)| &\leq |\Pi(t, 0) - a(t)| + |a(t)| \\
&= |\Pi(t, 0) - \Pi(t, a(t))| + |a(t)| \\
&\leq \sqrt{\frac{\|I(\cdot)\|_{L^{\infty}_{\mF}(-\theta, T;\mS^{k})}}{\delta}} |a(t)| + |a(t)|
\end{align*}
for almost every $(\omega, t) \in \Omega \times [-\theta, T]$. This leads to the conclusion that $\Pi(\cdot, 0) \in L_{\mF}^2(-\theta, T;\mathbb{R}^k) $.

Secondly, assume that $(\tilde{W}_{1}(\cdot),\tilde{W}_{2}(\cdot))$ and $(\bar{W}_{1}(\cdot),\bar{W}_{2}(\cdot))$ are two solution of  Eq.$(\pi_{3})$  ,where
 $$\tilde{W}_{1}(\cdot) := (\tilde{x}_{1}(\cdot)^{\top}, \tilde{y}_{1}(\cdot)^{\top}, \tilde{z}^{0}_{1}(\cdot)^{\top}, \tilde{z}^{1}_{1}(\cdot)^{\top})^{\top} \in \mM_{\mathbb{F}}(\mathbb{R}^{n+n+nd+nd}),$$  $$\tilde{W}_{2}(\cdot) := (\tilde{x}_{2}(\cdot)^{\top}, \tilde{y}_{2}(\cdot)^{\top}, \tilde{z}^{0}_{2}(\cdot)^{\top},\tilde{z}^{1}_{2}(\cdot)^{\top} )^{\top} \in \mM_{\mathbb{F}}(\mathbb{R}^{n+n+nd+nd})$$ and $$\bar{W}_{1}(\cdot) := (\bar{x}_{1}(\cdot)^{\top}, \bar{y}_{1}(\cdot)^{\top}, \bar{z}^{0}_{1}(\cdot)^{\top}, \bar{z}^{1}_{1}(\cdot)^{\top})^{\top} \in \mM_{\mathbb{F}}(\mathbb{R}^{n+n+nd+nd}),$$ $$\bar{W}_{2}(\cdot) := (\bar{x}_{2}(\cdot)^{\top}, \bar{y}_{2}(\cdot)^{\top}, \bar{z}^{0}_{2}(\cdot)^{\top}, \bar{z}^{1}_{2}(\cdot)^{\top})^{\top} \in \mM_{\mathbb{F}}(\mathbb{R}^{n+n+nd+nd}).$$
 Using the similar computation in \eqref{302}-\eqref{301+}, we can deduce that the associated monotonicity property holds for $\aleph_{1}(\cdot), \aleph_{2}(\cdot)$
with   adjoint functions provided below,

\begin{equation*}
\begin{cases}
    \bar{\chi}_{11}(v) =   \bar{\Pi}_0 (-A_{1}^{-1}v),\\
    \bar{\chi}_{12}(v) =   0,\\
    \bar{\chi}_{21}(v) =  \bar{\Pi}_0 (-A_{2}^{-1}v),\\
    \bar{\chi}_{22}(v) =  0,\\
     \ddot{\chi}_{1}(t, y_{1},y_{2}) = \Pi_{0} (t, -\tilde{R}^{-1}_{1}(t)(\tilde{F}_{1}(t+\theta)^{\top}y_{1}+\tilde{F}_{2}(t+\theta)^{\top}y_{2} )),t\in [-\theta,0),\\
     \ddot{\chi}_{2}(t, y_{1},y_{2}) = \Pi_{0} (t, -\tilde{R}^{-1}_{2}(t)(\tilde{F}_{2}(t+\theta)^{\top}y_{2} )),t\in [-\theta,0),\\
    \chi_{1}(t, u) = \Pi (t, -I^{-1}_{1}(t)u),,t\in [-\theta,T],\\
  \chi_{2}(t, u) = \Pi (t, -I^{-1}_{2}(t)u),t\in [-\theta,T].
\end{cases}
\end{equation*}

$\mathbf{Step\, 2}:$ $(\nu^*_{1},\nu^*_{2}, u^*_{1}, u^*_{2})$ is the only optimal control quartet.
Assume that the corresponding admissible sextet is $(\nu^*_{1},\nu^*_{2}, u^*_{1}, u^*_{2},x^{*}_{1}, x^{*}_{2}).$\\
Assume that  $(\nu_{1},\nu_{2}, u_{1}, u_{2}, x_{1},x_{2})$ is another arbitrary admissible sextet.  We examine the difference

\begin{equation*}
    \mJ(\nu_{1},\nu_{2}, u_{1}, u_{2})) - \mJ(\nu^*_{1},\nu^*_{2}, u^*_{1}, u^*_{2}) = \Re_{81} + \Re_{82}+\Re_{91} + \Re_{92}+\Re
\end{equation*}
where

\begin{align*}
    &\Re_{81}:= \frac{1}{2}\langle A_{1}(\xi_{1} - \xi^{*}_{1}), \xi_{1} - \xi^{*}_{1} \rangle \\
    &+\frac{1}{2} \mathbb{E}\bigg\{\bigg. \langle C_{1}(x_{1}(T)+x_{2}(T) - x^*_{1}(T)- x^*_{2}(T)),x_{1}(T)+x_{2}(T) - x^*_{1}(T))- x^*_{2}(T) \rangle\\ &+\int_{0}^{T}  \langle R_{1}(t)(x_{1}(t) - x^{*}_{1}(t)), x_{1}(t) - x^{*}_{1}(t) \rangle   \dif t\\
    &+ \int_{0}^{T}  \langle \bar{R}_{1}(t)\mE[x_{1}(t) - x^{*}_{1}(t)], x_{1}(t) - x^{*}_{1}(t) \rangle   \dif t\\
    &+\int_{0}^{T}  \langle \tilde{R}_{1}(t-\theta)[x_{1}(t-\theta) - x^{*}_{1}(t-\theta)], x_{1}(t-\theta) - x^{*}_{1}(t-\theta) \rangle   \dif t\no\\
    &+ \int^{T}_{0}\langle \ddot{R}_{1}(t)\mE^{\sF^{0}_{t}}[x_{1}(t)- x^{*}_{1}(t)], \mE^{\sF^{0}_{t}}[x_{1}(t)- x^{*}_{1}(t)]\rangle\dif t\no\\
    &+ \int_{0}^{T}\langle I_{1}(t)(u_{1}(t) - u^{*}_{1}(t)), u_{1}(t) - u^{*}_{1}(t) \rangle\dif t+ \int^{0}_{-\theta}\langle I_{1}(t)(u_{1}(t) - u^{*}_{1}(t)), u_{1}(t) - u^{*}_{1}(t) \rangle\dif t
 \bigg\}\bigg..
\end{align*}

\begin{align*}
\Re_{91} &:= \langle A_{1}\xi^{*}_{1}, \xi_{1}-\xi^{*}_{1} \rangle \\
&+  \mathbb{E}\bigg\{\bigg. \langle C_{1}( x^*_{1}(T)+x^*_{2}(T)),x_{1}(T)+x_{2}(T) - x^*_{1}(T))- x^*_{2}(T) \rangle\\
&+ \int_0^T \langle R_{1}(t)x^*_{1}(t), x_{1}(t) - x^*_{1}(t) \rangle\dif t+\int_0^T \langle \bar{R}_{1}(t)\mE[x^*_{1}(t)], x_{1}(t) - x^*_{1}(t) \rangle\dif t\no\\
&+\int_0^T \langle \tilde{R}_{1}(t-\theta)x^*_{1}(t-\theta), x_{1}(t-\theta) - x^*_{1}(t-\theta)\rangle\dif t\\
  &+ \int^{T}_{0}\langle \ddot{R}_{1}(t)\mE^{\sF^{0}_{t}}[x^{*}_{1}(t)], \mE^{\sF^{0}_{t}}[x_{1}(t)- x^{*}_{1}(t)]\rangle\dif t\no\\
  &+ \int_0^T\langle I_{1}(t)u^*_{1}(t), u_{1}(t) - u^*_{1}(t) \rangle ] \dif t
+\int_{-\theta}^{0}\langle I_{1}(t)u^*_{1}(t), u_{1}(t) - u^*_{1}(t) \rangle  \dif t\bigg\}\bigg..
\end{align*}

\begin{align*}
    \Re_{82} &:= \frac{1}{2}\langle A_{2}(\xi_{2}- \xi^{*}_{2}), \xi_{2}- \xi^{*}_{2} \rangle\\
    &+ \frac{1}{2} \mathbb{E}\bigg\{\bigg. \langle C_{2}(x_{1}(T)+x_{2}(T) - x^*_{1}(T))- x^*_{2}(T)),x_{1}(T)+x_{2}(T) - x^*_{1}(T))- x^*_{2}(T) \rangle\\  &+\int_{0}^{T}  \langle R_{2}(t)(x_{1}(t) - x^{*}_{1}(t)), x_{1}(t) - x^{*}_{1}(t) \rangle   \dif t\\
    &+ \int_{0}^{T}  \langle \bar{R}_{2}(t)\mE[x_{1}(t) - x^{*}_{1}(t)], x_{1}(t) - x^{*}_{1}(t) \rangle   \dif t\\
    &+ \int^{T}_{0}\langle \ddot{R}_{2}(t)\mE^{\sF^{0}_{t}}[x_{2}(t)- x^{*}_{2}(t)], \mE^{\sF^{0}_{t}}[x_{2}(t)- x^{*}_{2}(t)]\rangle\dif t\no\\
    &+\int_{0}^{T}  \langle \tilde{R}_{2}(t-\theta)[x_{1}(t-\theta) - x^{*}_{1}(t-\theta)], x_{1}(t-\theta) - x^{*}_{1}(t-\theta) \rangle   \dif t\no\\
    &+ \int_{0}^{T}\langle I_{2}(t)(u_{2}(t) - u^{*}_{2}(t)), u_{2}(t) - u^{*}_{2}(t) \rangle\dif t+ \int^{0}_{-\theta}\langle I_{2}(t)(u_{2}(t) - u^{*}_{2}(t)), u_{2}(t) - u^{*}_{2}(t) \rangle\dif t
    \bigg\}\bigg..
\end{align*}

\begin{align*}
\Re_{92} &:= \langle A_{2}\xi^{*}_{2}, \xi_{2} - \xi^{*}_{2} \rangle  \\
& + \mathbb{E}\bigg\{\bigg. \langle C_{2}( x^*_{1}(T)+ x^*_{2}(T)),x_{1}(T)+x_{2}(T) - x^*_{1}(T))- x^*_{2}(T) \rangle\\
&+ \int_0^T \langle R_{2}(t)x^*_{2}(t), x_{2}(t) - x^*_{2}(t) \rangle\dif t+\int_0^T \langle \bar{R}_{2}(t)\mE[x^*_{2}(t)], x_{2} - x^*_{2} \rangle\dif t\no\\
&+\int_0^T \langle \tilde{R}_{2}(t-\theta)x^*_{2}(t-\theta), x_{2}(t-\theta) - x^*_{2}(t-\theta)\rangle\dif t\\
   &+ \int^{T}_{0}\langle \ddot{R}_{2}(t)\mE^{\sF^{0}_{t}}[x^{*}_{2}(t)], \mE^{\sF^{0}_{t}}[x_{2}(t)- x^{*}_{2}(t)]\rangle\dif t\no\\
  &+ \int_0^T\langle I_{2}(t)u^*_{2}(t), u_{2}(t) - u^*_{2}(t) \rangle ] \dif t
+\int_{-\theta}^{0}\langle I_{2}(t)u^*_{2}(t), u_{2}(t) - u^*_{2}(t) \rangle  \dif t\bigg\}\bigg..
\end{align*}

\begin{align*}
\Re&:=\int^{T}_{T-\theta}\langle \tau \tilde{H}_{1}(t+\theta)^{\top}\mE^{\sF_{t}}[\eta_{1}(t+\theta)],  u_{1}(t)-u^{*}_{1}(t)\rangle\dif t\\
&+\int^{T}_{T-\theta}\langle \tau \tilde{H}_{2}(t+\theta)^{\top}\mE^{\sF_{t}}[\eta_{2}(t+\theta)],  u_{2}(t)-u^{*}_{2}(t)\rangle\dif t\\
&+\int^{T}_{T-\theta}\langle \tilde{F}_{2}(t+\theta)^{\top}\mE^{\sF_{t}}[\eta_{2}(t+\theta)]+\tilde{F}_{1}(t+\theta)^{\top}\mE^{\sF_{t}}[\eta_{1}(t+\theta)], x_{1}(t)-x^{*}_{1}(t) \rangle\dif t\\
&+\int^{T}_{T-\theta}\langle \tilde{F}_{2}(t+\theta)^{\top}\mE^{\sF_{t}}[\eta_{1}(t+\theta)], x_{2}(t)-x^{*}_{2}(t) \rangle\dif t.
\end{align*}
Initially, applying It\^{o}'s formula to $\langle y_{1}(\cdot), x_{1}(\cdot) - x^*_{1}(\cdot) \rangle$ and $\langle y_{2}(\cdot), x_{2}(\cdot) - x^*_{2}(\cdot) \rangle$  and using a calculation method similar to that used to obtain \eqref{34}, we derive

\begin{align*}
&\mathbb{E}\bigg\{\bigg. \langle C_{1}(x^*_{1}(T)+x^*_{2}(T)), x_{1}(T) - x^*_{1}(T) \rangle + \langle C_{2}(x^*_{1}(T)+x^*_{2}(T)), x_{1}(T) - x^*_{1}(T) \rangle \\
&+\int_0^T \langle R_{1}(t)x^*_{1}(t), x_{1}(t) - x^*_{1}(t) \rangle \dif t + \int_0^T \langle \bar{R}_{1}(t)\mE[x^*_{1}(t)], \mE[x_{1}(t) - x^*_{1}(t)] \rangle \dif t \\
&+ \int^{T}_{0}\langle \ddot{R}_{1}(t)\mE^{\sF^{0}_{t}}[x^{*}_{1}(t)], \mE^{\sF^{0}_{t}}[x_{1}(t)-x^{*}_{1}(t)]\rangle\dif t\\
 &+ \int_0^T \langle \tilde{R}_{1}(t-\theta)x^*_{1}(t-\theta), x_{1}(t-\theta) - x^*_{1}(t-\theta) \rangle \dif t\bigg\}\bigg.\\
&+\mathbb{E}\bigg\{\bigg. \langle C_{1}(x^*_{1}(T)+x^*_{2}(T)), x_{2}(T)) - x^*_{2}(T) \rangle+\langle C_{2}(x^*_{1}(T)+x^*_{2}(T)), x_{2}(T)) - x^*_{2}(T) \rangle\dif t\\
&+  \int_0^T \langle R_{2}(t)x^*_{2}(t), x_{2}(t) - x^*_{2}(t) \rangle+ \int_0^T \langle \bar{R}_{2}(t)\mE[x^*_{2}(t)], x_{2}(t) - x^*_{2}(t) \rangle \dif t \\
&+ \int^{T}_{0}\langle \ddot{R}_{2}(t)\mE^{\sF^{0}_{t}}[x^{*}_{2}(t)], \mE^{\sF^{0}_{t}}[x_{2}(t)-x^{*}_{2}(t)]\rangle\dif t\\
&+ \int_0^T \langle \tilde{R}_{2}(t-\theta)x^*_{2}(t-\theta), x_{2}(t-\theta) - x^*_{2}(t-\theta) \rangle \dif t \no\\
  &+\int^{T}_{T-\theta}\langle \tilde{F}_{2}(t+\theta)^{\top}\mE^{\sF_{t}}[\eta_{2}(t+\theta)], x_{1}(t)-x^{*}_{1}(t) \rangle\dif t\no\\
 &+\int^{T}_{T-\theta}\langle \tilde{F}_{1}(t+\theta)^{\top}\mE^{\sF_{t}}[\eta_{1}(t+\theta)], x_{1}(t)-x^{*}_{1}(t) \rangle\dif t\no\\
  &+\int^{T}_{T-\theta}\langle \tilde{F}_{2}(t+\theta)^{\top}\mE^{\sF_{t}}[\eta_{1}(t+\theta)], x_{2}(t)-x^{*}_{2}(t) \rangle\dif t\no\\
  &  + \int^{T}_{T-\theta}\langle \tau \tilde{H}_{1}(t+\theta)^{\top}\mE^{\sF_{t}}[\eta_{1}(t+\theta)],  u_{1}(t)-u^{*}_{1}(t)\rangle\dif t\no\\
&+ \int^{T}_{T-\theta}\langle \tau \tilde{H}_{2}(t+\theta)^{\top}\mE^{\sF_{t}}[\eta_{2}(t+\theta)],  u_{2}(t)-u^{*}_{2}(t)\rangle\dif t
\bigg\}\bigg.\\
&= \langle \sH^T y_{1}(0), \xi_{1} - \xi^*_{1} \rangle +\langle \sH^T y_{2}(0), \xi_{2} - \xi^*_{2} \rangle\\
 & +\mE\int^{T}_{0}\langle   H_{1}(t)^{\top}y_{1}(t)+ \tau\bar{H}_{1}(t)^{\top}\mE[y_{1}(t)]+ \tau\ddot{H}_{1}^{\top}(t)\mE^{\sF^{0}_{t}}[y_{1}(t)]\no\\
 &\quad\quad\quad\quad\quad+ \tau\tilde{H}_{1}(t+\theta)^{\top}\mE^{\sF_{t}}[y_{1}(t+\theta)]+D_{1}(t)^{\top}z^{0}_{1}(t)+\bar{D}_{1}(t)^{\top}z^{1}_{1}(t), u_{1}(t)-u^{*}_{1}(t)   \rangle\dif t\no\\
& +\mE\int^{T}_{0}\langle H_{2}(t)^{\top}y_{2}(t)+ \tau\bar{H}_{2}(t)^{\top}\mE[y_{2}(t)]+ \tau\ddot{H}_{2}^{\top}(t)\mE^{\sF^{0}_{t}}[y_{2}(t)]\no\\
&\quad\quad\quad\quad\quad+ \tau\tilde{H}_{2}(t+\theta)^{\top}\mE^{\sF_{t}}[y_{2}(t+\theta)]+D_{2}(t)^{\top}z^{0}_{2}(t)+\bar{D}_{2}(t)^{\top}z^{1}_{2}(t),   u_{2}(t)-u^{*}_{2}(t) \rangle\dif t\no\\
&+\int^{0}_{-\theta}\langle \tilde{R}_{1}(t)x^{*}_{1}(t)+\tilde{F}_{2}(t+\theta)^{\top}\mE[y_{2}(t+\theta)]+\tilde{F}_{1}(t+\theta)^{\top}\mE^{\sF_{t}}[y_{1}(t+\theta)],   x_{1}(t)-x^{*}_{1}(t)\rangle \dif t\no\\
&+\int^{0}_{-\theta}\langle \tilde{R}_{2}(t)x^{*}_{2}(t)+\tilde{F}_{2}(t+\theta)^{\top}\mE[y_{1}(t+\theta)],   x_{2}(t)-x^{*}_{2}(t)\rangle \dif t\no\\
&+\int^{0}_{-\theta}\langle\tau \tilde{H}_{1}(t+\theta)^{\top}\mE[y_{1}(t+\theta)],   u_{1}(t)-u^{*}_{1}(t)\rangle \dif t\no\\
&+\int^{0}_{-\theta}\langle\tau \tilde{H}_{2}(t+\theta)^{\top}\mE[y_{2}(t+\theta)],   u_{2}(t)-u^{*}_{2}(t)\rangle \dif t.
\end{align*}

Substituting the above calculations into $\Re_{91}+\Re_{92}+\Re$ gives
\begin{align*}
&\Re_{91}+\Re_{92}+\Re = \langle A_{1}\xi^{*}_{1}+H^{\top} y_{1}(0), \xi_{1}-\xi^{*}_{1} \rangle+\langle A_{2}\xi^{*}_{2}+H^{\top} y_{2}(0), \xi_{2} - \xi^{*}_{2} \rangle \\
&+\int^{0}_{-\theta}\langle \tilde{R}_{1}(t)x^{*}_{1}(t)+\tilde{F}_{2}(t+\theta)^{\top}\mE[y_{2}(t+\theta)]+\tilde{F}_{1}(t+\theta)^{\top}\mE^{\sF_{t}}[y_{1}(t+\theta)],   x_{1}(t)-x^{*}_{1}(t)\rangle \dif t\no\\
&+\int^{0}_{-\theta}\langle \tilde{R}_{2}(t)x^{*}_{2}(t)+\tilde{F}_{2}(t+\theta)^{\top}\mE[y_{1}(t+\theta)],   x_{2}(t)-x^{*}_{2}(t)\rangle \dif t\no\\
& + \mathbb{E} \int_0^T \langle I_{1}(t)u^*_{1}(t) + H_{1}(t)^{\top}y_{1}+ \tau\bar{H}_{1}(t)^{\top}\mE[y_{1}(t)]+ \tau\ddot{H}_{1}^{\top}(t)\mE^{\sF^{0}_{t}}[y_{1}(t)]\\
&\quad\quad\quad+ \tau\tilde{H}_{1}(t+\theta)^{\top}\mE^{\sF_{t}}[y_{1}(t+\theta)]+D_{1}(t)^{\top}z^{0}_{1}(t)+\bar{D}_{1}(t)^{\top}z^{1}_{1}(t), u_{1}(t) - u^*_{1}(t) \rangle dt\no\\
&+ \mathbb{E} \int_0^T \langle I_{2}(t)u^*_{2}(t)+ H_{2}(t)^{\top}y_{2} + \tau\bar{H}_{2}(t)^{\top}\mE[y_{2}(t)]+ \tau\ddot{H}_{1}^{\top}(t)\mE^{\sF^{0}_{t}}[y_{1}(t)]\\
&\quad\quad\quad+ \tau\tilde{H}_{2}(t+\theta)^{\top}\mE^{\sF_{t}}[y_{2}(t+\theta)]+D_{2}(t)^{\top}z^{0}_{2}(t)+\bar{D}_{2}(t)^{\top}z^{1}_{2}(t), u_{2}(t) - u^*_{2}(t) \rangle \dif t\\
&+\int^{0}_{-\theta}\langle I_{1}(t)u^{*}_{1}(t)+\tau \tilde{H}_{1}(t+\theta)^{\top}\mE[y_{1}(t+\theta)],   u_{1}(t)-u^{*}_{1}(t)\rangle \dif t\no\\
&+\int^{0}_{-\theta}\langle I_{2}(t)u^{*}_{2}(t)+\tau \tilde{H}_{2}(t+\theta)^{\top}\mE[y_{2}(t+\theta)],   u_{2}(t)-u^{*}_{2}(t)\rangle \dif t
\end{align*}
Because of \eqref{a} in Lemma \ref{c} and the definition of $\xi^*$ [refer to \eqref{44+-+} and \eqref{44++++-}], the sum of the first two terms of $\Re_{91}+\Re_{92}+\Re$  simplifies to
\begin{align*}
&\langle A_{1}\xi^{*}_{1}+H^{\top} y_{1}(0), \xi_{1}-\xi^{*}_{1} \rangle+\langle A_{2}\xi^{*}_{2}+H^{\top} y_{2}(0), \xi_{2} - \xi^{*}_{2} \rangle\\
&= - \langle -A_{1}^{-1} H^{\top} y_{1}(0) - \xi^*_{1}, \xi_{1} - \xi^*_{1} \rangle_{A_{1} }- \langle -A_{2}^{-1} H^{\top} y_{2}(0) - \xi^*_{2}, \xi_{2} - \xi^*_{2} \rangle_{A_{2} }\geq 0.
\end{align*}
 Similarly, other terms of $\Re_{91}+\Re_{92}+\Re $ are also found to be non-negative.
 Therefore, we deduce that $\Re_{91}+\Re_{92} +\Re \geq 0$.

On the other hand, by $\mathrm{Assumption\, 4\,(2})$ and $(\nu_{1},\nu_{2}, u_{1}, u_{2})\neq (\nu^*_{1},\nu^*_{2}, u^*_{1}, u^*_{2}),$ it yields

\begin{align*}
&\mathcal{J}(\nu_{1},\nu_{2}, u_{1}, u_{2}) - \mathcal{J}(\nu^*_{1},\nu^*_{2}, u^*_{1}, u^*_{2}) \geq \Re_{81}+\Re_{82}\\
&\geq \frac{\delta}{2} \left\{ |\xi_{1}- \xi^{*}_{1}|^{2} + \mathbb{E} \int_{0}^{T} |u_{1}(t) - u^{*}_{1}(t)|^{2} dt+ \mathbb{E} \int_{-\vartheta}^{0} |u_{1}(t) - u^{*}_{1}(t)|^{2} dt \right\} \\
&+\frac{\delta}{2} \left\{ |\xi_{2} - \xi^{*}_{2}|^{2} + \mathbb{E} \int_{0}^{T} |u_{2}(t) - u^{*}_{2}(t)|^{2} dt+ \mathbb{E} \int_{-\vartheta}^{0} |u_{1}(t) - u^{*}_{1}(t)|^{2} dt \right\} > 0
\end{align*}
 Consequently, $(\nu^*_{1},\nu^*_{1}, u^*_{1}, u^*_{2})$ defined by \eqref{44+-+} and \eqref{44++++-}  constitutes the unique optimal control  quadruple for Problem (LQ-IC). We complete the proof.

\end{proof}

\section{Appendix}

\subsection{Convexity and Uniform Convexity}
We list some auxiliary results in this section.
The following results can be found in \cite{Peypouquet2015,LiuNiuWangyu2026}.
\bl\label{A1}
 Under Assumption 3, both $\nabla g_{1i}(\cdot): \mathbb{R}^{m} \rightarrow \mathbb{R}^{m}, i=1,2$ and $\nabla g_{4i}(\cdot): \mathbb{R}^{k} \rightarrow \mathbb{R}^{k},i=1,2$ are bijective. Let $(\nabla g_{1i})^{-1}(\cdot), i=1,2$ and $(\nabla g_{4i})^{-1}(t, \cdot), i=1,2$ denote the four inverse mappings, respectively. Then they are Lipschitz continuous with Lipschitz constant $1 / \delta>0$. Moreover, $(\nabla g_{4i})^{-1}$ is $ \sP\times\sB(\mathbb{R}^{k})$-measurable.
\el

\bd \label{d1}
Suppose $D \subset \mathbb{R}^n$ is a nonempty and convex set. A function $g : D \rightarrow \mathbb{R}$ is convex if
$$
g(\lambda x + (1 - \lambda)y) \leq \lambda g(x) + (1 - \lambda) g(y)
$$
for any $\lambda \in (0, 1)$ and any $x, y \in D$. If the inequality holds strictly whenever $x \neq y$, then $f(\cdot)$ is called strictly convex. Furthermore, $g(\cdot)$ is uniformly convex (also known as strongly convex) with parameter $\delta > 0$ if
$$
g(\lambda x + (1 - \lambda)y) + \frac{\delta}{2} \lambda (1 - \lambda) |x - y|^2 \leq \lambda g(x) + (1 - \lambda) f(y)
$$
for any $\lambda \in (0, 1)$ and any $x, y \in D$.
\ed

 \bl \label{100}
 Let $D \subset \mathbb{R}^n$ be a nonempty, open, and convex set.  Suppose $g : D \rightarrow \mathbb{R}$ is differentiable.  Then, the following statements are equivalent.

\begin{enumerate}
    \item [(1)]$g(\cdot)$ is convex (respectively, uniformly convex with $\delta > 0$).
    \item[(2)] $g(x) - g(y) - \langle \nabla g(y), x - y \rangle \geq 0$ (respectively, $\geq (\delta/2) |x - y|^2$) for any $x, y \in D$.
    \item[(3)] $\langle \nabla g(x) - \nabla g(y), x - y \rangle \geq 0$ (respectively, $\geq \delta |x - y|^2$) for any $x, y \in D$.
\end{enumerate}

\el

\subsection{Projection Onto a Closed Convex Set}
In this section, we list some basic properties of a projection onto a closed convex subset of $\mathbb{R}^n$. More details can be found in \cite{Brezis2011}. With a slight abuse of notation, the inner product and the induced norm of $\mathbb{R}^n$ (not necessarily the Euclidean inner product and the Euclidean norm) are denoted by $\langle \cdot, \cdot \rangle$ and $|\cdot|$, respectively.

\bl\label{c}
 Let $K \subset \mathbb{R}^n$ be a nonempty closed convex set. Then, for each $x \in \mathbb{R}^n$, there exists a unique element $\Pi(x) \in K$ such that
\[
|x - \Pi(x)| = \min_{y \in K} |x - y|.
\]
\el
\noindent
The element $\Pi(x)$ is called the projection of $x$ onto $K$. Moreover, $\Pi(x) \in K$ is characterized by the property that
\begin{equation}\label{a}
\langle x - \Pi(x), y - \Pi(x) \rangle \leq 0
\end{equation}
for any $y \in K$. Furthermore,
\begin{equation}\label{b}
\begin{cases}
|\Pi(x) - \Pi(\bar{x})|^2 \leq \langle \Pi(x) - \Pi(\bar{x}), x - \bar{x} \rangle \\
|\Pi(x) - \Pi(\bar{x})| \leq |x - \bar{x}|
\end{cases}
\end{equation}
for any $x, \bar{x} \in \mathbb{R}^n$.

\section*{Declaration of competing interest}
The authors declare that they have no known competing financial interests or personal relationships that could have appeared
to influence the work reported in this paper.

\section*{Funding}
This research is supported by the National Natural Science Foundation of China (Grant no.  11626236), the Fundamental Research Funds for the Central Universities of South-Central Minzu University (Grant no. CZY15017).

\end{document}